\documentclass[11pt,a4paper,usenames,dvipsnames]{article}

\usepackage{url,amsmath,amssymb,latexsym,mathrsfs,comment,amsthm,enumerate,geometry,calc,ifthen,mathdots,textcomp,multicol, adjustbox}

\usepackage{cases}

\usepackage[noadjust]{cite}

\usepackage[colorlinks]{hyperref}

\usepackage[margin=10pt,font=small,labelfont=bf, labelsep=period]{caption}

\usepackage{makecell}

\usepackage[all,cmtip,2cell]{xy}

\usepackage[x11names, svgnames, rgb,table]{xcolor}
\usepackage{hhline}
\usepackage{multirow}

\usepackage{pst-node}
\usepackage{tikz-cd} 
\tikzcdset{scale cd/.style={every label/.append style={scale=#1},
    cells={nodes={scale=#1}}}}

\usepackage{enumitem}

\usepackage[utf8]{inputenc}
\usepackage[T1]{fontenc}

\usepackage{tikz}
\pgfdeclarelayer{background layer}
\pgfsetlayers{background layer,main}

\usetikzlibrary{decorations.markings}
  
\tikzset{->-/.style={decoration={
  markings,
  mark=at position #1 with {\arrow{{latex}}}},postaction={decorate}}}

\newcommand\nc\newcommand
\nc\rnc\renewcommand

\nc{\custpartn}[3]{{\lower1.4 ex\hbox{
\begin{tikzpicture}[scale=.3]
\foreach \x in {#1}
{ \uvert{\x}  }
\foreach \x in {#2}
{ \lvert{\x}  }
#3 \end{tikzpicture}
}}}
\nc{\uvert}[1]{\fill (#1,2)circle(.2);}
\rnc{\lvert}[1]{\fill (#1,0)circle(.2);}

\nc\cplab[2]{\node[above] () at (#1,2) {\tiny $#2$};}

\nc\too\rightsquigarrow
\nc\comm{\textup{c}}

\nc\oa\circledast
\nc\wh\widehat
\nc\di\diamondsuit

\nc\N{\mathbb N}
\nc\R{\operatorname{\mathscr R}}
\rnc\L{\operatorname{\mathscr L}}
\nc\J{\operatorname{\mathscr J}}
\nc\D{\operatorname{\mathscr D}}
\rnc\H{\operatorname{\mathscr H}}
\nc\K{\operatorname{\mathscr K}}
\nc\leqL{\leq_{\L}}
\nc\leqR{\leq_{\R}}
\nc\leqH{\leq_{\H}}
\nc\leqJ{\leq_{\J}}
\nc\leqK{\leq_{\K}}
\nc\Reg{\operatorname{Reg}}
\nc\Aut{\operatorname{Aut}}
\nc\PAut{\operatorname{PAut}}
\nc\End{\operatorname{End}}
\nc\PEnd{\operatorname{PEnd}}
\nc\Eq{\mathfrak{Eq}}

\nc\trans[1]{\left(\begin{smallmatrix}#1\end{smallmatrix}\right)}
\nc\sm\setminus

\nc{\mtlab}[1]{\mapstochar \xrightarrow {\ #1\ }}

\nc{\ssim}{\mathrel{\raise0.25 ex\hbox{\oalign{$\approx$\crcr\noalign{\kern-0.84 ex}$\approx$}}}}

\nc{\ldb}{[\hspace{-0.55truemm}[}
\nc{\rdb}{]\hspace{-0.55truemm}]}

\nc\larr{\mathrel{\begin{tikzpicture}\draw[>-](0,0)--(.5,0);\end{tikzpicture}}}
\nc\rarr{\mathrel{\begin{tikzpicture}\draw[->](0,0)--(.5,0);\end{tikzpicture}}}
\nc\llarr{\mathrel{\begin{tikzpicture}\draw[>-<](0,0)--(.5,0);\end{tikzpicture}}}
\nc\rrarr{\mathrel{\begin{tikzpicture}\draw[<->](0,0)--(.5,0);\end{tikzpicture}}}
\nc\lrarr{\mathrel{\begin{tikzpicture}\draw[>->](0,0)--(.5,0);\end{tikzpicture}}}
\nc\arrl{\mathrel{\begin{tikzpicture}\draw[-<](0,0)--(.5,0);\end{tikzpicture}}}
\nc\arrr{\mathrel{\begin{tikzpicture}\draw[<-](0,0)--(.5,0);\end{tikzpicture}}}

\nc\Z{\mathbb Z}

\nc\dda{{\scalebox{0.5}{\begin{tikzpicture}\draw[thick,>->](0,.5)--(0,0);\end{tikzpicture}}}}

\nc\pre\preceq

\newcommand{\uvx}[2]{\fill (#1,2)circle(#2);}
\newcommand{\lvx}[2]{\fill (#1,0)circle(#2);}
\newcommand{\uv}[1]{\fill (#1,2)circle(.17);}
\newcommand{\lv}[1]{\fill (#1,0)circle(.17);}

\newcommand{\uvs}[1]{{\foreach \x in {#1} { \uv{\x}}}}
\newcommand{\lvs}[1]{{\foreach \x in {#1} { \lv{\x}}}}
\newcommand{\darcx}[3]{\draw(#1,0)arc(180:90:#3) (#1+#3,#3)--(#2-#3,#3) (#2-#3,#3) arc(90:0:#3);}
\newcommand{\darc}[2]{\darcx{#1}{#2}{.4}}
\newcommand{\uarcx}[3]{\draw(#1,2)arc(180:270:#3) (#1+#3,2-#3)--(#2-#3,2-#3) (#2-#3,2-#3) arc(270:360:#3);}
\newcommand{\uarc}[2]{\uarcx{#1}{#2}{.4}}
\newcommand{\darcxx}[4]{\draw[#4](#1,0)arc(180:90:#3) (#1+#3,#3)--(#2-#3,#3) (#2-#3,#3) arc(90:0:#3);}
\newcommand{\uarcxx}[4]{\draw[#4](#1,2)arc(180:270:#3) (#1+#3,2-#3)--(#2-#3,2-#3) (#2-#3,2-#3) arc(270:360:#3);}
\newcommand{\stline}[2]{\draw(#1,2)--(#2,0);}
\newcommand{\uline}[2]{\draw(#1,2)--(#2,2);}

\nc{\uarcs}[1]{{\foreach \x/\y in {#1}{ \uarc{\x}{\y} }}}
\nc{\darcs}[1]{{\foreach \x/\y in {#1}{ \darc{\x}{\y} }}}
\newcommand{\stlines}[1]{{\foreach \x/\y in {#1} { \stline{\x}{\y} }}}
\nc\udotted[2]{\draw[dotted](#1+.5,2)--(#2-.5,2);}
\nc\ldotted[2]{\draw[dotted](#1+.5,0)--(#2-.5,0);}
\nc\ddotted[2]{\draw[dotted](#1+.5,0)--(#2-.5,0);}
\nc\vertlab[2]{\node () at (#1,2.7) {\tiny $#2$};}

\newcommand{\stlinex}[3]{\draw[#3](#1,2)--(#2,0);}
\newcommand{\ulinex}[3]{\draw[#3](#1,2)--(#2,2);}
\newcommand{\dlinex}[3]{\draw[#3](#1,0)--(#2,0);}

\nc\Ptwo[3]{{\lower1.0 ex\hbox{
\begin{tikzpicture}[scale=.23]
\uvx1{.25}\uvx2{.25}\lvx1{.25}\lvx2{.25}
\foreach \x/\y in {#1} {\stlinex\x\y{thick}}
\foreach \x/\y in {#2} {\ulinex\x\y{thick}}
\foreach \x/\y in {#3} {\dlinex\x\y{thick}}
\end{tikzpicture}
}}}

\nc\RedPtwo[3]{{\lower1.0 ex\hbox{
\begin{tikzpicture}[red,scale=.23]
\uvx1{.25}\uvx2{.25}\lvx1{.25}\lvx2{.25}
\foreach \x/\y in {#1} {\stlinex\x\y{thick}}
\foreach \x/\y in {#2} {\ulinex\x\y{thick}}
\foreach \x/\y in {#3} {\dlinex\x\y{thick}}
\end{tikzpicture}
}}}

\nc\BluePtwo[3]{{\lower1.0 ex\hbox{
\begin{tikzpicture}[blue,scale=.23]
\uvx1{.25}\uvx2{.25}\lvx1{.25}\lvx2{.25}
\foreach \x/\y in {#1} {\stlinex\x\y{thick}}
\foreach \x/\y in {#2} {\ulinex\x\y{thick}}
\foreach \x/\y in {#3} {\dlinex\x\y{thick}}
\end{tikzpicture}
}}}

\nc\mot[2]{{\lower1.0 ex\hbox{
\begin{tikzpicture}[scale=.23]
\uv1\uv2\uv3\lv1\lv2\lv3
\foreach \x in {#1} {\stline\x\x}
\foreach \x/\y in {#2} {\uarc\x\y \darc\x\y}
\end{tikzpicture}
}}}

\nc\mott[3]{{\lower1.0 ex\hbox{
\begin{tikzpicture}[scale=.23]
\uv1\uv2\uv3\lv1\lv2\lv3
\foreach \x/\y in {#1} {\stline\x\y}
\foreach \x/\y in {#2} {\uarc\x\y}
\foreach \x/\y in {#3} {\darc\x\y}
\end{tikzpicture}
}}}

\nc\moot[2]{{\lower1.4 ex\hbox{
\begin{tikzpicture}[scale=.3]
\uv1\uv2\uv3\uv4\lv1\lv2\lv3\lv4
\foreach \x in {#1} {\stline\x\x}
\foreach \x/\y in {#2} {\uarc\x\y \darc\x\y}
\end{tikzpicture}
}}}

\newcounter{ncols}
\newcounter{incols}
\newenvironment{partn}[1]{
  \setcounter{ncols}{#1} \setcounter{incols}{\thencols - 1}\setlength{\arraycolsep}{1pt}
  \Bigl( \hspace{-1.5truemm}\scriptsize 
    \begin{array}{@{\hskip 3pt} c *{\theincols}{|c} @{\hskip 3pt}  }
}{
     \end{array}
     \normalsize \hspace{-1.5truemm}\Bigr)\setlength{\arraycolsep}{6pt}
}

\nc\la\langle
\nc\ra\rangle

\nc\ol\overline
\nc\ul\underline
\nc\da\downarrow
\nc\pr\bullet
\nc\ccirc\circledcirc
\nc\bc\odot

\nc\RSS{{\bf RSS}}
\nc\BS{{\bf BS}}
\nc\RBS{{\bf RBS}}
\nc\SBS{{\bf SBS}}
\nc\RSBS{{\bf RSBS}}
\nc\PA{{\bf PA}}
\nc\OG{{\bf OG}}
\nc\CPG{{\bf CPG}}
\nc\Set{{\bf Set}}
\nc\Sgp{{\bf Sgp}}
\nc\IS{{\bf IS}}
\nc\IG{\operatorname{\textup{\textsf{IG}}}}
\nc\RIG{\operatorname{\textup{\textsf{RIG}}}}
\nc\bG{{\bf G}}
\nc\bS{{\bf S}}
\nc\bE{{\bf E}}
\nc\bP{{\bf P}}
\nc\bEE{{\bf E}}
\nc\bPP{{\bf P}}
\nc\bB{{\bf B}}
\nc\bC{{\bf C}}
\nc\bD{{\bf D}}
\nc\bbG{\mathbb G}
\nc\bbS{\mathbb S}

\nc\PG{\operatorname{\textup{\textsf{PG}}}}
\nc\PGPz{\PG(P_0)}

\nc\ds\displaystyle
\nc\ts\textstyle

\nc\bn{{\bf n}}

\nc\al\alpha
\nc\be\beta
\nc\ga\gamma
\nc\de\delta
\nc\ve\varepsilon
\nc\lam\lambda
\nc\om\omega
\nc\ze\zeta
\nc\si\sigma
\nc\Si\Sigma
\nc\Ga\Gamma
\nc\De\Delta
\nc\Om\Omega
\nc\io\iota
\nc\ka\kappa

\nc\nab\nabla

\nc\vt\vartheta
\rnc\th\theta
\nc\Th\Theta

\nc\BY{\qquad\text{by}\qquad}
\nc\GIVENBY{\qquad\text{given by}\qquad}
\nc\ISGIVENBY{\qquad\text{is given by}\qquad}
\nc\OR{\qquad\text{or}\qquad}
\nc\Or{\quad\text{or}\quad}
\nc\AND{\qquad\text{and}\qquad}
\nc\ANDSIM{\qquad\text{and similarly}\qquad}
\nc\ANDSIm{\quad\text{and similarly}\quad}
\nc\ANDSO{\qquad\text{and so}\qquad}
\nc\ANd{\quad\text{and}\quad}
\nc\COMMA{,\qquad}
\nc\COMMa{,\quad}
\nc\WHERE{\qquad\text{where}\qquad}

\rnc\iff{\ \Leftrightarrow\ }
\nc\IFf{\quad \Leftrightarrow\quad }
\nc\Iff{\ \ \Leftrightarrow\ \ }
\nc\IFF{\qquad \Leftrightarrow\qquad }
\rnc\implies{\ \Rightarrow\ }
\nc\IMPLIES{\qquad \Rightarrow\qquad }

\nc\set[2]{\{#1:#2\}}
\nc\bigset[2]{\big\{#1:#2\big\}}
\nc\pres[2]{\la#1:#2\ra}

\nc\bit{\begin{itemize}[label=\textbullet, leftmargin=5mm]}
\nc\eit{\end{itemize}}
\nc\ben{\begin{enumerate}[label=\textup{(\roman*)},leftmargin=10mm]}
\nc\bena{\begin{enumerate}[label=\textup{(\alph*)},leftmargin=10mm]}
\nc\een{\end{enumerate}}
\nc\bmc{\begin{multicols}}
\nc\emc{\end{multicols}}

\nc\pf{\begin{proof}}
\nc\epf{\end{proof}}
\nc\pfclaim{\begin{quote}\begin{proof}}
\nc\epfclaim{\end{proof}\end{quote}}
\nc\epfres{\hfill\qed}
\nc\epfreseq{\tag*{\qed}}
\let\oldproofname=\proofname
\renewcommand{\proofname}{\rm\bf{\oldproofname}}
\nc{\pfitem}[1]{\medskip \noindent #1.}
\nc{\firstpfitem}[1]{#1.}
\nc{\pfcase}[1]{\medskip\noindent {\bf Case #1.}}
\nc\aftercases{\medskip\noindent}

\nc\rH{\mathrel{\H}}
\nc\rL{\mathrel{\L}}
\nc\rR{\mathrel{\R}}
\nc\rD{\mathrel{\D}}
\nc\rJ{\mathrel{\J}}
\nc\rK{\mathrel{\K}}
\nc\rsi{\mathrel{\si}}

\nc\leqF{\leq_{\F}}
\nc\geqF{\geq_{\F}}

\newcommand{\id}{\operatorname{id}}
\newcommand{\dom}{\operatorname{dom}}
\newcommand{\codom}{\operatorname{codom}}
\newcommand{\coker}{\operatorname{coker}}
\newcommand{\rank}{\operatorname{rank}}
\newcommand{\non}{\operatorname{non}}

\nc\mr\mathrel
\nc\pc[2]{(#1,#2)^\sharp}

\nc\U{{\bf U}}
\nc\bF{{\bf F}}
\nc\GG{{\bf G}}

\nc\V{\mathcal V}
\nc\G{\mathcal G}
\rnc\iff{\ \Leftrightarrow\ }
\rnc\implies{\ \Rightarrow\ }
\nc\Implies{\quad \Rightarrow\quad }
\nc\F{\mathrel{\mathscr F}}
\nc\C{\mathscr C}
\nc\M{\mathcal M}
\nc\CC{\mathcal C}
\nc\CP{\CC_P}
\nc\CPz{\CC_{P_0}}
\nc\EP{\bE(P)}
\nc\EPd{\bE(P')}
\nc\PE{\bP(E)}
\nc\PEd{\bP(E')}
\nc\DD{\mathcal D}
\nc\I{\mathcal I}
\rnc\O{\mathcal O}
\rnc\S{\mathcal S}
\rnc\P{\mathcal P}
\nc\TL{\mathcal T\!\mathcal L}
\nc\PP{\mathscr P\mathcal P}
\nc\T{\mathcal T}
\nc\p{\mathfrak p}
\nc\q{\mathfrak q}
\rnc\r{\mathfrak r}
\nc\s{\mathfrak s}
\rnc\t{\mathfrak t}
\nc\bd{{\bf d}}
\nc\br{{\bf r}}
\nc\op\oplus
\nc\lra{\mathrel\leftrightarrow}
\nc\es\varnothing
\nc\rev{\textup{rev}}
\nc\corestt{{\upharpoonleft}}
\nc\restt{{\upharpoonright}}
\nc\corest{{\downharpoonleft}}
\nc\rest{{\upharpoonright}}
\nc\WHERe{\quad\text{where}\quad}
\nc\ot\leftarrow
\rnc\a{\mathfrak a}
\rnc\b{\mathfrak b}
\rnc\c{\mathfrak c}
\rnc\d{\mathfrak d}
\nc\sub\subseteq
\nc\mt\mapsto
\nc\im{\operatorname{im}}
\nc\B{\mathcal B}
\nc\PB{\P\B}
\nc\PT{\P\T}
\nc\E{\mathbb E}
\nc\Ef{\E^\flat}
\nc\BP{\operatorname{\textup{\textsf{BP}}}}
\rnc\SS{\operatorname{\textup{\textsf{S}}}}
\nc\MM{\operatorname{\textup{\textsf{M}}}}

\numberwithin{equation}{section}

\newtheorem{thm}[equation]{Theorem}
\newtheorem{lemma}[equation]{Lemma}
\newtheorem{cor}[equation]{Corollary}
\newtheorem{prop}[equation]{Proposition}

\theoremstyle{definition}

\newtheorem{rem}[equation]{Remark}
\newtheorem{eg}[equation]{Example}

\newcounter{caseco}

\newcounter{subcaseco}

\newcounter{stepco}

\newcounter{stageco}

\makeatletter
\def\blfootnote{\gdef\@thefnmark{}\@footnotetext}
\makeatother

\begin{document}

\title{Involutions of (twisted) diagram monoids}

\date{}

\author{James East\footnote{Partially supported by ARC Future Fellowship FT190100632.} \ and P.A.~Azeef Muhammed\\[3mm]
{\it\small Centre for Research in Mathematics and Data Science,}\\
{\it\small Western Sydney University, Locked Bag 1797, Penrith NSW 2751, Australia.}\\[3mm]
{\tt\small J.East@WesternSydney.edu.au}, {\tt\small A.ParayilAjmal@WesternSydney.edu.au}}

\maketitle

\vspace{-5mm}

\maketitle

\begin{abstract}
\noindent
We classify the involutions of all of the most well-studied diagram monoids---namely the partition, planar partition, partial Brauer, Motzkin, Brauer and Temperley--Lieb monoids---and characterise those that give rise to star-regular or regular star-monoid structures.  We then complete the same program for the associated twisted diagram monoids, with respect to both the canonical float-counting twisting, and the recently-discovered rank-based twisting.  This necessitates developing a general theory of involutions of twisted products.  Some of our results were quite unexpected.  For example, a Brauer monoid is star-regular with respect to many of its involutions, but only a regular star-monoid for one of them.  We will also see that twisted diagram monoids over the integers are always star-regular, thereby providing new and very natural examples of star-regular monoids.  Along the way, we also obtain (by necessity) a number of results of independent interest; specifically, we classify the automorphisms of the Motzkin and Temperley--Lieb monoids (and hence also of the planar partition monoids), and we show that all of our diagram monoids generically have trivial centre.

\medskip

\noindent
\emph{Keywords}: Twistings, twisted products, diagram monoids, twisted diagram monoids, involutions, star-regular monoids, regular star-monoids.
\medskip

\noindent
MSC: 
20M10,  
20M17,    
20M20.  
\end{abstract}

\tableofcontents

\section{Introduction}\label{sect:intro}

An \emph{involution} of a multiplicative structure (such as a ring or semigroup) is a self-inverse self-map that reverses the product, i.e.~a map $a\mt a^*$ obeying the laws $a^{**}=a$ and $(ab)^* = b^*a^*$.  Canonical examples include inversion in groups, or transposition of matrices.  Involutions provide an important left-right symmetry, and have been exploited in many algebraic studies, one notable instance being Graham and Lehrer's work on \emph{cellular algebras} \cite{GL1996}.  Key motivating examples in~\cite{GL1996} included the \emph{Brauer} and \emph{Temperley--Lieb algebras}.
These are important classes of \emph{diagram algebras}, which are fundamental objects of study in many branches of mathematics and science, including representation theory, topology, invariant theory, logic, category theory, statistical mechanics and phylogenetics \cite{Jones1994_2,Martin1994,Brauer1937,TL1971,Abramsky2008,LZ2015,HR2005,Kauffman1990,FJ2022}.
The standard involution in a diagram algebra is given by a vertical reflection of diagrams.  (Examples of the kinds of diagrams in play can be seen in Figure \ref{fig:P6}.)  The larger \emph{partition algebras}~\cite{Jones1994_2,Martin1994} were shown to be cellular by Xi \cite{Xi1999}, before a unified framework was provided by Wilcox in \cite{Wilcox2007}.  The key idea in \cite{Wilcox2007} was to construct diagram algebras as \emph{twisted semigroup algebras} of underlying diagram \emph{monoids}, building on ideas implicit in the work of Jones and Kauffman \cite{Jones1983_2,Jones1987,Kauffman1990,Kauffman1997}, and extending techniques introduced in \cite{East2006}.

The canonical involution on a diagram monoid is still given by vertical reflection, but it has an additional important property.  Namely, it maps a diagram to a \emph{(von Neumann) inverse}, meaning that we also have the laws $a=aa^*a$ and $a^*=a^*aa^*$.  (In a diagram \emph{algebra}, $aa^*a$ is generally a scalar multiple of $a$.)  Consequently, diagram monoids are examples of \emph{regular $*$-monoids} in the sense of Nordahl and Scheiblich \cite{NS1978}.  This extra layer of structure has been instrumental in many studies, going back to \cite{EF2012,FL2011}; for some recent examples see \cite{CEM2026,EPA2024,EGPR2025a,EGPR2025,EMRT2018}.

Diagram algebras are also (traditional) semigroup algebras of associated \emph{twisted diagram monoids}.  These monoids are typically infinite, and include information about the number of floating components formed when diagrams are stacked together (again see Figure \ref{fig:P6}).  The twisted monoids still have an involution, but it no longer maps elements to (von Neumann) inverses, as noted in \cite[Remark 1]{KV2023}.  In fact, the most well-studied families of twisted monoids in the literature \cite{BDP2002,DE2018,DE2017,ER2022c,ER2022b,KV2023,LF2006,KV2020,CHKLV2019, ACHLV2015} are not even regular, meaning that inverses do not exist for some elements.  In recent work on pseudovarieties and identities, Volkov and his collaborators \cite{KV2023,KV2020,CHKLV2019, ACHLV2015,AV2020} have found it convenient to embed these monoids into larger \emph{regular} twisted monoids, via the introduction of so-called \emph{negative loops}.  Intriguing properties of the involutions of these monoids formed one of the original motivations for the current work, and among other things we will show that although they are not regular $*$-monoids, they are \emph{$*$-regular}.  

The class of $*$-regular semigroups was introduced by Drazin \cite{Drazin1979}; such a semigroup has an involution $a\mt a^*$ with the property that $a^*$ is \emph{$\H$-related} to an inverse $a^+$ of $a$.  (Here $\H$ is one of Green's divisibility relations \cite{Green1951}; see below for precise definitions.)  As a special case, the regular $*$-semigroups of \cite{NS1978} are the $*$-regular semigroups in which $a^+=a^*$.  An original motivating example of a $*$-regular semigroup is the monoid of $n\times n$ real or complex matrices, with involution given by matrix transposition, and where $a^+$ is the \emph{Moore--Penrose inverse} of~$a$~\cite{Penrose1955,Moore1920}.  

Regular twisted partition monoids featured crucially in the recent paper \cite{EGPR2025}, where they arose naturally (and somewhat unexpectedly) in relation to the \emph{free idempotent-generated semigroup} over the ordinary partition monoid.  The key connection was that the two (non-isomorphic) monoids have isomorphic \emph{biordered sets of idempotents}.  This observation, and a number of other structural parallels that also hold in other families of twisted diagram monoids, then led to the article~\cite{EGPR2026}, where the idea of an abstract \emph{twisted product} was introduced.  Such a product has the form~${M\times_\Phi^qS}$, for monoids $M$ and $S$, and includes the direct product $M\times S$ as a degenerate case; full definitions are given below.  Traditional twisted diagram monoids arise when $S$ is an ordinary diagram monoid, $M=\N$ is the additive monoid of natural numbers, and $\Phi:S\times S\to\N$ is the `float-counting' twisting map alluded to above.  In this realisation, Volkov's regular extensions are obtained by replacing~$\N$ by~$\Z$, the additive group of integers.  More examples may be obtained by allowing~$M$ to vary further and/or replacing $\Phi$ by other twistings, such as the so-called \emph{rigid} twistings introduced in \cite{EGPR2026}.

As mentioned above, the original motivation for the current article was to explore the properties of the naturally-occurring involutions of twisted diagram monoids.  Such involutions typically `extend' an involution of the underlying un-twisted monoid, but it turns out that multiple extensions of the same involution can exist.  It also transpires that un-twisted diagram monoids have multiple involutions, with varying regularity properties.  Thus, two parallel programs emerge:
\begin{enumerate}[label=\textup{(\Roman*)},leftmargin=10mm]
\item \label{I} Classify the involutions of diagram monoids, and characterise those giving rise to star-regular or regular star-monoid structures.\footnote{We use the term `star-regular' to mean `$*$-regular with respect to some involution ${}^*$', and similarly for `regular star-'.  It is also linguistically convenient to use the terms `star-regularity' and `regular-starity'.}
\item \label{II} Do the same for the corresponding twisted diagram monoids.  
\een
As we will see, Program \ref{I} led to a number of very interesting results, some of which were quite unexpected.  It also naturally fed into Program \ref{II}, for which we additionally needed to develop some general theory concerning involutions of arbitrary twisted products.  The strongest results of \cite{EGPR2026} concern so-called \emph{tight} twistings.  While the canonical float-counting twisting is tight for many of the diagram monoids considered in the literature, key examples are lacking this property, including Motzkin and partial Brauer monoids.  Thus, a key part of our theoretical development is the identification of appropriate properties weaker than tightness that are satisfied by \emph{all} of our diagram monoids, but still allow us to prove strong general results.  The key such property is what we call \emph{unitality}, and this is an assumption in many of our main theoretical results.

The paper is organised as follows.  
We begin in Section \ref{sect:prelim} with preliminary material on semigroups and involutions, star-regular and regular star-semigroups, twisted products, and diagram monoids, including the definitions of the partition monoid $\P_n$, the planar partition monoid~$\PP_n$, the partial Brauer monoid $\PB_n$, the Motzkin monoid $\M_n$, the Brauer monoid $\B_n$, and the Temperley--Lieb monoid $\TL_n$ (a.k.a.~the Jones monoid $\mathcal J_n$).

Section \ref{sect:D} concerns the above-mentioned (un-twisted) diagram monoids, with the main results being the classification of involutions in Theorem \ref{thm:invD}; and the characterisation of star-regularity and regular-starity in Theorems~\ref{thm:Pn*}--\ref{thm:PPn*}.  Several results are given for these regularity properties, as the situations are rather different from monoid to monoid.  For example,~$\P_n$ is only star-regular with respect to the canonical involution (given by vertical reflection), which as we have already mentioned gives it the stronger structure of a regular star-monoid.  By contrast, $\B_n$ is star-regular for \emph{many} (but not all) involutions, yet regular-starity holds only for the canonical one.  The proof of the classification involves an application of a result from~\cite{EN2016}, which in turn requires knowledge of the automorphisms of the monoids.  These were classified for $\P_n$, $\PB_n$ and~$\B_n$ by Mazorchuk~\cite{Maz2002}.  For $\PP_n$, $\M_n$ and $\TL_n$ we provide the classification here in Theorems~\ref{thm:AutTLn} and~\ref{thm:AutPPnMn}.  For later application, we also show that our diagram monoids generically have trivial centres in Theorem \ref{thm:Z}.  We have not seen this result stated elsewhere in the literature, and we believe it is of independent interest, as is the classification of the automorphisms of the planar monoids.

Sections~\ref{sect:*} and~\ref{sect:*reg} contain the main theoretical developments of the paper, giving a sequence of general results on involutions, star-regularity and regular-starity of an arbitrary twisted product $T = M\times_\Phi^qS$.  These results require mild assumptions on $\Phi$ that are introduced and discussed in Section \ref{subsect:T}; each is a consequence of tightness.  
In Section \ref{sect:*} we identify a special kind of involution of $T$, called \emph{pure}, and classifiy these in Theorem \ref{thm:inv1}.  Arbitrary involutions are classified in Theorem \ref{thm:inv2} in terms of four maps $S\to S$, $M\to M$, $S\to M$ and $M\to S$.  The first two of these are involutions of~$S$ and~$M$ in the case of pure involutions of~$T$, but need not be in general; see Examples \ref{eg:nonext} and \ref{eg:nonext2} for some non-pure involutions (of a specific twisted product).  Corollary~\ref{cor:inv} gives a simple sufficient condition for all involutions of $T$ to be pure, namely that~$S$ has trivial centre, hence our earlier interest in this property for diagram monoids.
In Section~\ref{sect:*reg} we characterise the involutions of $T=M\times_\Phi^qS$ (as classified in the above-mentioned theorems) that give it the structure of a star-regular or regular star-monoid.  Star-regularity is covered in Theorems \ref{thm:*reggen} and \ref{thm:*reg}, which treat the general and pure cases, respectively; both theorems include a formula for Moore--Penrose inverses.  Regular-starity is treated in Theorem~\ref{thm:reg*}; there is only one theorem here, as it turns out that~$T$ can only be a regular star-monoid with respect to a pure involution.

Finally, we apply all of this to twisted diagram monoids, treating the canonical float-counting twisting $\Phi$ in Section \ref{sect:T}, and an alternative rank-based twisting $\Psi$ in Section \ref{sect:Psi}.  
We first prove a result (Theorem~\ref{thm:ia*}) about an important special family of involutions of a general twisted product $M\times_\Phi^qS$ (with~$S$ being one of~$\P_n$,~$\PP_n$, $\PB_n$, $\M_n$,~$\B_n$ or~$\TL_n$), and then restrict our attention to the case that $M=\Z$ and $q=1$, where we obtain Volkov's regular twisted diagram monoids $T = \Z\times_\Phi^1S$.  We classify the involutions of $T$ in Theorem \ref{thm:invT}, showing in particular that these are all pure, except in the case of $S=\B_2$, when $T$ has two non-pure involutions (and two pure ones).  Star-regularity and regular-starity are characterised in Theorem \ref{thm:*regZD}.  In general,~$T$ is star-regular with respect to involutions extending star-regular involutions of $S$, but regular-starity of $T$ is very rare, only holding in three specific low-degree cases.  Intriguingly,~${\Z\times_\Phi^1\B_2}$ is star-regular with respect to both of its non-pure involutions.
The same program is carried out in Section \ref{sect:Psi} for products of the form $M\times_\Psi^qS$, arising from the rank-based twisting $\Psi$.  In general there are more involutions of such a product (see Theorem \ref{thm:invTPsi}) than for float-counting products.  Star-regularity and regular-starity are also more common (see Theorem \ref{thm:*regZDPsi}), and in particular we obtain new infinite families of star-regular and regular star-monoids.

\section{Preliminaries}\label{sect:prelim}

Here we gather the preliminary definitions and results we need, concerning general semigroups (Section \ref{subsect:M}), star-regular and regular-star semigroups (Section \ref{subsect:*}), twisted products (Section~\ref{subsect:T}) and diagram monoids (Section \ref{subsect:D}).  For more details, and for proofs of the various assertions, see \cite{Howie1995,CPbook,EGPR2026,EG2017,NS1978,Drazin1979,NP1985}.

\subsection{Semigroups and involutions}\label{subsect:M}

Let $S$ be a semigroup, and denote by $S^1$ the \emph{monoid completion} of $S$; so $S^1=S$ if $S$ is a monoid, and otherwise $S^1=S\cup\{1\}$ where $1$ is a symbol not in $S$, acting as an adjoined identity.  \emph{Green's (equivalence) relations} $\L$, $\R$ and $\J$ are defined in terms of principal left, right or two-sided ideals.  For $a,b\in S$ we have
\[
a\L b \iff S^1a=S^1b \COMMA
a\R b \iff aS^1=bS^1 \AND
a\J b \iff S^1aS^1=S^1bS^1.
\]
These can also be thought of as divisibility relations, as for example $a\L b$ is equivalent to having $a=b$ or else $a=sb$ and $b=ta$ for some $s,t\in S$.  Green's $\H$ and $\D$ relations are defined by
\[
{\H} = {\L}\cap{\R} \AND {\D} = {\L}\vee{\R},
\]
where the latter is the \emph{join} of $\L$ and $\R$ in the lattice of equivalences on $S$.  That is, $\D$ is the least equivalence containing both $\L$ and $\R$.  It is well known that ${\D}={\L}\circ{\R}={\R}\circ{\L}$, and that ${\J}={\D}$ if $S$ is finite.  If $S$ is commutative, then all of Green's relations coincide.  If $\K$ is any of Green's relations, and if $a\in S$, we write $K_a = \set{b\in S}{b\K a}$ for the $\K$-class of~$a$ in~$S$.  We will sometimes use superscripts to emphasise the underlying semigroup, so write $\L^S$ for the~$\L$ relation on $S$, $L_a^S$ for the $\L^S$-class of $a\in S$, and so on; see for example Lemma \ref{lem:Green}.

An \emph{idempotent} of $S$ is an element $e\in S$ satisfying $e^2=e$.  For any subset $A\sub S$, we write~$E(A)$ for the set of all idempotents in $A$.  An element $a\in S$ is \emph{(von Neumann) regular} if there exists $b\in S$ such that $a=aba$.  Note then that $ab$ and $ba$ are idempotents, and that $a\R ab$ and $a\L ba$.  Note also that with $c=bab$, we have $a=aca$ and $c=cac$; such elements $a$ and $c$ are said to be mutual \emph{inverses}.  Of course idempotents are regular.

If $S$ is a monoid, then an element $a\in S$ is a \emph{unit} if $ab=1=ba$ for some $b\in S$.  This element~$b$ is necessarily unique, and is denoted $a^{-1}$.  The set $G_S$ of all units is a subgroup of $S$, and in fact $G_S = H_1$ is the $\H$-class of~$1$.  Units are also clearly regular.

For a subset $U$ of the semigroup $S$, we denote the \emph{centraliser} of $U$ in $S$ by
\[
C_S(U) = \set{a\in S}{au=ua\text{ for all }u\in U}.
\]
The \emph{centre} of $S$ is $Z(S) = C_S(S)$.  

We write $\Aut(S)$ and $\Aut^-(S)$ for the sets of all automorphisms and anti-automorphisms of~$S$.  The sets $\Aut(S)$ and ${\Aut^\pm(S) = \Aut(S)\cup\Aut^-(S)}$ are both groups under composition.  Note that $\Aut^-(S)$ could be empty, whereas $\Aut(S)$ always contains at least the identity map~$\id_S$.  Note also that when $S$ is commutative, we have $\Aut(S)=\Aut^-(S)=\Aut^\pm(S)$.

As we have already noted, an \emph{involution} of $S$ is an anti-automorphism of order at most $2$, i.e.~a map $\io:S\to S$ satisfying
\[
(ab)\io = (b\io)(a\io) \AND a\io\io = a \qquad\text{for all $a,b\in S$.}
\]
As is customary, we will often denote involutions using a symbol in a superscript, so for example~${{}^*:S\to S:a\mt a^*}$, in which case the defining properties become $a^{**}=a$ and $(ab)^*=b^*a^*$.  At other times it will be more convenient to use mapping notation.  It is easy to show that every involution of a \emph{monoid} necessarily fixes the identity element.  Note also that the identity map~$\id_S$ is an involution if and only if $S$ is commutative.

If $\io$ is any fixed involution of $S$, then it is easy to see that
\[
\Aut^-(S) = \io\Aut(S) = \Aut(S)\io.
\]
Thus, every involution of $S$ has the form $\io\al$ for some automorphism $\al$, and to have $(\io\al)^2=\id$, we need $\io\al\io = \al^{-1}$; in particular, if $\io$ and $\al$ commute, then we need $\al=\al^{-1}$.  This demonstrates the following, which is part of \cite[Proposition 3.4]{EN2016}:

\begin{lemma}\label{lem:EN}
Let $S$ be a semigroup, and suppose $\io$ is an involution of $S$ that commutes with every automorphism of $S$.  Then the involutions of $S$ are precisely the compositions $\io\al$, where $\al\in\Aut(S)$ has order $\leq2$.  \epfres
\end{lemma}

Given a semigroup $S$, and a fixed element $a\in S$, there is an associative operation ${\star} = {\star_a}$ on~$S$ defined by $s\star t = sat$ for $s,t\in S$.  The semigroup $S^a = (S,\star)$ is called the \emph{variant} of~$S$ with respect to the \emph{sandwich element}~$a$.  These have been studied extensively; see for example \cite{Hickey1986,KL2001}.  Of particular interest to us is the fact that $S^a \cong S$ if $S$ is a group.  Indeed, if $S$ is a monoid and $a\in S$ a unit, then the map $s\mt sa$ is easily seen to be an isomorphism~${S^a\to S}$, in which case the identity element of $S^a$ is $a^{-1}$.

\subsection{Regular-star and star-regular semigroups}\label{subsect:*}

Let $S$ be a semigroup with a fixed involution ${}^*$, so that $(ab)^*=b^*a^*$ and $a^{**}=a$ for all $a,b\in S$.
\bit
\item We say that $S$ is a \emph{regular $*$-semigroup} if $a^*$ is an inverse of $a$, for each $a\in S$.
\item We say that $S$ is a \emph{$*$-regular semigroup} if $a^*$ is $\H$-related to an inverse of $a$, for each $a\in S$.  This inverse (in the $\H$-class of $a^*$) is necessarily unique, and is denoted $a^+$, and called the \emph{Moore--Penrose inverse} of $a$.
\eit
These classes were introduced by Nordahl and Scheiblich \cite{NS1978} and Drazin \cite{Drazin1979}, respectively.  Note that a regular $*$-semigroup is $*$-regular (with $a^+=a^*$), but not conversely in general.  Indeed, we will see several examples of $*$-regular semigroups that are not regular $*$-semigroups.  

A regular $*$-semigroup can be thought of as an algebra $(S,\cdot,{}^*)$, where $(S,\cdot)$ is a semigroup, and where ${}^*$ is a \emph{unary operation} satisfying the identities
\[
(ab)^* = b^*a^* \COMMA a^{**}=a \COMMA a=aa^*a \AND a^*=a^*aa^*.
\]
(Note that in fact the fourth identity follows from the second and third.)  Similarly, a $*$-regular semigroup is an algebra $(S,\cdot,{}^*,{}^+)$, where ${}^*$ and ${}^+$ are both unary operations satisfying the identities
\[
(ab)^* = b^*a^* \COMMa a^{**}=a \COMMa a=aa^+a \COMMa a^+=a^+aa^+ \COMMa (aa^+)^*=aa^+ \ANd (a^+a)^*=a^+a.
\]
The first two laws say that ${}^*$ is an involution, and the other four abstractly axiomatise the Moore--Penrose inverse as the unique $x\in S$ satisfying
\begin{equation}\label{eq:axa}
a=axa \COMMA x=xax \COMMA (ax)^*=ax \AND (xa)^*=xa.
\end{equation}
There are a number of other consequences of the above laws, e.g.~$a^{++}=a$ and $a^{*+} = a^{+*}$, although ${}^+$ is typically not itself an involution.  As noted above, a regular $*$-semigroup is a $*$-regular semigroup satisfying the identity $a^+=a^*$.

There are several other equivalent definitions of star-regular semigroups, one of which will be useful in the proof of Theorem \ref{thm:Bn*}.  To describe this, let~$S$ be a semigroup with an involution~${}^*$.  An element $p$ of $S$ is a \emph{$*$-projection} (or just a \emph{projection} if ${}^*$ is understood from context) if $p^2=p=p^*$, i.e.~if $p$ is an idempotent of $S$ that is fixed by the involution.  Then $S$ is $*$-regular if and only if every $\R$-class contains a $*$-projection, which is in turn equivalent to every $\L$-class containing a $*$-projection.  Viewing $S$ as the algebra $(S,\cdot,{}^*,{}^+)$, any element $a\in S$ satisfies $a\R aa^+$ and $a\L a^+a$, for the projections $aa^+$ and $a^+a$.

\newpage

\subsection{Twisted products}\label{subsect:T}

Throughout this section we fix the following:
\begin{itemize}
\item
a monoid $S$, written multiplicatively, with identity $1$,
\item
a commutative monoid $M$, written additively, with identity $0$, and
\item
an arbitrary element $q\in M$.
\end{itemize}
A \emph{twisting} of $S$ is a map $\Phi:S\times S\to \N$ (where $\N=\{0,1,2,\ldots\}$ is the additive monoid of natural numbers) satisfying the following:
\begin{enumerate}[label=\textup{\textsf{(T\arabic*)}},leftmargin=10mm]
\item \label{T1} $\Phi(a,b) + \Phi(ab,c) = \Phi(a,bc) + \Phi(b,c)$ for all $a,b,c\in S$.
\end{enumerate}
Given such a twisting $\Phi$, the \emph{twisted product} $M\times_\Phi^qS$ is the semigroup with underlying set $M\times S$, and operation
\[
(i,a)(j,b) = (i+j+\Phi(a,b)q,ab) \qquad\text{for $i,j\in M$ and $a,b\in S$.}
\]
Associativity follows from \ref{T1}.  If we write $\Phi(a,b,c)$ for the common value in \ref{T1}, then we have
\[
(i,a)(j,b)(k,c) = (i+j+k+\Phi(a,b,c)q,abc) \qquad\text{for all $i,j,k\in M$ and $a,b,c\in S$.}
\]

We say a twisting $\Phi:S\times S\to\N$ is \emph{tight} if it also satisfies:
\begin{enumerate}[label=\textup{\textsf{(T\arabic*)}},leftmargin=10mm]\addtocounter{enumi}{1}
\item \label{T2} For all $a,b\in S$ there exist $a',b'\in S$ such that $ab = a'b = ab'$ and $\Phi(a',b) = \Phi(a,b') = 0$.
\end{enumerate}
We then say the twisted product $M\times_\Phi^qS$ is itself tight.  One of the motivations for studying tight twisted products is that Green's relations on such a product are built from the corresponding relations on $M$ and $S$ in a very simple way.

With this in mind, we say a twisted product $T = M\times_\Phi^qS$ is \emph{Green-decomposable} if for each of Green's relations $\K$ we have
\[
(i,a)\K^T(j,b) \IFF i\K^Mj \ANd a\K^Sb \qquad\text{for all $i,j\in M$ and $a,b\in S$.}
\]
Note that Green-decomposability is a property of the product $T$ itself, not just of $\Phi$.  

\begin{lemma}\label{lem:Green}
Let $T = M\times_\Phi^qS$ be a twisted product.  If $\Phi$ is tight, or if $q$ is a unit, then $T$ is Green-decomposable.
\end{lemma}

\pf
When $\Phi$ is tight, the result is part of \cite[Theorem 5.1]{EGPR2026}, so for the rest of the proof we assume that $q$ is a unit.  We just give the details for ${\K}={\R}$.  The argument for $\L$ is symmetrical, and $\J$ is analogous.  For $\H$ and $\D$ we combine the cases of $\L$ and $\R$.

Suppose first that $(i,a)\R^T(j,b)$.  If $(i,a)=(j,b)$ then certainly $i\R^Mj$ and $a\R^Sb$.  Otherwise there exist $k,l\in M$ and $c,d\in S$ such that
\[
(i,a)=(j,b)(k,c) = (j+k+\Phi(b,c)q,bc) \AND (j,b)=(i,a)(l,d)=(i+l+\Phi(a,d)q,ad).
\]
Looking at the first coordinates we obtain $i\R^M j$; the second coordinates give $a\R^Sb$.

Conversely, suppose $i\R^Mj$ and $a\R^Sb$, so that
\[
i=j+k, \ \ j=i+l,\ \ a=bc\ \ \text{and}\ \ b=ad \qquad\text{for some $k,l\in M$ and $c,d\in S$.}
\]
A routine calculation then gives
\[
(j,b)(k-\Phi(b,c)q,c) 
= (i,a)
\AND
(i,a)(l-\Phi(a,d)q,d) = (j,b),
\]
so that $(i,a)\R^T(j,b)$.
\epf

\begin{rem}\label{rem:Green}
In the above proof, the assumption that $q$ is a unit was not used for the forward implication.  That is, in any twisted product $T=M\times_\Phi^qS$, and for each of Green's relations $\K$, we have
\[
(i,a)\K^T(j,b) \IMPLIES i\K^Mj \ANd a\K^Sb \qquad\text{for all $i,j\in M$ and $a,b\in S$.}
\]
\end{rem}

So Green-decomposability of a product $M\times_\Phi^qS$ is a consequence of tightness of $\Phi$, but is not equivalent to it.  
Throughout the paper, we will be interested in a number of other properties of twistings that are implied by tightness.  Our main motivation for introducing these properties is that they are satisfied by all of our twisted diagram monoids, even though some of them lack tightness.  Beginning with the most prominent, we say a twisting~$\Phi$ is \emph{unital} if it satisfies:
\begin{enumerate}[label=\textup{\textsf{(U\arabic*)}},leftmargin=10mm]
\item \label{U1} For all $a,b\in S$, if $a$ or $b$ is a unit then $\Phi(a,b)=0$.
\end{enumerate}
We say $\Phi$ is \emph{semi-unital} if it satisfies:
\begin{enumerate}[label=\textup{\textsf{(U\arabic*)}$^\flat$},leftmargin=10mm]
\item \label{U1'} For all $a\in S$, we have $\Phi(a,1)=0=\Phi(1,a)$.
\end{enumerate}
When $\Phi$ is (semi-)unital, we say that the twisted product $M\times_\Phi^qS$ is itself (semi-)unital.  It is easy to check that a semi-unital twisted product is a monoid, with identity $(0,1)$.  Of course~\ref{U1} implies \ref{U1'}, but these conditions are not equivalent.  For example, if $S=\{1,x\}$ is a cyclic group of order $2$ (so $x^2=1$), then there is a twisting $\Phi$ given by
\[
\Phi(1,1)=\Phi(x,1)=\Phi(1,x)=0 \AND \Phi(x,x)=1.
\]
This satisfies \ref{U1'} but not \ref{U1}.

We say a twisting $\Phi:S\times S\to\N$ is \emph{Green-invariant} if it satisfies the following two conditions:
\begin{enumerate}[label=\textup{\textsf{(G\arabic*)}},leftmargin=10mm]
\item \label{G1} For all $a,b,c\in S$, we have $a\L^S b \implies \Phi(a,c) = \Phi(b,c)$.
\item \label{G2} For all $a,b,c\in S$, we have $a\R^S b \implies \Phi(c,a) = \Phi(c,b)$.
\end{enumerate}

Among other things, the next result shows that unitality and Green-invariance are indeed consequences of tightness.

\begin{lemma}\label{lem:TUGA}\leavevmode
\ben
\item \label{TUGA1}  Any Green-invariant semi-unital twisting is unital.
\item \label{TUGA2}  Any tight twisting is Green-invariant and unital.
\een
\end{lemma}

\pf
\firstpfitem{\ref{TUGA1}}  Let $\Phi$ be a Green-invariant semi-unital twisting of $S$.  If $a,b\in S$, and if $a$ is a unit, then $a\L^S1$, so it follows from \ref{G1} and \ref{U1'} that $\Phi(a,b)=\Phi(1,b)=0$.  A similar argument shows that $\Phi(a,b)=0$ if $b$ is a unit.

\pfitem{\ref{TUGA2}}  Lemma 3.4 of \cite{EGPR2026} says that any tight twisting is Green-invariant and semi-unital.  This part now follows from part \ref{TUGA1}.
\epf

\begin{rem}\label{rem:TUGA}
Any twisting $\Phi$ of a monoid $S$ restricts to a twisting of any submonoid $U$ of $S$, which we denote here by $\Psi=\Phi|_{U\times U}:U\times U\to \N$.  If $\Phi$ is (semi-)unital, then so too is $\Psi$, since any unit of $U$ is a unit of $S$.  The same is true of Green-invariance, as for example, if $a,b,c\in U$ then
\[
a\L^Ub \implies a\L^Sb\implies \Phi(a,c)=\Phi(b,c) \implies \Psi(a,c)=\Psi(b,c).
\]
Tightness of $\Phi$ does \emph{not} necessarily imply tightness of $\Psi$, however.  Indeed, it was observed in \cite[Proposition 4.2]{EGPR2026} that tightness of the canonical twisting on a partition monoid $\P_n$ is not inherited by some of its key submonoids; we will discuss this in more detail in Section \ref{subsect:D}.
\end{rem}

We now gather some basic results on twistings that we will need in later sections.  The first concerns some useful consequences of unitality.

\begin{lemma}\label{lem:unit}
Let $\Phi$ be a unital twisting of $S$, and let $a,b,g\in S$, where $g$ is a unit.  Then
\begin{enumerate}[label=\textup{\textsf{(U\arabic*)}},leftmargin=10mm]\addtocounter{enumi}{1}
\item \label{U2} $\Phi(a,bg)=\Phi(a,b)=\Phi(ga,b)$,
\item \label{U3} $\Phi(a,b,g)=\Phi(a,b)=\Phi(g,a,b)$,
\item \label{U4} $\Phi(ag,b) = \Phi(a,g,b) = \Phi(a,gb)$,
\item \label{U5} $\Phi(g^{-1}ag,g^{-1}bg) = \Phi(a,b)$.
\end{enumerate}
\end{lemma}

\pf
\firstpfitem{\ref{U2} and \ref{U3}}  By definition, and using \ref{T1}, we have
\[
\Phi(a,b,g) = \Phi(a,b)+\Phi(ab,g) = \Phi(a,bg)+\Phi(b,g).
\]
Since $g$ is a unit, \ref{U1} gives $\Phi(ab,g)=0=\Phi(b,g)$, so the above reduces to 
\[
\Phi(a,b,g) = \Phi(a,b) = \Phi(a,bg).
\]
We obtain $\Phi(g,a,b) = \Phi(a,b) = \Phi(ga,b)$ in similar fashion.

\pfitem{\ref{U4}}  By definition, and using \ref{U1}, we have $\Phi(a,g,b) = \Phi(a,g) + \Phi(ag,b) = \Phi(ag,b)$.  An analogous calculation gives $\Phi(a,g,b) = \Phi(a,gb)$.

\pfitem{\ref{U5}}  By \ref{U4} we have $\Phi(ag,g^{-1}b) = \Phi(a,gg^{-1}b) = \Phi(a,b)$.  But by \ref{U2} we also have $\Phi(ag,g^{-1}b) = \Phi(g^{-1}ag,g^{-1}bg)$.
\epf

The domain of a twisting $\Phi$ (which need not be tight, unital, etc.)~can be extended to any number of coordinates.  For $k\geq0$ and $a_1,\ldots,a_k\in S$, we inductively define
\[
\Phi(a_1,\ldots,a_k) = \begin{cases}
0 &\text{if $k\leq1$}\\
\Phi(a_1,\ldots,a_{k-1}) + \Phi(a_1\cdots a_{k-1},a_k) &\text{if $k\geq2$.}
\end{cases}
\]
When $k=2$ or $3$, it is easy to see that this agrees with the numbers $\Phi(a,b)$ and $\Phi(a,b,c)$ we have already encountered.  In the twisted product~${M\times_\Phi^qS}$, we have
\[
(i_1,a_1)\cdots(i_k,a_k) = (i_1+\cdots+i_k+\Phi(a_1,\ldots,a_k)q,a_1\cdots a_k).
\]

\begin{lemma}\label{lem:a1ak}
Let $\Phi$ be a twisting of $S$, and let $a_1,\ldots,a_m\in S$.
\ben
\item \label{a1ak1} For any $1\leq k\leq m$, we have
\[
\Phi(a_1,\ldots,a_m) = \Phi(a_1,\ldots,a_k) + \Phi(a_{k+1},\ldots,a_m) + \Phi(a_1\cdots a_k,a_{k+1}\cdots a_m).
\]
\item \label{a1ak2} For any $1\leq k\leq l\leq m$, we have
\begin{align*}
\Phi(a_1,\ldots,a_m) = \Phi(a_1,\ldots,a_k) + \Phi(a_{k+1},\ldots,a_l) &+ \Phi(a_{l+1},\ldots,a_m) \\
&+ \Phi(a_1\cdots a_k,a_{k+1}\cdots a_l,a_{l+1}\cdots a_m).
\end{align*}
\een
\end{lemma}

\pf
We just prove \ref{a1ak1}, as the proof of \ref{a1ak2} is analogous.  We can prove \ref{a1ak1} directly, using the definition.  Alternatively, we can prove it by noting that in $\N\times_\Phi^1S$ we have
\begin{align*}
(\Phi(a_1,\ldots,a_m),a_1\cdots a_m) &= (0,a_1)\cdots(0,a_k) \cdot (0,a_{k+1})\cdots(0,a_m)\\
&= (\Phi(a_1,\ldots,a_k),a_1\cdots a_k) \cdot (\Phi(a_{k+1},\ldots,a_l),a_{k+1}\cdots a_l).
\end{align*}
The result now follows by expanding this last product, and comparing first coordinates.
\epf

Later on we will make use of the following special cases, which can all be deduced from Lemma \ref{lem:a1ak}:
\begin{align}
\nonumber \Phi(a,b,c,d) &= \Phi(a,bcd) + \Phi(b,c,d) = \Phi(a,b,c) + \Phi(abc,d)\\
\nonumber &= \Phi(a,b) + \Phi(ab,c,d) = \Phi(a,b,cd) + \Phi(c,d)\\
\label{eq:abcd} &= \Phi(b,c) + \Phi(a,bc,d) = \Phi(a,b) + \Phi(ab,cd) + \Phi(c,d) &&\text{for any $a,b,c,d\in S$.}
\end{align}
For the next statement, recall that $a,b\in S$ are mutual inverses if $a=aba$ and $b=bab$.

\begin{lemma}\label{lem:aba}
Let $\Phi$ be a twisting of $S$, and suppose $a,b\in S$ are mutual inverses.  Then
\ben
\item \label{aba1} $\Phi(ab,a) = \Phi(ab,ab) = \Phi(b,ab)$,
\item \label{aba2} $\Phi(a,b,a)=\Phi(b,a,b)$,
\item \label{aba3} $\Phi(a,b,a)=\Phi(b,a,b) = \Phi(a,b) + \Phi(b,a)$ if $\Phi$ is Green-invariant.
\een
\end{lemma}

\pf
\firstpfitem{\ref{aba1}}  Define the idempotent $e=ab$.  Using \ref{T1}, and noting that $ea=a$, we calculate
\[
\Phi(e,e)+\Phi(e,a) = \Phi(e,e)+\Phi(ee,a) = \Phi(e,ea)+\Phi(e,a) = \Phi(e,a)+\Phi(e,a),
\]
which yields $\Phi(e,e)=\Phi(e,a)$.  An analogous calculation gives $\Phi(e,e)=\Phi(b,e)$.  

\pfitem{\ref{aba2}}  Using \ref{aba1} we have 
\[
\Phi(a,b,a) = \Phi(a,b) + \Phi(ab,a) = \Phi(a,b) + \Phi(b,ab) = \Phi(b,ab) + \Phi(a,b) = \Phi(b,a,b).
\]

\pfitem{\ref{aba3}}  If $\Phi$ is Green-invariant, then since $ab\L^S b$, it follows from \ref{G1} that
\[
\Phi(a,b,a) = \Phi(a,b) + \Phi(ab,a) = \Phi(a,b) + \Phi(b,a).  \qedhere
\]
\epf

\subsection{Diagram monoids and their canonical involutions and twistings}\label{subsect:D}

Fix a non-negative integer $n$, and write $\bn=\{1,\ldots,n\}$, $\bn'=\{1',\ldots,n'\}$ and $\bn''=\{1'',\ldots,n''\}$.  The \emph{partition monoid} $\P_n$ consists of all set partitions of $\bn\cup\bn'$, under a product described below.  A partition $a\in\P_n$ is identified with any (simple) graph on vertex set $\bn\cup\bn'$ whose connected components are the blocks of~$a$.  Such graphs are drawn with un-dashed vertices on a top row, ordered $1<\cdots<n$, dashed vertices on a bottom row, ordered $1'<\cdots<n'$, and edges inside the rectangle bounded by the vertices.  

To define the product of partitions $a,b\in\P_n$, we first define three additional graphs:
\bit
\item $a^\downarrow$ on vertex set $\bn\cup\bn''$, obtained by changing every lower vertex $x'$ of $a$ to $x''$,
\item $b^\uparrow$ on vertex set $\bn''\cup\bn'$, obtained by changing every upper vertex $x$ of $b$ to $x''$,
\item $\Pi(a,b)$ on vertex set $\bn\cup\bn''\cup\bn'$, whose edge set is the union of the edge sets of $a^\downarrow$ and $b^\uparrow$.
\eit
We call $\Pi(a,b)$ the \emph{product graph} of $a$ and $b$; it is drawn with vertices $1''<\cdots<n''$ on a middle row.  The product $ab\in\P_n$ is then the unique partition with the property that $x,y\in\bn\cup\bn'$ belong to the same block of $ab$ if and only if they belong to the same connected component of~$\Pi(a,b)$.  An example calculation is shown in Figure \ref{fig:P6}, for the partitions
\begin{align*}
a &= \big\{\{1,4\},\{2,3,4',5'\},\{5,6\},\{1',2',6'\},\{3'\}\big\} ,\\
b &= \big\{\{1,2\},\{3,4,1'\},\{5,4',5',6'\},\{6\},\{2',3'\}\big\} , \\
ab &= \big\{\{1,4\},\{2,3,1',4',5',6'\},\{5,6\},\{2',3'\}\big\},
\end{align*}
all from $\P_6$.
The identity of $\P_n$ is $1 = \custpartn{1,2,4}{1,2,4}{
\stline11
\stline22
\stline44
\udotted24
\ddotted24
}$, and the group of units is (isomorphic to) the symmetric group $\S_n$.  A permutation $f\in\S_n$ is identified with the partition ${\bigset{\{x,(xf)'\}}{x\in\bn}}$.

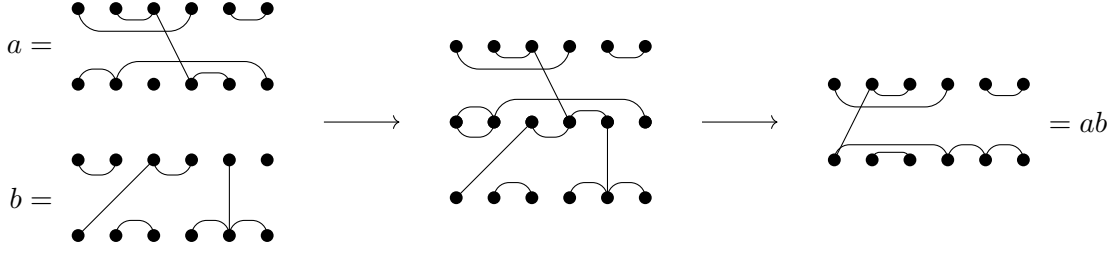
\begin{figure}[ht]
\begin{center}
\begin{tikzpicture}[scale=.5]

\begin{scope}[shift={(0,0)}]	
\uvs{1,...,6}
\lvs{1,...,6}
\uarcx14{.6}
\uarcx23{.3}
\uarcx56{.3}
\darc12
\darcx26{.6}
\darcx45{.3}
\stline34
\draw(0.6,1)node[left]{$a=$};
\draw[->](7.5,-1)--(9.5,-1);
\end{scope}

\begin{scope}[shift={(0,-4)}]	
\uvs{1,...,6}
\lvs{1,...,6}
\uarc12
\uarc34
\darc45
\darc56
\darc23
\stline31
\stline55
\draw(0.6,1)node[left]{$b=$};
\end{scope}

\begin{scope}[shift={(10,-1)}]	
\uvs{1,...,6}
\lvs{1,...,6}
\uarcx14{.6}
\uarcx23{.3}
\uarcx56{.3}
\darc12
\darcx26{.6}
\darcx45{.3}
\stline34
\draw[->](7.5,0)--(9.5,0);
\end{scope}

\begin{scope}[shift={(10,-3)}]	
\uvs{1,...,6}
\lvs{1,...,6}
\uarc12
\uarc34
\darc45
\darc56
\stline31
\stline55
\darc23
\end{scope}

\begin{scope}[shift={(20,-2)}]	
\uvs{1,...,6}
\lvs{1,...,6}
\uarcx14{.6}
\uarcx23{.3}
\uarcx56{.3}
\darc14
\darc45
\darc56
\stline21
\darcx23{.2}
\draw(6.4,1)node[right]{$=ab$};
\end{scope}

\end{tikzpicture}
\caption{Multiplication of partitions $a,b\in\P_6$, with the product graph $\Pi(a,b)$ in the middle.  Here we have $\Phi(a,b)=1$ because there is a single floating component in $\Pi(a,b)$, namely $\{1'',2'',6''\}$.}
\label{fig:P6}
\end{center}
\end{figure}

An \emph{upper block} of a partition $a\in\P_n$ is a block containing only un-dashed vertices.  \emph{Lower blocks} are defined analogously.  A \emph{transversal} of $a$ is a block containing both dashed and un-dashed vertices.  The upper and lower blocks are collectively referred to as \emph{non-transversals}.  The \emph{rank} of $a$, denoted $\rank(a)$, is defined to be the number of transversals; this can be any integer from $0$ to $n$.  We also define the (\emph{co})\emph{domain} and (\emph{co})\emph{kernel} of $a$ to be
\begin{align}
\nonumber \dom(a) &= \set{x\in\bn}{x\text{ belongs to a transversal of }a},\\
\nonumber \codom(a) &= \set{x\in\bn}{x'\text{ belongs to a transversal of }a},\\
\nonumber \ker(a) &= \set{(x,y)\in\bn\times\bn}{x\text{ and }y\text{ belong to the same block of }a},\\
\label{eq:dom} \coker(a) &= \set{(x,y)\in\bn\times\bn}{x'\text{ and }y'\text{ belong to the same block of }a}.
\end{align}
So $\dom(a)$ and $\codom(a)$ are subsets of $\bn$, while $\ker(a)$ and $\coker(a)$ are equivalence relations.  For example, with $a\in\P_6$ as in Figure \ref{fig:P6} we have $\rank(a)=1$, $\dom(a)=\{2,3\}$, $\codom(a)=\{4,5\}$, ${\bf6}/\ker(a)=\big\{\{1,4\},\{2,3\},\{5,6\}\big\}$ and ${\bf6}/\coker(a)=\big\{\{1,2,6\},\{3\},\{4,5\}\big\}$.

The above parameters can be used to characterise Green's relations on $\P_n$ \cite{FL2011,Wilcox2007}.  For example, we have
\begin{equation}\label{eq:RJ}
a \R b \iff \dom(a)=\dom(b) \text{ and }\ker(a)=\ker(b),
\AND
a \J b \iff \rank(a) = \rank(b).
\end{equation}

The partition monoid $\P_n$ has a canonical involution ${}^*:\P_n\to\P_n:a\mt a^*$ obtained by interchanging dashed and un-dashed vertices.  Diagrammatically, $a^*$ is obtained by reflecting $a$ in a horizontal axis.  For example, with $a= \custpartn{1,...,6}{1,...,6}{
\uarcx14{.6}
\uarcx23{.3}
\uarcx56{.3}
\darc12
\darcx26{.6}
\darcx45{.3}
\stline34}\in\P_6$ as in Figure \ref{fig:P6}, we have ${a^* = \custpartn{1,...,6}{1,...,6}{
\darcx14{.6}
\darcx23{.3}
\darcx56{.3}
\uarc12
\uarcx26{.6}
\uarcx45{.3}
\stline43}}$.  We will see in later sections that $\P_n$ has many more involutions, but ${}^*$ is the unique one giving it a star-regular structure.  In fact, it is well known that $\P_n$ is a regular star-monoid with respect to this involution, and it is easy to check that the defining identities hold:
\[
a^{**}=a=aa^*a \AND (ab)^* = b^*a^* \qquad\text{for $a,b\in\P_n$.}
\]

Consider partitions $a,b\in\P_n$.  A connected component of the product graph $\Pi(a,b)$ containing only double-dashed vertices is called a \emph{floating component}.  The number of such components is denoted $\Phi(a,b)$.  So $\Phi$ is a mapping $\P_n\times\P_n\to\N$, and it turns out that this is a tight twisting, which we call the \emph{canonical float-counting twisting}; see \cite[Lemma 4.1]{FL2011} and \cite[Lemma~2.9(i)]{ER2022b}, which prove properties \ref{T1} and \ref{T2}, respectively.  For example, we have $\Pi(a,b)=1$ for $a,b\in\P_6$ from Figure \ref{fig:P6}, with the unique floating component of the product graph $\Pi(a,b)$ being~$\{1'',2'',6''\}$.
There is another natural \emph{rank-based} twisting of $\P_n$, whose definition we will recall in Section \ref{sect:Psi}.

The next two lemmas record some simple interactions between the canonical involution and twisting.

\begin{lemma}\label{lem:*symPn}
For all $a,b\in\P_n$ we have $\Phi(a,b) = \Phi(b^*,a^*)$.
\end{lemma}

\pf
This follows quickly by considering the relationships between the product graphs $\Pi(a,b)$ and $\Pi(b^*,a^*)$
\epf

\begin{lemma}\label{lem:Phiaa*a}
If $a\in\P_n$, then 
\ben
\item \label{Phiaa*a1} $\Phi(a,a^*)$ is equal to the number of lower blocks of $a$,
\item \label{Phiaa*a2} $\Phi(a^*,a)$ is equal to the number of upper blocks of $a$,
\item \label{Phiaa*a3} $\Phi(a,a^*,a)$ is equal to the number of non-transversals of $a$.
\een
\end{lemma}

\pf
\firstpfitem{\ref{Phiaa*a1} and \ref{Phiaa*a2}}  These are easily checked, given the form of $\Pi(a,a^*)$ and $\Pi(a^*,a)$.  

\pfitem{\ref{Phiaa*a3}}  We have already noted that $\Phi$ is tight, and hence Green-invariant by Lemma \ref{lem:TUGA}\ref{TUGA2}.  So by Lemma \ref{lem:aba}\ref{aba3} we have ${\Phi(a,a^*,a) = \Phi(a,a^*) + \Phi(a^*,a)}$.  This part now follows from~\ref{Phiaa*a1} and~\ref{Phiaa*a2}.
\epf

The numbers $\Phi(a,a^*,a)$ will play an important role in later sections, so we introduce the following notation for them:
\[
\non(a) = \Phi(a,a^*,a) = {} \text{the number of non-transversals of } a\in\P_n.
\]

Next we recall the definitions of five important submonoids of $\P_n$.
\bit
\item
The \emph{planar partition monoid} $\PP_n$ consists of all \emph{planar} partitions, by which we mean those corresponding to some graph that can be drawn with no edge crossings, remembering that edges are always drawn within the rectangle bounded by the vertices.  For the partitions $a,b\in\P_6$ from Figure \ref{fig:P6} we have $b\in\PP_6$ but $a\not\in\PP_6$.
\item 
The \emph{partial Brauer monoid}~$\PB_n$ consists of all partitions whose blocks have size at most $2$.  Note that each partial Brauer partition $a\in\PB_n$ is represented by a \emph{unique} (simple) graph on vertex set $\bn\cup\bn'$.
\item 
The \emph{Brauer monoid}~$\B_n$ consists of all partitions whose blocks have size (exactly) $2$.
\item 
The \emph{Motzkin monoid}~$\M_n = \PB_n\cap\PP_n$ consists of all planar partial Brauer partitions.
\item 
The \emph{Temperley--Lieb monoid}~$\TL_n = \B_n\cap\PP_n$ consists of all planar Brauer partitions.
\eit
Note that the symmetric group $\S_n$ is contained in $\B_n$, and hence also in $\PB_n$, and is consequently the group of units of these monoids.  On the other hand, $\PP_n$, $\M_n$ and $\TL_n$ are all $\H$-trivial, and hence have trivial group of units.  Figure \ref{fig:submonoids} shows the containments among the submonoids we have discussed, and gives representative examples of each.
Each of $\PP_n$, $\PB_n$, $\B_n$, $\M_n$ and~$\TL_n$ is clearly closed under the involution ${}^*$, so these are all regular $*$-monoids.  Moreover, the canonical twisting $\Phi$ restricts to a twisting of each submonoid, which we will also denote by~$\Phi$.  These twistings are tight for~$\PP_n$,~$\B_n$ and~$\TL_n$ (i.e.~\ref{T2} holds in these three submonoids), but \emph{not} for $\PB_n$ and $\M_n$ when $n\geq2$; see~\cite[Proposition 4.2]{EGPR2026}.  On the other hand since the twisting of $\P_n$ is tight, it is also unital and Green-invariant (by Lemma \ref{lem:TUGA}), and hence the twistings on \emph{all} of the above submonoids are unital and Green-invariant (see Remark \ref{rem:TUGA}).

\begin{figure}[ht]
\begin{center}
\begin{tikzpicture}[scale=1.05]
\node[rounded corners,rectangle,draw] (P) at (3,8) {$\P_n$};
\node[rounded corners,rectangle,draw] (PB) at (3,6) {$\PB_n$};
\node[rounded corners,rectangle,draw] (B) at (3,4) {$\B_n$};
\node[rounded corners,rectangle,draw] (S) at (3,2) {$\S_n$};
\node[rounded corners,rectangle,draw] (PP) at (6,6) {$\PP_n$};
\node[rounded corners,rectangle,draw] (M) at (6,4) {$\M_n$};
\node[rounded corners,rectangle,draw] (TL) at (6,2) {$\TL_n$};
\node[rounded corners,rectangle,draw] (1) at (6,0) {$\{1\}$};
\draw (P)--(PB)--(B)--(TL)--(M)--(PP)--(P) (PB)--(M) (B)--(S)--(1)--(TL);
\end{tikzpicture}
\qquad\qquad\qquad\qquad
\begin{tikzpicture}[scale=1]
\node[rounded corners,rectangle,draw] (P) at (3,8) {$\custpartn{1,2,3,4,5,6}{1,2,3,4,5,6}{\uarcx14{.6}\uarcx23{.3}\uarcx56{.3}\darc12\darcx26{.6}\darcx45{.3}\stline34}$};
\node[rounded corners,rectangle,draw] (PB) at (3,6) {$\custpartn{1,2,3,4,5,6}{1,2,3,4,5,6}{\uarcx13{.5}\uarc56\darc45\stline23}$};
\node[rounded corners,rectangle,draw] (B) at (3,4) {$\custpartn{1,2,3,4,5,6}{1,2,3,4,5,6}{\stline13\stline42\uarc23\uarc56\darcx16{.8}\darc45}$};
\node[rounded corners,rectangle,draw] (PP) at (6,6) {$\custpartn{1,2,3,4,5,6}{1,2,3,4,5,6}{\uarcx14{.6}\uarcx23{.2}\darc45\darc12\darc23\stline43\stline11\stline56}$};
\node[rounded corners,rectangle,draw] (M) at (6,4) {$\custpartn{1,2,3,4,5,6}{1,2,3,4,5,6}{\uarcx13{.5}\uarc56\darc45\stline43}$};
\node[rounded corners,rectangle,draw] (TL) at (6,2) {$\custpartn{1,2,3,4,5,6}{1,2,3,4,5,6}{\uarcx12{.4}\uarc45\darc34\darcx25{.8}\stlines{3/1,6/6}}$};
\node[rounded corners,rectangle,draw] (S) at (3,2) {$\custpartn{1,2,3,4,5,6}{1,2,3,4,5,6}{\stline13\stline21\stline32\stline44\stline56\stline65}$};
\node[rounded corners,rectangle,draw] (1) at (6,0) {$\custpartn{1,2,3,4,5,6}{1,2,3,4,5,6}{\stline11\stline22\stline33\stline44\stline55\stline66}$};
\draw (P)--(PB)--(B)--(TL)--(M)--(PP)--(P) (PB)--(M) (B)--(S)--(1)--(TL);
\end{tikzpicture}
\caption{Submonoids of $\P_n$ (left) and representative elements from each submonoid (right).}
\label{fig:submonoids}
\end{center}
\end{figure}
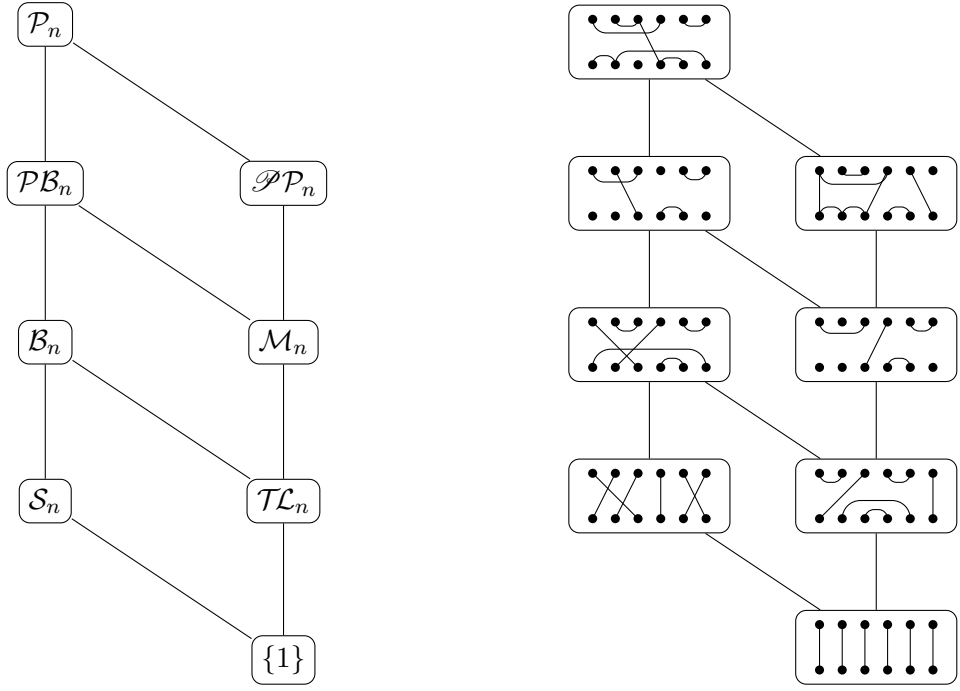

It is well known that the planar partition monoid $\PP_n$ is isomorphic to the even-degree Temperley--Lieb monoid $\TL_{2n}$ \cite{HR2005,Jones1994_2}, with an explicit isomorphism $\PP_n\to\TL_{2n}:a\mt\ol a$ given by `tracing' the blocks of $a$, as shown by example in Figure \ref{fig:PPnTL2n}; for more details see \cite[p.~989]{EMRT2018}.

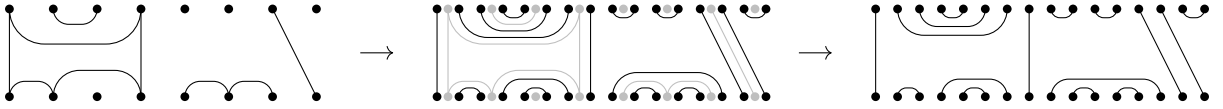
\begin{figure}[ht]
\begin{center}
\begin{tikzpicture}[scale=0.58]
\uarcx14{.8}
\uarcx23{.35}
\darcx12{.35}
\darcx24{.6}
\darcx56{.35}
\darcx67{.35}
\stlines{1/1,4/4,7/8}
\foreach \x in {1,...,8} {
\fill (\x,2)circle(.1); \fill (\x,0)circle(.1);
}
\draw[->] (9,1)--(9.75,1);
\begin{scope}[shift={(10,0)}]
\uarcxx{1.25}{3.75}{.65}{black}
\uarcxx{1.75}{3.25}{.5}{black}
\uarcxx{2.25}{2.75}{.2}{black}
\uarcxx{4.75}{5.25}{.2}{black}
\uarcxx{5.75}{6.25}{.2}{black}
\uarcxx{7.75}{8.25}{.2}{black}
\uarcxx14{.8}{lightgray}
\uarcxx23{.35}{lightgray}
\darcxx{2.75}{3.25}{.2}{black}
\darcxx{1.25}{1.75}{.2}{black}
\darcxx{5.25}{5.75}{.2}{black}
\darcxx{6.25}{6.75}{.2}{black}
\darcxx{2.25}{3.75}{.4}{black}
\darcxx{4.75}{7.25}{.55}{black}
\stlinex{.75}{.75}{black}
\stlinex{4.25}{4.25}{black}
\stlinex{6.75}{7.75}{black}
\stlinex{7.25}{8.25}{black}
\darcxx12{.35}{lightgray}
\darcxx24{.6}{lightgray}
\darcxx56{.35}{lightgray}
\darcxx67{.35}{lightgray}
\stlinex11{lightgray}
\stlinex44{lightgray}
\stlinex78{lightgray}
\foreach \x in {1,...,8} {
\fill[lightgray] (\x,2)circle(.1); \fill[lightgray] (\x,0)circle(.1);
\fill (\x+.25,2)circle(.1);
\fill (\x-.25,2)circle(.1);
\fill (\x+.25,0)circle(.1);
\fill (\x-.25,0)circle(.1);
}
\draw[->] (9,1)--(9.75,1);
\end{scope}
\begin{scope}[shift={(20,0)}]
\uarcxx{1.25}{3.75}{.6}{black}
\uarcxx{1.75}{3.25}{.4}{black}
\uarcxx{2.25}{2.75}{.2}{black}
\uarcxx{4.75}{5.25}{.2}{black}
\uarcxx{5.75}{6.25}{.2}{black}
\uarcxx{7.75}{8.25}{.2}{black}
\darcxx{2.75}{3.25}{.2}{black}
\darcxx{1.25}{1.75}{.2}{black}
\darcxx{5.25}{5.75}{.2}{black}
\darcxx{6.25}{6.75}{.2}{black}
\darcxx{2.25}{3.75}{.4}{black}
\darcxx{4.75}{7.25}{.4}{black}
\stlinex{.75}{.75}{black}
\stlinex{4.25}{4.25}{black}
\stlinex{6.75}{7.75}{black}
\stlinex{7.25}{8.25}{black}
\foreach \x in {1,...,8} {
\fill (\x+.25,2)circle(.1);
\fill (\x-.25,2)circle(.1);
\fill (\x+.25,0)circle(.1);
\fill (\x-.25,0)circle(.1);
}
\end{scope}
\end{tikzpicture}
\caption{A planar partition $a\in\PP_8$ (left), with its corresponding Temperley--Lieb partition $\ol a\in\TL_{16}$ (right).}
\label{fig:PPnTL2n}
\end{center}
\end{figure}

The isomorphism $\PP_n\to\TL_{2n}:a\mt\ol a$ will be used on several occasions in the current paper, but we note that it does \emph{not} preserve the twistings, in the sense that we do not necessarily have $\Phi(a,b) = \Phi(\ol a,\ol b)$ for $a,b\in\PP_n$.  For example, with
\[
a = \custpartn{1,2,3,4}{1,2,3,4}{\stline11\stline44\uarc12\uarc23\uarc34\darc12\darc23\darc34} \in \PP_4
\AND
\ol a = \custpartn{1,2,3,4,5,6,7,8}{1,2,3,4,5,6,7,8}{\stline11\stline88\uarc23\uarc45\uarc67\darc23\darc45\darc67} \in \TL_8,
\]
we have $\Phi(a,a)=0$ but $\Phi(\ol a,\ol a) = 3$.  This has the consequence that twisted products $M\times_\Phi^q\PP_n$ and $M\times_\Phi^q\TL_{2n}$ are typically not isomorphic.

At certain times, it will be convenient to refer to families of specific elements of $\P_n$.  For $1\leq i\leq n-1$ and $1\leq j\leq n$, we define the partitions
\begin{equation}\label{eq:elements}
s_i = \custpartn{1,3,4,5,6,8}{1,3,4,5,6,8}{\stlines{1/1,3/3,4/5,5/4,6/6,8/8}\udotted13\udotted68\ddotted13\ddotted68\vertlab11\vertlab4i\vertlab8n} 
\COMMa
e_i = \custpartn{1,3,4,5,6,8}{1,3,4,5,6,8}{\stlines{1/1,3/3,6/6,8/8}\uarc45\darc45\udotted13\udotted68\ddotted13\ddotted68\vertlab11\vertlab4i\vertlab8n} 
\COMMa
h_i = \custpartn{1,3,4,5,6,8}{1,3,4,5,6,8}{\stlines{1/1,3/3,4/4,5/5,6/6,8/8}\uarc45\darc45\udotted13\udotted68\ddotted13\ddotted68\vertlab11\vertlab4i\vertlab8n} 
\COMMa
t_j = \custpartn{1,3,4,5,7}{1,3,4,5,7}{\stlines{1/1,3/3,5/5,7/7}\udotted13\udotted57\ddotted13\ddotted57\vertlab11\vertlab4j\vertlab7n} .
\end{equation}
Note that each $s_i\in\S_n$, $e_i\in\TL_n$, $h_i\in\PP_n$ and $t_j\in\M_n$.
For $1\leq i<j\leq n$, we define
\begin{equation}\label{eq:elements2}
e_{i,j} = e_{j,i} = \custpartn{1,3,4,5,7,8,9,11}{1,3,4,5,7,8,9,11}{\stlines{1/1,3/3,5/5,7/7,9/9,11/11}\uarc48\darc48\udotted13\udotted57\udotted9{11}\ddotted9{11}\ddotted13\ddotted57\vertlab11\vertlab4i\vertlab8j\vertlab{11}n} ,
\end{equation}
noting that each $e_{i,j}\in\B_n$.  We also have $e_i = e_{i,i+1}$ for $1\leq i\leq n-1$.

\section{Diagram monoids}\label{sect:D}

In this section we classify the involutions of our diagram monoids---$\P_n$, $\PP_n$, $\PB_n$, $\M_n$, $\B_n$ and~$\TL_n$---and characterise star-regularity and regular-starity.  
To understand the involutions, we must first understand the automorphisms.  In Section \ref{subsect:AutS} we recall the classification of the automorphisms of~$\P_n$, $\PB_n$ and~$\B_n$ from \cite{Maz2002}, and then provide the classifications for $\PP_n$, $\M_n$ and~$\TL_n$; see Theorems \ref{thm:AutPnBn}, \ref{thm:AutTLn} and \ref{thm:AutPPnMn}.  We then use these automorphisms, combined with a result from \cite{EN2016} to classify the involutions in Section \ref{subsect:InvS}; see Theorem~\ref{thm:invD}.  In Section~\ref{subsect:*regS} we characterise the involutions that give rise to star-regular or regular-star structures; see Theorems~\ref{thm:Pn*}--\ref{thm:PPn*}.  For example, each of our diagram monoids is a regular star-monoid with respect to the canonical involution ${}^*$.  This is the only involution for which $\P_n$, $\PB_n$, $\PP_n$ or $\M_n$ is star-regular, whereas~$\B_n$ is star-regular for some of its other involutions (but regular-starity holds only for ${}^*$).  The situation for $\TL_n$ is more delicate, with exceptional behaviour when $n=3$.
For applications later in the paper, we conclude by describing the centres of our diagram monoids in Section \ref{subsect:ZS}, which in the generic cases are trivial; see Theorem \ref{thm:Z}.

Throughout this section, we make use of the special partitions listed in \eqref{eq:elements} and \eqref{eq:elements2}.

\subsection{Automorphisms}\label{subsect:AutS}

Given any monoid $S$, and any unit $g$ of $S$, there is an \emph{inner} automorphism ${\chi_g=\chi_g^S \in \Aut(S)}$ given by
\[
a\chi_g = g^{-1}ag \qquad\text{for all $a\in S$.}
\]
(We typically drop the superscript ${}^S$, as the monoid should be clear from context.)  The automorphisms of $\P_n$, $\PB_n$ and $\B_n$ were classified by Mazorchuk \cite{Maz2002}:

\begin{thm}[{see \cite[Theorems 2.1, 3.1 and 4.1]{Maz2002}}]\label{thm:AutPnBn}
If $S$ is any of $\P_n$, $\PB_n$ or $\B_n$, then
\[
\Aut(S) = \set{\chi_g}{g\in\S_n}.  \epfreseq
\]
\end{thm}

As far as we are aware, the automorphism groups of the planar monoids $\PP_n$, $\M_n$ and~$\TL_n$ have not been described in the literature, so we turn to this now.

Let $\mu$ be the automorphism of $\PP_n$ corresponding to a reflection in a vertical axis.  For example, with $b = \custpartn{1,...,6}{1,...,6}{
\uarc12
\uarc34
\darc45
\darc56
\darc23
\stline31
\stline55}\in\PP_6$ as in Figure \ref{fig:P6}, we have $b\mu = \custpartn{1,...,6}{1,...,6}{
\uarc56
\uarc34
\darc12
\darc23
\darc45
\stline46
\stline22}$.
So $\mu$ is the restriction to $\PP_n$ of the inner automorphism $\chi_g\in\Aut(\P_n)$ for the unique order-reversing permutation ${g = (1,n)(2,n-1)\cdots}\in\S_n$.  We also denote by $\mu$ the corresponding automorphisms of $\M_n$ and~$\TL_n$, with context making it clear which monoid $\mu$ is acting on.  It is known (see for example \cite[Theorem 1.11(b)]{HR2005}) that $\PP_n$ is generated as a monoid by the set ${\{e_1,\ldots,e_{n-1}\}\cup\{t_1,\ldots,t_n\}}$, where these partitions are as in \eqref{eq:elements}, and we note that
\[
e_i\mu = e_{n-i} \AND t_j\mu = t_{n-j+1} \qquad\text{for all $1\leq i\leq n-1$ and $1\leq j\leq n$.}
\]

We are now ready describe the automorphism groups of $\PP_n$, $\M_n$ and $\TL_n$.  For convenience, we state two separate results, beginning with $\TL_n$.  For both proofs, recall that the $\J$-classes of a monoid $S$ are partially ordered by
\begin{align*}
J_1 \leq J_2 &\iff x \in SyS \text{ for some $x\in J_1$ and $y\in J_2$}\\
&\iff x \in SyS \text{ for all $x\in J_1$ and $y\in J_2$.}
\end{align*}

\newpage

\begin{thm}\label{thm:AutTLn}\leavevmode
\ben
\item \label{AutTLn} For $n\not=3$ we have $\Aut(\TL_n) = \{\id,\mu\}$, but note that $\mu=\id$ for $n\leq2$.
\item \label{AutTL3} We have $\Aut(\TL_3) = \{\id,\mu,\al,\be\}$, where the automorphisms $\al$ and $\be$ are given by:
\[
\al:\begin{cases}
e_1\mt e_1e_2,\\
e_2\mt e_2e_1,
\end{cases}
\AND\ \ 
\be:\begin{cases}
e_1\mt e_2e_1,\\
e_2\mt e_1e_2.
\end{cases}
\]
\een
\end{thm}

\pf
\firstpfitem{\ref{AutTLn}}  This is clear for $n\leq2$, as $\TL_n$ is trivial for $n\leq1$, while $\TL_2$ is isomorphic to the two-element multiplicative semilattice $\{0,1\}$.

We now assume that $n\geq4$.  It is known (see for example \cite[Lemma 2]{BDP2002}) that $\TL_n$ is generated as a monoid by $e_1,\ldots,e_{n-1}$.  Thus, we must show that an arbitrary ${\phi\in\Aut(\TL_n)}$ agrees with~$\id$ or~$\mu$ on these generators.  The $\J$-classes of $\TL_n$ form a chain: ${J_n > J_{n-2} > \cdots}$, where each ${J_r = \set{a\in \TL_n}{\rank(a)=r}}$; see for example \cite[Theorem 3.5]{FL2011}.  It follows that~$\phi$ induces a bijection of each $J_r$, in particular of $J_{n-2}$.  Since $\phi$ also maps idempotents to idempotents, it induces a bijection of $E(J_{n-2})$.  It is well known (see \cite[Lemma~9.8]{EG2017}) that
\[
E(J_{n-2}) = \set{e_i}{1\leq i\leq n-1} \cup \set{e_ie_{i+1},e_{i+1}e_i}{1\leq i\leq n-2}.
\]
The idempotents from the right-hand subset have the form
\[
e_ie_{i+1} = \custpartn{1,3,4,5,6,7,9}{1,3,4,5,6,7,9}{\stlines{1/1,3/3,6/4,7/7,9/9}\uarc45\darc56\udotted13\udotted79\ddotted13\ddotted79\vertlab11\vertlab4i\vertlab9n} 
\AND
e_{i+1}e_i = \custpartn{1,3,4,5,6,7,9}{1,3,4,5,6,7,9}{\stlines{1/1,3/3,4/6,7/7,9/9}\uarc56\darc45\udotted13\udotted79\ddotted13\ddotted79\vertlab11\vertlab4i\vertlab9n} .
\]
Now let $\Ga_n$ be the graph with vertex set $E(J_{n-2})$, and an edge $e-f$ for distinct $e,f\in E(J_{n-2})$ if and only if $e=efe$ and $f=fef$, i.e.~if $e$ and $f$ are mutual inverses.  This graph is shown in Figure \ref{fig:Ga} in the case $n=6$.  Since $\phi$ maps mutual inverses to mutual inverses, it follows that it restricts to an automorphism of $\Ga_n$.  Examining the graph, the idempotents $e_i$ have the following degrees:
\[
\deg(e_i) = \begin{cases}
3 &\text{if $i=1$ or $n-1$}\\
6 &\text{if $2\leq i\leq n-2$}
\end{cases}
\]
All other vertices have degree $4$ or $5$.  Thus, it follows that $\phi$ permutes the sets $\{e_1,e_{n-1}\}$ and $\{e_2,\ldots,e_{n-2}\}$.  Since $\phi$ also preserves the adjacency relation among the $e_i$, it follows that
\[
e_i\phi = e_i \text{ for all $1\leq i\leq n-1$} \OR e_i\phi = e_{n-i} \text{ for all $1\leq i\leq n-1$.}
\]
That is, $\phi=\id$ or $\phi=\mu$.

\pfitem{\ref{AutTL3}}  This can be proved computationally using GAP \cite{Semigroups, GAP}.  It can also be proved directly, using the fact that $\TL_3$ is a $2\times2$ rectangular band with an identity attached.
\epf

\begin{rem}
The graph $\Ga_n$ constructed in the proof of Theorem \ref{thm:AutTLn}\ref{AutTLn} has four automorphisms for $n\geq4$, but only two correspond to automorphisms of $\TL_n$.  The graph $\Ga_3$ has $24$ automorphisms (see Figure \ref{fig:Ga}), but only four correspond to automorphisms of $\TL_3$.

Since $\mu^2=\al^2=\be^2=\id$, we see that $\Aut(\TL_3)$ is the Klein $4$-group.
\end{rem}

\begin{figure}[ht]
\begin{center}\scalebox{0.9}{
\begin{tikzpicture}[scale=2.5]
\tikzstyle{vertex}=[circle,draw=black, fill=white, inner sep = 0.07cm]
\node[vertex] (1) at (1,0) {\Large $1$};
\node[vertex] (2) at (2,0) {\Large $2$};
\node[vertex] (3) at (3,0) {\Large $3$};
\node[vertex] (4) at (4,0) {\Large $4$};
\node[vertex] (5) at (5,0) {\Large $5$};
\node[vertex] (12) at (1.5,.5) {$12$};
\node[vertex] (32) at (2.5,.5) {$32$};
\node[vertex] (34) at (3.5,.5) {$34$};
\node[vertex] (54) at (4.5,.5) {$54$};
\node[vertex] (21) at (1.5,-.5) {$21$};
\node[vertex] (23) at (2.5,-.5) {$23$};
\node[vertex] (43) at (3.5,-.5) {$43$};
\node[vertex] (45) at (4.5,-.5) {$45$};
\draw
(1)--(2)--(3)--(4)--(5)
(12)--(32)--(34)--(54)
(21)--(23)--(43)--(45)
(12)--(21) (23)--(32) (34)--(43) (45)--(54)
(1)--(12)--(2)--(32)--(3)--(34)--(4)--(54)--(5)
(1)--(21)--(2)--(23)--(3)--(43)--(4)--(45)--(5)
;
\begin{scope}[shift = {(5,0)}]
\node[vertex] (1) at (1,0) {\Large $1$};
\node[vertex] (2) at (2,0) {\Large $2$};
\node[vertex] (12) at (1.5,.5) {$12$};
\node[vertex] (21) at (1.5,-.5) {$21$};
\draw
(1)--(2)--(12)--(1)--(21)--(2) (12)--(21)
;
\end{scope}
\end{tikzpicture}
}
\caption{The graphs $\Ga_n$ from the proof of Theorem \ref{thm:AutTLn}, shown in the cases $n=6$ (left) and $n=3$ (right).  A vertex $e_i$ or $e_ie_j$ is denoted $i$ or $ij$, respectively.  }
\label{fig:Ga}
\end{center}
\end{figure}
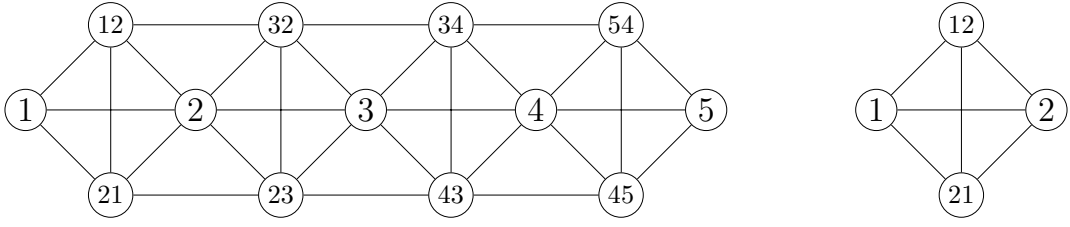

\begin{thm}\label{thm:AutPPnMn}
If $S$ is either of $\PP_n$ or $\M_n$, then $\Aut(S) = \{\id,\mu\}$, but note that $\mu=\id$ for~$n\leq1$.
\end{thm}

\pf
The $S=\PP_n$ case follows immediately from Theorem \ref{thm:AutTLn}, given the isomorphism $\PP_n\cong\TL_{2n}$, so for the rest of the proof we assume that $S=\M_n$.

Again the $\J$-classes of $\M_n$ form a chain: ${J_n>J_{n-1}>\cdots>J_0}$, where this time each ${J_r=\set{a\in\M_n}{\rank(a)=r}}$; see for example \cite[Proposition~2.6]{DEG2017}.  It is known (see for example \cite[Remark 4.14]{DEG2017}) that $\M_n$ is generated as a monoid by the set $J_{n-1}\cup\{e_1,\ldots,e_{n-1}\}$, where again the elements $e_i$ are as in \eqref{eq:elements}.  Thus, we can prove this part by showing that an arbitrary automorphism $\phi\in\Aut(\M_n)$ acts as $\id$ or $\mu$ on this generating set, i.e.~that we have
\begin{equation}\label{eq:AutMn}
a\phi = g^{-1}ag \qquad\text{for all $a\in J_{n-1}\cup\{e_1,\ldots,e_{n-1}\}$,}
\end{equation}
where $g\in\S_n$ is either the identity permutation, or else the unique order-reversing permutation $(1,n)(2,n-1)\cdots$.

Now, $\phi$ induces a bijection of each $\J$-class, and hence in particular of $J_{n-1}$.  It is easy to see that $J_{n-1} = \set{u_{i,j}}{i,j\in\bn}$, where $u_{i,j}$ is the unique element of $\M_n$ with domain $\bn\sm\{i\}$ and codomain $\bn\sm\{j\}$.  For example, when $1\leq i<j\leq n$ we have
\[
u_{i,j} = \custpartn{1,3,4,5,7,8,10}{1,3,4,6,7,8,10}{\stlines{1/1,3/3,5/4,7/6,8/8,10/10}\udotted13\udotted57\udotted8{10}\ddotted8{10}\ddotted13\ddotted46\vertlab11\vertlab4i\vertlab7j\vertlab{10}n}.
\]
Note that $u_{i,i} = t_i$ is the idempotent pictured in \eqref{eq:elements}.  Since $\phi$ also restricts to a bijection of $E(J_{n-1}) = \{t_1,\ldots,t_n\}$, it follows that there is a permutation $g\in\S_n$ such that $t_i\phi=t_{ig}$ for all $i\in\bn$.  Now let $i,j\in\bn$ be arbitrary.  Since $\phi$ maps $J_{n-1}$ into $J_{n-1}$, we have $u_{i,j}\phi=u_{k,l}$ for some $k,l\in\bn$.  But then $u_{k,l} = u_{i,j}\phi = (t_iu_{i,j})\phi = t_{ig}u_{k,l}$, and it follows that $k=ig$; an analogous argument gives $l=jg$.  So in fact
\[
u_{i,j}\phi = u_{ig,jg} \qquad\text{for all $i,j\in\bn$.}
\]
Next note that if $i,j\in\bn$ are such that $|i-j|=1$ (i.e.~if $i$ and $j$ are adjacent in the usual order on $\bn$), then $u_{i,j}^2 = t_it_j$ is an idempotent of rank $n-2$.  On the other hand, if $|i-j|\geq2$, then~$u_{i,j}^2$ is not an idempotent.  For example, if $j\geq i+2$, then
\[
u_{i,j}^2 = \custpartn{1,3,4,5,6,8,9,11}{1,3,4,6,7,8,9,11}{\stlines{1/1,3/3,6/4,8/6,9/9,11/11}\udotted13\udotted68\udotted9{11}\ddotted9{11}\ddotted13\ddotted46\vertlab11\vertlab4i\vertlab8j\vertlab{11}n}.
\]
It follows that for all $i,j\in\bn$,
\begin{align*}
|i-j|=1 &\implies u_{i,j}^2 \text{ is an idempotent}\\
&\implies u_{ig,jg}^2 = u_{i,j}^2\phi\text{ is an idempotent}\\
&\implies |ig-jg|=1.
\end{align*}
(We cannot have $ig=jg$ as $|i-j|=1\implies i\not=j\implies ig\not=jg$.)  This says that $g\in\S_n$ preserves the adjacency relation in the chain $\bn=\{1<\cdots<n\}$, and it follows that $g$ is either the identity or else the unique order-reversing permutation.  In both of these cases, we also have
\[
u_{i,j}\phi = u_{ig,jg} = g^{-1}u_{i,j}g \qquad\text{for all $i,j\in\bn$.}
\]
This completes the proof of \eqref{eq:AutMn} for $a\in J_{n-1}$.

It therefore remains to show that $e_i\phi=g^{-1}e_ig$ for all $1\leq i\leq n-1$.  To do so, fix some such~$i$, and first consider the element
\[
a = t_it_{i+1}e_i = \custpartn{1,3,4,5,6,8}{1,3,4,5,6,8}{\stlines{1/1,3/3,6/6,8/8}\darc45\udotted13\udotted68\ddotted13\ddotted68\vertlab11\vertlab4i\vertlab8n} \in E(J_{n-2}).
\]
Recalling that $ig$ and $(i+1)g$ are adjacent, we have $\{ig,(i+1)g\} = \{k,k+1\}$ for some ${1\leq k\leq n-1}$.  Since
\[
a\phi = (t_it_{i+1}e_i)\phi = t_kt_{k+1}(e_i\phi),
\]
it follows that $a\phi$ contains the blocks $\{k\}$ and $\{k+1\}$.  Since $\phi$ induces a bijection on $E(J_{n-2})$, we have $a\phi\in E(J_{n-2})$.  There are only two Motzkin idempotents of rank $n-2$ containing the blocks $\{k\}$ and $\{k+1\}$, namely
\[
t_kt_{k+1}e_k = \custpartn{1,3,4,5,6,8}{1,3,4,5,6,8}{\stlines{1/1,3/3,6/6,8/8}\darc45\udotted13\udotted68\ddotted13\ddotted68\vertlab11\vertlab4k\vertlab8n}
\AND
t_kt_{k+1} = \custpartn{1,3,4,5,6,8}{1,3,4,5,6,8}{\stlines{1/1,3/3,6/6,8/8}\udotted13\udotted68\ddotted13\ddotted68\vertlab11\vertlab4k\vertlab8n}.
\]
We cannot have $a\phi=t_kt_{k+1}$, since $(t_it_{i+1})\phi=t_kt_{k+1}$, and since $\phi$ is injective, so in fact $a\phi=t_kt_{k+1}e_k$.  All of the above shows that
\[
(t_it_{i+1}e_i)\phi = t_kt_{k+1}e_k.
\]
A symmetrical argument gives $(e_it_it_{i+1})\phi = e_kt_kt_{k+1}$ (for the same $k$).  It then follows that
\[
e_i\phi = (e_it_it_{i+1}\cdot t_it_{i+1}e_i)\phi = e_kt_kt_{k+1} \cdot t_kt_{k+1}e_k = e_k.
\]
But remembering that $\{k,k+1\} = \{ig,(i+1)g\}$, it follows that $e_i\phi=e_k=g^{-1}e_ig$, as required.
\epf

\subsection{Involutions}\label{subsect:InvS}

We are now ready to classify the involutions of our diagram monoids.  As we will show in Theorem~\ref{thm:invD}, there are many such involutions of $\P_n$, $\PB_n$ and $\B_n$ (one for each self-inverse permutation of $\bn$), but only two for the planar monoids $\PP_n$, $\M_n$ and $\TL_n$.

In what follows, we write $\io$ for the canonical involution of $\P_n$, given by $a\io = a^*$.\footnote{Because we will be composing various involutions and (anti)automorphisms, it is more convenient here to use `functional' notation for our involutions, rather than `index' notation such as ${}^*$ (although the latter will be used at times).}  For a permutation $g\in\S_n$, we denote by $\io_g = \io\chi_g$ the anti-automorphism of $\P_n$ given by
\[
a\io_g = g^{-1}a^*g = g^*a^*g \qquad\text{for $a\in\P_n$.}
\]
In particular we have $\io_1 = \io$.  It is easy to check that
\begin{equation}\label{eq:iogh}
\chi_g\chi_h = \chi_{gh} = \io_g\io_h \AND \io_g\chi_h = \io_{gh} = \chi_g\io_h \qquad\text{for any $g,h\in\S_n$.}
\end{equation}
We also write $\io$ and $\io_g$ for the restrictions of the above maps to the monoids $\PB_n$ and $\B_n$, with context again making it clear which monoid these are acting on.

The planar partition monoid $\PP_n$ is not closed under conjugation by arbitrary permutations, and hence under the maps $\io_g\in\Aut^-(\P_n)$ in general.  We observed above, however, that the automorphism $\mu\in\Aut(\PP_n)$ is the restriction to $\PP_n$ of $\chi_g\in\Aut(\P_n)$ for the unique order-reversing permutation $g = (1,n)(2,n-1)\cdots\in\S_n$.  We will denote by $\rho = \io\mu \in \Aut^-(\PP_n)$ the composition of~$\io$ with~$\mu$.  This is defined by the rule
\[
a\rho = a^*\mu = (a\mu)^* \qquad\text{for $a\in\PP_n$,}
\]
and corresponds to a $180^\circ$ rotation.  For example, with $b = \custpartn{1,...,6}{1,...,6}{
\uarc12
\uarc34
\darc45
\darc56
\darc23
\stline31
\stline55}\in\PP_6$ as in Figure \ref{fig:P6}, we have $b\rho = \custpartn{1,...,6}{1,...,6}{
\uarc12
\uarc23
\uarc45
\darc34
\darc56
\stline22
\stline64}$.
We will also denote the corresponding $180^\circ$ rotation maps on $\M_n$ and $\TL_n$ by $\rho$.  The involution $\rho$ has been considered for $\TL_n$ in \cite{AV2020,ADV2012_2}.

\begin{thm}\label{thm:invD}\leavevmode
\ben
\item \label{invD1} If $S$ is any of $\P_n$, $\PB_n$ or $\B_n$, then the involutions of $S$ are the maps~$\io_g$, for $g\in\S_n$ of order $|g|\leq2$.
\item \label{invD2} If $S$ is any of $\PP_n$, $\M_n$ or $\TL_n$, then the involutions of $S$ are the maps  $\io$ and $\rho$.
\een
\end{thm}

\pf
\firstpfitem{\ref{invD1}}  We first note that $\io$ commutes with every automorphism of $S$.  Indeed, by Theorem~\ref{thm:AutPnBn}, any such automorphism has the form $\chi_g$ for some $g\in\S_n$, and it follows from \eqref{eq:iogh} that this commutes with $\io_1 = \io$.  It then follows from Lemma \ref{lem:EN} that the involutions of $S$ are precisely the maps obtained by composing $\io$ with an automorphism of order $\leq2$.  The result now follows from the easily-checked fact that $|\chi_g|=|g|$, with a single exception; when $S=\B_2$, we have $\chi_{(1,2)} = \chi_{\id_{{\bf2}}}$, which has order $1$.  

\pfitem{\ref{invD2}}  It is easy to check that $\io$ commutes with the automorphisms $\id$ and $\mu$, and so $\io(=\io\id)$ and $\rho(=\io\mu)$ are involutions of $S$.  If $S\not=\TL_3$, then $\Aut(S)=\{\id,\mu\}$ by Theorems \ref{thm:AutTLn} and~\ref{thm:AutPPnMn}, so these are all the involutions of $S$.  We now consider $S=\TL_3$, in which case $\Aut(\TL_3)=\{\id,\mu,\al,\be\}$ in the notation of Theorem \ref{thm:AutTLn}\ref{AutTL3}.  Here one can check that $\io\al\io=\be\not=\al^{-1}$, so neither~$\io\al$ nor~$\io\be$ is an involution of~$\TL_3$.
\epf

The classification of involutions allows us to deduce the following result concerning the canonical float-counting twisting $\Phi$, which will be used in Section \ref{sect:T}.

\begin{cor}\label{cor:*symD}
If $S$ is any of $\P_n$, $\PP_n$, $\PB_n$, $\M_n$, $\B_n$ or $\TL_n$, and if $\tau$ is any involution of~$S$, then for all $a,b\in S$ we have $\Phi(a,b) = \Phi(b\tau,a\tau)$.
\end{cor}

\pf
By Theorem \ref{thm:invD}, the involution $\tau$ is given by $a\tau = g^{-1}a^*g$ for some fixed $g\in\S_n$.  Recalling that $\Phi$ is unital, it now follows from \ref{U5} and Lemma \ref{lem:*symPn} that
\[
\Phi(b\tau,a\tau) = \Phi(g^{-1}b^*g,g^{-1}a^*g) = \Phi(b^*,a^*) = \Phi(a,b).  \qedhere
\]
\epf

We also record the following. 

\begin{cor}\label{cor:rankatau}
If $S$ is any of $\P_n$, $\PP_n$, $\PB_n$, $\M_n$, $\B_n$ or $\TL_n$, and if $\tau$ is any involution of~$S$, then for all $a\in S$ we have $\rank(a\tau)=\rank(a)$.
\end{cor}

\pf
This again follows from the fact that $a\tau = g^{-1}a^*g$ for some $g\in\S_n$, so $a\tau \J a^* \J a$.
\epf

\subsection{Star-regularity and regular-starity}\label{subsect:*regS}

Now that we have classified the involutions of our diagram monoids (Theorem \ref{thm:invD}), we characterise those that give rise to star-regular or regular star-structures.
The situations for the various monoids are somewhat different, so it is convenient to give four separate statements in Theorems~\ref{thm:Pn*}--\ref{thm:PPn*}, grouping together the monoids with similar behaviour.

\begin{thm}\label{thm:Pn*}
If $S=\P_n$ or $\PB_n$, and if $g\in\S_n$ has order $\leq2$, then 
\[
\text{$S$ is an $\io_g$-regular monoid} \IFf \text{$S$ is a regular $\io_g$-monoid} \IFf g=1.
\]
\end{thm}

\pf
The backward implications have all been discussed previously, so suppose $S$ is $\io_g$-regular; we must show that $g=1$.  Suppose to the contrary that $ig=j$ for distinct $i,j\in\bn$.  We see then that $t_i\io_g = t_j$ (where the $t_i$ are as in \eqref{eq:elements}).  If $S$ was $\io_g$-regular, then the $\H$-class of~$t_j$ would contain an inverse~$a$ of $t_i$.  By \eqref{eq:RJ} every such element $a$ contains the upper block~$\{j\}$, and so too therefore does $t_iat_i$, which gives $t_iat_i\not=t_i$.
\epf

\begin{thm}\label{thm:Bn*}
If $g\in\S_n$ has order $\leq2$, then 
\ben
\item \label{Bn*1} $\B_n$ is an $\io_g$-regular monoid$\IFf g=1$ or $g$ is a transposition,
\item \label{Bn*2} $\B_n$ is a regular $\io_g$-monoid$\IFf g=1$ \emph{(}or $n=2$ and $g=(1,2)$, where $\io_{(1,2)}=\io=\id$\emph{)}.
\een
\end{thm}

\pf
\firstpfitem{\ref{Bn*1}}  First suppose $\B_n$ is $\io_g$-regular.  Since $\io_g$ is an involution, $g$ is a product of disjoint transpositions.  Aiming for a contradiction, suppose this product involves at least two transpositions, say $(x,y)$ and $(u,v)$.  Then $e_{x,u}\io_g = e_{y,v}$ is not $\H$-related to any inverse of $e_{x,u}$ (where the $e_{i,j}$ are as in \eqref{eq:elements2}).  Indeed, any element~$a$ in the $\H$-class of $e_{y,v}$ contains the block $\{y,v\}$, and so too therefore does $e_{x,u}ae_{x,u}$, which gives $e_{x,u}ae_{x,u}\not=e_{x,u}$.

Conversely, since we know that $\B_n$ is a regular $\io$-monoid, it suffices to show that $\B_n$ is $\io_g$-regular when $g=(x,y)$ for distinct $x,y\in\bn$.  As explained at the end of Section \ref{subsect:*}, we can do this by showing that every $\R$-class~$R$ of $\B_n$ contains an $\io_g$-projection, i.e.~an element $e$ such that $e^2=e=e\io_g$.  Let $p\in R$ be the unique $\io$-projection in $R$, which exists because of $\io$-regularity.  It is well known (see \cite[Lemma 4]{EF2012}) that $p$ has the form
\[
p =  \begin{partn}{6}i_1&\cdots&i_r&u_1,v_1&\cdots&u_t,v_t \\ \hhline{~|~|~|-|-|-} i_1&\cdots&i_r&u_1,v_1&\cdots&u_t,v_t\end{partn}.
\]
(This notation indicates that $p$ has transversals $\{i_1,i_1'\},\ldots,\{i_r,i_r'\}$, upper blocks $\{u_1,v_1\},\ldots,\{u_t,v_t\}$, and lower blocks $\{u_1',v_1'\},\ldots,\{u_t',v_t'\}$.)
If $x,y\in\{i_1,\ldots,i_r\}$ or if $\{x,y\} = \{u_s,v_s\}$ for some $s$, then $p\io_g = p$, and $p$ is itself an $\io_g$-projection (in $R$).  Up to symmetry, there are only two other possibilities:
\bena
\item \label{Bn*a} $x=i_1$ and $y=u_1$, or
\item \label{Bn*b} $x=u_1$ and $y=u_2$.
\een
In each case, we claim that $e = pg\in R$ is an $\io_g$-projection.  Since $g$ is a unit, certainly $e\in R$, and we also have $e\io_g = g(pg)^*g = g(gp)g = pg = e$, since $p^*=p$, $g^*=g$ and $g^2=1$.  Finally, it is easy to check that $e$ is an idempotent, given that in cases \ref{Bn*a} and \ref{Bn*b} we respectively have
\[
e = \begin{partn}{8}
i_1&i_2&\cdots&i_r&u_1,v_1&u_2,v_2&\cdots&u_t,v_t \\ \hhline{~|~|~|~|-|-|-|-} 
u_1&i_2&\cdots&i_r&i_1,v_1&u_2,v_2&\cdots&u_t,v_t
\end{partn}
\OR
e = \begin{partn}{8}
i_1&\cdots&i_r&u_1,v_1&u_2,v_2&u_3,v_3&\cdots&u_t,v_t \\ \hhline{~|~|~|-|-|-|-|-} 
i_1&\cdots&i_r&u_2,v_1&u_1,v_2&u_3,v_3&\cdots&u_t,v_t
\end{partn}.
\]

\pfitem{\ref{Bn*2}}  Since we know that $\B_n$ is a regular $\io$-monoid, and since $\io_{(1,2)} = \io$ when $n=2$, it remains by symmetry to show that $\B_n$ ($n\geq3$) is not a regular $\io_g$-monoid for the transposition $g=(1,2)$, i.e.~that there exists $a\in\B_n$ for which~$a\io_g$ is not an inverse of $a$.  But this is easy to check for $a = (1,2,3) \in \S_n \sub \B_n$.
\epf

For the planar monoids $\PP_n$, $\M_n$ and $\TL_n$, there are only two involutions to consider, $\io$ and $\rho$, but we recall that $\rho=\io=\id$ for sufficiently small $n$.

\newpage

\begin{thm}\label{thm:TLn*}\leavevmode
\ben
\item \label{TLn*1} $\TL_n$ is a regular $\io$-monoid for any $n\geq0$.
\item \label{TLn*2} $\TL_3$ is a regular $\rho$-monoid.
\item \label{TLn*3} $\TL_n$ is not $\rho$-regular for $n\geq4$.
\een
\end{thm}

\pf
\firstpfitem{\ref{TLn*1}}  We have already noted this.

\pfitem{\ref{TLn*2}}  This follows from the fact that $\rho$ maps the two $\J$-classes $J_1$ and $J_3$ into themselves (in the notation of the proof of Theorem \ref{thm:AutTLn}), and that these are both rectangular bands.  (Any pair of elements in a rectangular band are mutual inverses.)

\pfitem{\ref{TLn*3}}  Since $\TL_n$ is $\H$-trivial, it is enough to note that $e_1$ are $e_1\rho = e_{n-1}$ are not mutual inverses when~$n\geq4$.
\epf

\begin{thm}\label{thm:PPn*}
If $S=\PP_n$ or $\M_n$, then
\ben
\item \label{PPn*1} $S$ is a regular $\io$-monoid for any $n\geq0$,
\item \label{PPn*2} $S$ is not $\rho$-regular for $n\geq2$.
\een
\end{thm}

\pf
We have already discussed \ref{PPn*1}.  For \ref{PPn*2}, it is again sufficient to note that $t_1$ is not an inverse of $t_1\rho=t_n$ for $n\geq2$.
\epf

\subsection{The centre}\label{subsect:ZS}

It will transpire that in order to apply the above results to twisted diagram monoids in Section~\ref{sect:T}, we will need to understand the centres of our (un-twisted) diagram monoids, and we turn to this task now.  In particular we will show in Theorem \ref{thm:Z} that they all have trivial centre for sufficiently large $n$.  This is possibly known, but we are unaware of an explicit statement or proof in the literature, and we believe that the result is of independent interest.

So that our proof of Theorem \ref{thm:Z} works for all diagram monoids at once, we first prove a result about the centraliser of the Temperley--Lieb and Motzkin monoids in the partition monoid.  For the proof we need the following simple fact.

\begin{lemma}\label{lem:ijkl}
If $n\geq3$, and if $i,j\in\bn$ are distinct, then there exist $k\in\{i,j\}$ and $l\in\bn\sm\{i,j\}$ such that ${k\not\equiv l\pmod 2}$.  \epfres
\end{lemma}

The following statement and proof uses the partitions from \eqref{eq:elements} and \eqref{eq:elements2}.  The proof also uses the (co)domain and (co)kernel parameters from \eqref{eq:dom}.

\begin{prop}\label{prop:C}\leavevmode
\ben
\item \label{C1} If $n\geq3$, then $C_{\P_n}(\TL_n) = \{1\}$,
\item \label{C2} $C_{\P_2}(\TL_2) = \{1,s_1,e_1,h_1\}$,
\item \label{C3} If $n\geq2$, then $C_{\P_n}(\M_n) = \{1\}$,
\een
\end{prop}

\pf
\firstpfitem{\ref{C1}}  Let $a\in C_{\P_n}(\TL_n)$, where $n\geq3$.  
We first claim that
\[
\ker(a)=\De=\coker(a),
\]
where here $\De$ denotes the trivial equivalence on $\bn$.  To prove this, it suffices by symmetry to show that $\ker(a)=\De$.  To do so, suppose to the contrary that ${(i,j)\in\ker(a)}$ for distinct $i,j\in\bn$.  Let $k\in\{i,j\}$ and $l\in\bn\sm\{i,j\}$ be as in Lemma~\ref{lem:ijkl}, and choose any $b\in\TL_n$ containing the upper block $\{k,l\}$.  (It is easy to see that such $b$ exists, given $k\not\equiv l\pmod2$.)  Then $\{k,l\}$ is also a block of $ba$, which implies that $(i,j)\not\in\ker(ba)$.  But $(i,j)\in\ker(a)\sub\ker(ab)$, so it follows that $ab\not=ba$, a contradiction.

Next we claim that
\[
\dom(a)=\bn=\codom(a).
\]
Again it suffices to prove the statement about $\dom(a)$.  For this, suppose to the contrary that $\dom(a)\not=\bn$.  By the previous claim we have $\ker(a)=\De$, so $a$ has an upper singleton block,~$\{i\}$ say.  Choose any $j\in\bn\sm\{i\}$ of opposite parity to $i$, and any $b\in\TL_n$ containing the upper block $\{i,j\}$.  Then $\{i\}$ and $\{i,j\}$ are blocks of $ab$ and $ba$, respectively, so $ab\not=ba$.

It follows from the two claims that $a\in\S_n$ is a permutation.  For any $1\leq i\leq n-1$, the element $a$ commutes with $e_i\in\TL_n$, and so $e_{i,i+1} = e_i = a^{-1}e_ia = e_{ia,(i+1)a}$, so it follows that $\{i,i+1\}a = \{i,i+1\}$.  Since this holds for all $1\leq i\leq n-1$ (and since $n\geq3$), it quickly follows that $a=1$ is the identity permutation.

\pfitem{\ref{C2}}  It is easy to see that the four stated elements commute with the two elements of ${\TL_2 = \{1,e_1\}}$.  Conversely, suppose $a\in\P_2$ commutes with  $e_1\in\TL_2$.  Note that $e_1a$ contains the upper block $\{1,2\}$, while $ae_1$ contains the lower block $\{1',2'\}$; since $ae_1=e_1a$ therefore contains both these blocks, it follows that $ae_1=e_1a=e_1$.  If $a$ had a singleton block, then so too would one of $ae_1$ or $e_1a$, contradicting $ae_1=e_1a=e_1$.  So $a$ has no singlteton blocks, which leaves precisely the four stated elements.

\pfitem{\ref{C3}}  Since $\TL_n\sub\M_n$, we have $C_{\P_n}(\M_n) \sub C_{\P_n}(\TL_n)$, so the $n\geq3$ case follows immediately from part \ref{C1}.  For $n=2$, it is easy to check that the three non-identity elements of $C_{\P_2}(\TL_2)$ listed in part \ref{C2} do not commute with $t_1\in\M_2$.
\epf

\begin{thm}\label{thm:Z}\leavevmode
\ben
\item \label{Z1} For $n\geq2$ we have $Z(\P_n) = Z(\PP_n) = Z(\PB_n) = Z(\M_n) = \{1\}$.
\item \label{Z2} For $n\geq3$ we have $Z(\B_n) = Z(\TL_n) = \{1\}$.
\een
\end{thm}

\pf
Part \ref{Z2} and the $n\geq3$ case of part \ref{Z1} follow from Proposition \ref{prop:C}\ref{C1}.  The $n=2$ case of part \ref{Z1} follows from Proposition \ref{prop:C}\ref{C2}.
\epf

\begin{rem}
Each of $\P_n$, $\PP_n$, $\PB_n$, $\M_n$ ($n\leq1$), and $\B_n$,~$\TL_n$~($n\leq2$) is commutative, so for each of these we have $Z(S)=S$.
\end{rem}

\section{Involutions of twisted products}\label{sect:*}

We have now classified the involutions of our diagram monoids---$\P_n$, $\PP_n$, $\PB_n$, $\M_n$, $\B_n$ and~$\TL_n$---and characterised those giving rise to star-regular or regular-star structures.  The corresponding results for \emph{twisted} diagram monoids are harder to obtain, and we must first prove a number of results for general twisted products.  Most of these results require the assumption of unitality or semi-unitality.  In the current section we classify involutions, and in the next we treat star-regularity and regular-starity.  

In Section \ref{subsect:inv1} we identify a special kind of involution of a twisted product that we call \emph{pure}, and we classify these in Theorem \ref{thm:inv1} in the semi-unital case.  Arbitrary involutions are treated in Section~\ref{subsect:inv2}, under the stronger assumption of unitality, with the classification given in Theorem~\ref{thm:inv2}.  
We then derive some useful consequences of the latter theorem in Section \ref{subsect:consequences}, and give some examples of \emph{non-pure} involutions of a specific unital (indeed, tight) twisted product in Section \ref{subsect:examples}.
Along the way, we also consider an important case in which $(i,a)^\oa=(i,a^*)$ determines an involution of a twisted product (see Corollary \ref{cor:*sym}).

\subsection{Pure involutions}\label{subsect:inv1}

Consider an involution $(i,a)\mt(i,a)^\oa$ of a twisted product $M\times_\Phi^qS$.  We say ${}^\oa$ is \emph{pure} if the second coordinate of $(i,a)^\oa$ depends only on $a$.  This means that we have a map $S\to S:a\mt a^*$ such that for all $i\in M$ and $a\in S$,
\[
(i,a)^\oa = (j,a^*) \qquad\text{for some $j = j(i,a)\in M$.}
\]
In this case we say ${}^\oa$ \emph{extends} ${}^*$.  The next simple result will be used frequently without explicit reference.  Note that it does not place any restrictions on the twisting/product, such as \mbox{(semi-)unitality}.

\begin{lemma}\label{lem:pure}
If ${}^\oa$ is a pure involution of a twisted product $M\times_\Phi^qS$ extending ${}^*$, then~${}^*$ is an involution of $S$.
\end{lemma}

\pf
Let $a,b\in S$.  By definition there exist $i,j\in M$ such that
\[
(0,a) = (0,a)^{\oa\oa} = (i,a^*)^\oa = (j,a^{**}),
\]
and so $a=a^{**}$.  Similarly, there exist $k,l,m\in M$ such that
\begin{align*}
(k,(ab)^*) = (\Phi(a,b)q,ab)^\oa = ((0,a)(0,b))^\oa = (0,b)^\oa(0,a)^\oa &= (l,b^*)(m,a^*) \\ &= (l+m+\Phi(b^*,a^*)q,b^*a^*),
\end{align*}
which gives $(ab)^* = b^*a^*$.
\epf

The next result contains the classification of the pure involutions of semi-unital products.  In the statement, note that the forward direction does not assume semi-unitality, in which case the product $M\times_\Phi^qS$ might not even be a monoid.

\begin{thm}\label{thm:inv1}
Let $T = M\times_\Phi^qS$ be a twisted product.  Suppose we are given:
\ben
\item \label{1inv1} involutions $S\to S:a\mt a^*$ and $M\to M:i\mt i^\circ$, and
\item \label{1inv2} a map $S\to M:a\mt k_a$ such that for all $a,b\in S$:
\bena
\item \label{inv2a} $k_a^\circ+k_{a^*} = 0$,
\item \label{inv2b} $k_a+k_b + \Phi(b^*,a^*)q = k_{ab} + \Phi(a,b)q^\circ$. 
\een
\een
Then $(i,a)^\oa = (i^\circ+k_a,a^*)$ determines a pure involution of $T$.  Moreover, if $T$ is semi-unital, then any pure involution of~$T$ arises in this way.
\end{thm}

\pf
First suppose we are given the mappings in \ref{1inv1} and \ref{1inv2}.  Also fix $i,j\in M$ and $a,b\in S$.  Then
\begin{align*}
(i,a)^{\oa\oa} = (i^\circ+k_a,a^*)^\oa &= ((i^\circ+k_a)^\circ+k_{a^*},a^{**}) \\
&= (i^{\circ\circ}+k_a^\circ+k_{a^*},a) &&\hspace{-1cm}\text{since ${}^\circ$ and ${}^*$ are involutions}\\
&= (i+0,a) = (i,a) &&\hspace{-1cm}\text{by \ref{inv2a}, and since ${}^\circ$ is an involution,}
\intertext{and}
((i,a)(j,b))^\oa &= (i+j+\Phi(a,b)q,ab)^\oa \\
&= ((i+j+\Phi(a,b)q)^\circ+k_{ab},(ab)^*) \\
&= (i^\circ+j^\circ+\Phi(a,b)q^\circ+k_{ab},b^*a^*) &&\text{since ${}^\circ$ and ${}^*$ are involutions}\\
&= (i^\circ+j^\circ+k_a+k_b + \Phi(b^*,a^*)q,b^*a^*) &&\text{by \ref{inv2b}}\\
&= (j^\circ+k_b,b^*)(i^\circ+k_a,a^*) \\
&= (j,b)^\oa(i,a)^\oa.
\end{align*}
This all shows that $(i,a)^\oa = (i^\circ+k_a,a^*)$ does indeed determine an involution of $T$, and it is pure by definition.

For the rest of the proof we assume that $T$ is semi-unital.  In what follows, we use \ref{U1'} extensively without explicit reference.  We also fix a pure involution ${T\to T:(i,a)\mt(i,a)^\oa}$, extending an involution $S\to S:a\mt a^*$.  For any $i\in M$ and $a\in S$, we have
\[
(0,a)^\oa = (k_a,a^*)  \AND (i,1)^\oa = (i^\circ,1) \qquad\text{for some $k_a,i^\circ\in M$.}
\]
(For the second of these, note that $1^*=1$ since ${}^*$ is an involution.)  We then note that
\[
(i,a)^\oa = ((0,a)(i,1))^\oa = (i,1)^\oa(0,a)^\oa = (i^\circ,1)(k_a,a^*) = (i^\circ+k_a,a^*) \qquad\text{for $i\in M$ and $a\in S$.}
\]
For $i,j\in M$ we have
\[
(i,1) = (i,1)^{\oa\oa} = (i^\circ,1)^\oa = (i^{\circ\circ},1),
\]
and
\[
((i+j)^\circ,1) = (i+j,1)^\oa = ((j,1)(i,1))^\oa = (i,1)^\oa(j,1)^\oa = (i^\circ,1)(j^\circ,1) = (i^\circ+j^\circ,1),
\]
which together tell us that $i^{\circ\circ}=i$ and $(i+j)^\circ = i^\circ + j^\circ$, i.e.~that ${}^\circ$ is an involution of $M$.  For $a,b\in S$ we have
\[
(0,a) = (0,a)^{\oa\oa} = (k_a,a^*)^\oa = (k_a^\circ + k_{a^*},a),
\]
and
\begin{align*}
(\Phi(a,b)q^\circ+k_{ab},b^*a^*) = (\Phi(a,b)q,ab)^\oa = ((0,a)(0,b))^\oa = (0,b)^\oa&(0,a)^\oa = (k_b,b^*)(k_a,a^*) \\
&= (k_a+k_b+\Phi(b^*,a^*)q,b^*a^*),
\end{align*}
showing that \ref{inv2a} and \ref{inv2b} hold.
\epf

An important and simple special case is as follows.  Given a fixed involution ${}^*$ of $S$, we say that a twisting $\Phi$ is \emph{$*$-symmetric} if $\Phi(a,b) = \Phi(b^*,a^*)$ for all $a,b\in S$.  For example, Lemma \ref{lem:*symPn} says that the canonical float-counting twisting of $\P_n$ is $*$-symmetric for the canonical involution~${}^*$.  Corollary~\ref{cor:*symD} extended this to any involution of any of our diagram monoids.  Again, the next statement does not place any restriction on the twisting/product.

\begin{cor}\label{cor:*sym}
If $T = M\times_\Phi^qS$ is a twisted product, and if $\Phi$ is $*$-symmetric with respect to an involution ${}^*$ of $S$, then ${(i,a)^\oa = (i,a^*)}$ determines a pure involution of $T$.
\end{cor}

\pf
It is routine to check that the conditions of Theorem \ref{thm:inv1} are satisfied by ${}^*$ and the maps
\[
M\to M:i\mt i^\circ=i \AND S\to S:a\mt k_a = 0.
\]
Applying that theorem yields the involution ${(i,a)^\oa = (i^\circ+k_a,a^*) = (i,a^*)}$.
\epf

\subsection{General involutions}\label{subsect:inv2}

The situation of general involutions of $T = M\times_\Phi^qS$ (which need not be pure) is more complicated, and we must strengthen our underlying assumption to unitality of $T$.  We note that the properties listed in the following result contain certain redundancies, but are stated as they are for clarity.  For example, given~\ref{inv9}, we only need the first equalities in~\ref{inv3} and~\ref{inv5}; since $M$ is commutative, we only need the first equalities in \ref{inv1} and \ref{inv2}.

\newpage

\begin{thm}\label{thm:inv2}
Let $T = M\times_\Phi^qS$ be a unital twisted product.  Suppose we are given maps 
\[
S\to S:a\mt a^* \COMMA
M\to M:i\mt i^\circ \COMMA
S\to M:a\mt k_a \AND
M\to S:i\mt c_i ,
\]
such that for all $i,j\in M$ and $a,b\in S$:
\begin{enumerate}[label=\textup{\textsf{(I\arabic*)}},leftmargin=10mm]
\item 
\label{inv1} 
$c_ic_j=c_{i+j}=c_jc_i$, 
\item 
\label{inv2} 
$(i+j)^\circ = i^\circ+j^\circ+\Phi(c_i,c_j)q = i^\circ+j^\circ+\Phi(c_j,c_i)q$, 
\item 
\label{inv3} 
$c_i^*c_{i^\circ} = 1 = c_{i^\circ}c_i^*$,
\item 
\label{inv4} 
$i^{\circ\circ}+k_{c_i} = i$,
\item
\label{inv5} 
$a^{**}c_{k_a} = a =c_{k_a}a^{**}$,
\item 
\label{inv6} 
$k_a^\circ+k_{a^*}+\Phi(c_{k_a},a^{**})q = 0 = k_a^\circ+k_{a^*}+\Phi(a^{**},c_{k_a})q$,
\item 
\label{inv7} 
$(ab)^*c_{\Phi(a,b)q} = b^*a^*$,
\item
\label{inv8} 
$k_a+k_b+\Phi(b^*,a^*)q = k_{ab} + (\Phi(a,b)q)^\circ + \Phi((ab)^*,c_{\Phi(a,b)q})q$,
\item 
\label{inv9} 
$ac_i = c_ia$,
\item
\label{inv10} 
$i^\circ+k_a+ \Phi(a^*,c_i)q = i^\circ+k_a+\Phi(c_i,a^*)q$.
\een
Then $(i,a)^\oa = (i^\circ+k_a+\Phi(a^*,c_i)q,a^*c_i)$ determines an involution of $T$.  Moreover, any involution of $T$ arises in this way.
\end{thm}

To keep the proof manageable, we break it into two separate results, Lemmas \ref{lem:inv21} and \ref{lem:inv22}.  For their statements and proofs we keep the notation of the theorem, including the assumption of unitality.

\begin{lemma}\label{lem:inv21}
If the maps
\[
S\to S:a\mt a^* \COMMA
M\to M:i\mt i^\circ \COMMA
S\to M:a\mt k_a \AND
M\to S:i\mt c_i 
\]
satisfy conditions \ref{inv1}--\ref{inv10}, then $(i,a)^\oa = (i^\circ+k_a+\Phi(a^*,c_i)q,a^*c_i)$ determines an involution of $T$.
\end{lemma}

\pf
We must show that
\begin{equation}\label{eq:inv}
(i,a)^{\oa\oa} = (i,a)
\AND
((i,a)(j,b))^\oa = (j,b)^\oa(i,a)^\oa
\qquad\text{for all $i,j\in M$ and $a,b\in S$.}
\end{equation}
Before we do so, we first establish a consequence of \ref{inv1} and \ref{inv2}, namely that
\begin{equation}\label{eq:ijko}
(i+j+k)^\circ = i^\circ+j^\circ+k^\circ+\Phi(c_i,c_j,c_k)q \qquad\text{for $i,j,k\in M$.}
\end{equation}
Indeed, we have
\begin{align*}
((i+j)+k)^\circ &= (i+j)^\circ + k^\circ + \Phi(c_{i+j},c_k)q &&\text{by \ref{inv2}}\\
&= i^\circ+j^\circ+\Phi(c_i,c_j)q + k^\circ + \Phi(c_ic_j,c_k)q &&\text{by \ref{inv2} and \ref{inv1}}\\
&= i^\circ+j^\circ+ k^\circ +(\Phi(c_i,c_j) + \Phi(c_ic_j,c_k))q \\
&= i^\circ+j^\circ+ k^\circ +\Phi(c_i,c_j,c_k)q .
\end{align*}

We now work towards verifying the first equality in \eqref{eq:inv}.  To do so, fix $i\in M$ and $a\in S$, and also let
\[
(j,b) = (i,a)^\oa = (i^\circ+k_a+\Phi(a^*,c_i)q,a^*c_i).
\]
Then
\[
(i,a)^{\oa\oa} = (j,b)^\oa = (j^\circ+k_b+\Phi(b^*,c_j)q,b^*c_j),
\]
and we need to show that this is equal to $(i,a)$, i.e.~that
\begin{equation}\label{eq:iacaca}
j^\circ+k_b+\Phi(b^*,c_j)q = i \AND b^*c_j = a.
\end{equation}
For the latter we have
\begin{align*}
b^*c_j &= (a^*c_i)^*c_{\Phi(a^*,c_i)q}c_{i^\circ}c_{k_a} &&\text{by \ref{inv1}}\\
&= c_i^*a^{**}c_{i^\circ}c_{k_a} &&\text{by \ref{inv7}}\\
&= c_i^*c_{i^\circ}a^{**}c_{k_a} &&\text{by \ref{inv9}}\\
&= 1\cdot a = a &&\text{by \ref{inv3} and \ref{inv5}.}
\end{align*}
For the former we need to work a little harder.  First we note that
\begin{align}
\nonumber j^\circ &= (i^\circ+\Phi(a^*,c_i)q+k_a)^\circ \\
\nonumber &= i^{\circ\circ}+(\Phi(a^*,c_i)q)^\circ +k_a^\circ+ \Phi(c_{i^\circ},c_{\Phi(a^*,c_i)q},c_{k_a})q &&\text{by \eqref{eq:ijko}}\\
\label{eq:jo} &= i^{\circ\circ}+k_a^\circ+(\Phi(a^*,c_i)q)^\circ + \Phi(c_{\Phi(a^*,c_i)q},c_{k_a})q &&\text{by \ref{U3}, as $c_{i^\circ}$ is a unit by \ref{inv3}.}
\end{align}
We then note that $c_j = c_{\Phi(a^*,c_i)q}c_{k_a}c_{i^\circ}$ by \ref{inv1}.  Since $c_{i^\circ}$ is a unit, it then follows from~\ref{U2} that
\begin{equation}\label{eq:PhiB*cj}
\Phi(b^*,c_j) = \Phi((a^*c_i)^*,c_{\Phi(a^*,c_i)q}c_{k_a}c_{i^\circ}) = \Phi((a^*c_i)^*,c_{\Phi(a^*,c_i)q}c_{k_a}). 
\end{equation}
Combining \eqref{eq:jo} and \eqref{eq:PhiB*cj}, it follows that
\begin{align}
\nonumber j^\circ+k_b+\Phi(b^*,c_j)q = i^{\circ\circ}+k_a^\circ+ k_b &+ (\Phi(a^*,c_i)q)^\circ  \\
\label{eq:jokb} & + (\Phi(c_{\Phi(a^*,c_i)q},c_{k_a})
+ \Phi((a^*c_i)^*,c_{\Phi(a^*,c_i)q}c_{k_a}))q.
\end{align}
Looking at the bracketed sum of $\Phi$ terms in \eqref{eq:jokb}, we note that
\begin{align*}
\Phi(c_{\Phi(a^*,c_i)q},c_{k_a}) + {} &\Phi((a^*c_i)^*,c_{\Phi(a^*,c_i)q}c_{k_a}) \\ 
&= \Phi((a^*c_i)^*,c_{\Phi(a^*,c_i)q}c_{k_a}) + \Phi(c_{\Phi(a^*,c_i)q},c_{k_a}) \\
&= \Phi((a^*c_i)^*,c_{\Phi(a^*,c_i)q}) + \Phi((a^*c_i)^*c_{\Phi(a^*,c_i)q},c_{k_a}) &&\text{by \ref{T1}}\\
&= \Phi((a^*c_i)^*,c_{\Phi(a^*,c_i)q}) + \Phi(c_i^*a^{**},c_{k_a}) &&\text{by \ref{inv7}.}
\end{align*}
Substituting this into \eqref{eq:jokb}, we then continue
\begin{align*}
j^\circ+k_b+\Phi(b^*,c_j)q &= i^{\circ\circ} +k_a^\circ+ k_b +(\Phi(a^*,c_i)q)^\circ + \Phi((a^*c_i)^*,c_{\Phi(a^*,c_i)q})q + \Phi(c_i^*a^{**},c_{k_a})q \\
&= i^{\circ\circ} +k_a^\circ+ [k_{a^*c_i}+(\Phi(a^*,c_i)q)^\circ  + \Phi((a^*c_i)^*,c_{\Phi(a^*,c_i)q})q] + \Phi(c_i^*a^{**},c_{k_a})q \\
&= i^{\circ\circ} +k_a^\circ+[k_{a^*}+k_{c_i}+\Phi(c_i^*,a^{**})q]  + \Phi(c_i^*a^{**},c_{k_a})q &&\hspace{-2cm}\text{by \ref{inv8}.}
\end{align*}
Next we note that 
\[
\Phi(c_i^*,a^{**}) = 0 \AND \Phi(c_i^*a^{**},c_{k_a}) = \Phi(a^{**},c_{k_a}).
\]
Indeed, these follow from \ref{U1} and \ref{U2}, as $c_i^*$ is a unit by \ref{inv3}.  Continuing from above, we then conclude that
\begin{align*}
j^\circ+k_b+\Phi(b^*,c_j)q &= i^{\circ\circ} +k_a^\circ+k_{a^*}+k_{c_i}+\Phi(c_i^*,a^{**})q  + \Phi(c_i^*a^{**},c_{k_a})q\\
&= i^{\circ\circ}+k_{c_i} +k_a^\circ+k_{a^*}+ \Phi(a^{**},c_{k_a})q\\
&= i+0 = i &&\text{by \ref{inv4} and \ref{inv6}.}
\end{align*}
This completes the proof of \eqref{eq:iacaca}, and hence of the first equality in \eqref{eq:inv}.

We now move on to the second equality in \eqref{eq:inv}.  For this let $i,j\in M$ and $a,b\in S$.  We first claim that
\begin{equation}\label{eq:abcdacbd}
j^\circ+k_a+ \Phi(b^*,a^*,c_j,c_i)q = j^\circ+k_a+ \Phi(b^*,c_j,a^*,c_i)q,
\end{equation}
in terms of the four-parameter $\Phi$ numbers introduced at the end of Section \ref{subsect:T}.
Indeed, we have
\begin{align*}
j^\circ+k_a+ \Phi(b^*,a^*,c_j,c_i)q &= j^\circ+k_a + \Phi(a^*,c_j)q + \Phi(b^*,a^*c_j,c_i)q &&\text{by \eqref{eq:abcd}}\\
&= j^\circ+k_a + \Phi(c_j,a^*)q + \Phi(b^*,c_ja^*,c_i)q &&\text{by \ref{inv9} and \ref{inv10}}\\
&= j^\circ+k_a+ \Phi(b^*,c_j,a^*,c_i)q &&\text{by \eqref{eq:abcd} again.}
\end{align*}
With the claim established, we now write $m = i+j+\Phi(a,b)q$, and note that
\begin{equation}\label{eq:ab*}
((i,a)(j,b))^\oa = (m,ab)^\oa = (m^\circ+k_{ab}+\Phi((ab)^*,c_m)q,(ab)^*c_m),
\end{equation}
and
\begin{align}
\nonumber (j,b)^\oa(i,a)^\oa &= (j^\circ+k_b+\Phi(b^*,c_j)q,b^*c_j)(i^\circ+k_a+\Phi(a^*,c_i)q,a^*c_i) \\
\nonumber &= (j^\circ+k_b+\Phi(b^*,c_j)q+i^\circ+k_a+\Phi(a^*,c_i)q+\Phi(b^*c_j,a^*c_i)q,b^*c_ja^*c_i) \\
\nonumber &= (i^\circ+j^\circ+k_a+k_b+\Phi(b^*,c_j,a^*,c_i)q,b^*a^*c_jc_i) &&\hspace{-3.5cm}\text{by \eqref{eq:abcd} and \ref{inv10}}\\
\label{eq:b*a*} &= (i^\circ+j^\circ+k_a+k_b+\Phi(b^*,a^*,c_j,c_i)q,b^*a^*c_jc_i) &&\hspace{-3.5cm}\text{by \eqref{eq:abcdacbd}.}
\end{align}
To establish equality of the second components in~\eqref{eq:ab*} and~\eqref{eq:b*a*} we use \ref{inv1} and \ref{inv7}:
\[
(ab)^*c_m = (ab)^*c_{i+j+\Phi(a,b)q} = (ab)^* c_{\Phi(a,b)q}c_jc_i = b^*a^*c_jc_i.
\]
For the first components in~\eqref{eq:ab*} and~\eqref{eq:b*a*}, first note that by \eqref{eq:ijko} we have
\[
m^\circ= (\Phi(a,b)q+j+i)^\circ = (\Phi(a,b)q)^\circ+j^\circ+i^\circ + \Phi(c_{\Phi(a,b)q},c_j,c_i)q.
\]
Next, keeping in mind that $c_m = c_{\Phi(a,b)q}c_jc_i$ by \ref{inv1}, we have
\begin{align*}
 m^\circ&+k_{ab}+\Phi((ab)^*,c_m)q \\
 &= [(\Phi(a,b)q)^\circ+j^\circ+i^\circ + \Phi(c_{\Phi(a,b)q},c_j,c_i)q] + k_{ab} + \Phi((ab)^*,c_{\Phi(a,b)q}c_jc_i)q \\
 &= i^\circ+j^\circ+ k_{ab} + (\Phi(a,b)q)^\circ + [\Phi((ab)^*,c_{\Phi(a,b)q}c_jc_i)+\Phi(c_{\Phi(a,b)q},c_j,c_i)]q \\
 &= i^\circ+j^\circ+ k_{ab} + (\Phi(a,b)q)^\circ + [\Phi((ab)^*,c_{\Phi(a,b)q})+\Phi((ab)^*c_{\Phi(a,b)q},c_j,c_i)]q &&\text{by \eqref{eq:abcd}}\\
 &= i^\circ+j^\circ+ [k_{ab} + (\Phi(a,b)q)^\circ + \Phi((ab)^*,c_{\Phi(a,b)q})q] + \Phi(b^*a^*,c_j,c_i)q &&\text{by \ref{inv7}}\\
 &= i^\circ+j^\circ+ [k_a+k_b+\Phi(b^*,a^*)q] + \Phi(b^*a^*,c_j,c_i)q &&\text{by \ref{inv8}}\\
 &= i^\circ+j^\circ+ k_a+k_b+ \Phi(b^*,a^*,c_j,c_i)q &&\text{by \eqref{eq:abcd},}
\end{align*}
as required.
\epf

\begin{lemma}\label{lem:inv22}
Any involution ${}^\oa$ of $T$ has the form $(i,a)^\oa = (i^\circ+k_a+\Phi(a^*,c_i)q,a^*c_i)$ for maps
\[
S\to S:a\mt a^* \COMMA
M\to M:i\mt i^\circ \COMMA
S\to M:a\mt k_a \AND
M\to S:i\mt c_i 
\]
satisfying conditions \ref{inv1}--\ref{inv10}.
\end{lemma}

\pf
Throughout the proof we fix an involution ${}^\oa$ of $T$.  For any $a\in S$ and $i\in M$ we have 
\begin{align}
\label{eq:kaa*} (0,a)^\oa &= (k_a,a^*) &&\text{for some $k_a\in M$ and $a^*\in S$, and}\\
\label{eq:i0ci} (i,1)^\oa &= (i^\circ,c_i) &&\text{for some $i^\circ\in M$ and $c_i\in S$.}
\end{align}
For any $i\in M$ and $a\in S$, we use the fact that ${}^\oa$ is an involution to calculate
\begin{equation}\label{eq:iaca1}
(i,a)^\oa = ((i,1)(0,a))^\oa = (0,a)^\oa(i,1)^\oa = (k_a,a^*)(i^\circ,c_i) = (i^\circ+k_a+\Phi(a^*,c_i)q,a^*c_i).
\end{equation}
(In the first step, note that $\Phi(1,a)=0$ by \ref{U1}.)  
An analogous calculation gives
\begin{equation}\label{eq:iaca2}
(i,a)^\oa = ((0,a)(i,1))^\oa = (i^\circ+k_a+\Phi(c_i,a^*)q,c_ia^*).
\end{equation}
Comparing the first coordinates in \eqref{eq:iaca1} and \eqref{eq:iaca2}, we see that \ref{inv10} holds.  The second coordinates give:
\begin{equation}\label{eq:a*ci}
a^*c_i = c_ia^* \qquad\text{for all $i\in M$ and $a\in S$.}
\end{equation}
We now consider the remaining conditions \ref{inv1}--\ref{inv9} in turn.

\pfitem{\ref{inv1} and \ref{inv2}}  For $i,j\in M$ we use \eqref{eq:i0ci} to calculate
\begin{align}
\nonumber ((i+j)^\circ,c_{i+j}) = (i+j,1)^\oa = ((j,1)(i,1))^\oa = (i,1)^\oa(j,1)^\oa &= (i^\circ,c_i)(j^\circ,c_j) \\
\label{eq:i+jo} &= (i^\circ+j^\circ+\Phi(c_i,c_j)q,c_ic_j).
\end{align}
Looking at the second coordinates gives $c_ic_j = c_{i+j}$, which is the first part of \ref{inv1}.  The first coordinates give the $(i+j)^\circ = i^\circ+j^\circ+\Phi(c_i,c_j)q$ part of \ref{inv2}.  The other parts of \ref{inv1} and \ref{inv2} follow from commutativity of $M$.

\pfitem{\ref{inv3} and \ref{inv4}}  For $i\in M$ we use \eqref{eq:i0ci}, \eqref{eq:iaca1} and \eqref{eq:iaca2} to calculate
\begin{align*}
(i,1) = (i,1)^{\oa\oa} = (i^\circ,c_i)^\oa &= (i^{\circ\circ} + k_{c_i} + \Phi(c_i^*,c_{i^\circ})q,c_i^*c_{i^\circ})\\
&= (i^{\circ\circ} + k_{c_i} + \Phi(c_{i^\circ},c_i^*)q,c_{i^\circ}c_i^*).
\end{align*}
The second coordinates give \ref{inv3}.  In particular, $c_{i^\circ}$ and $c_i^*$ are units, so it follows from \ref{U1} that $\Phi(c_i^*,c_{i^\circ}) = 0$.  Combining this with the first coordinates in the above calculation yields \ref{inv4}.

\pfitem{\ref{inv5} and \ref{inv6}}  For $a\in S$ we use \eqref{eq:kaa*}, \eqref{eq:iaca1} and \eqref{eq:iaca2} to calculate
\begin{align*}
(0,a) = (0,a)^{\oa\oa} = (k_a,a^*)^\oa &= (k_a^\circ+k_{a^*}+\Phi(a^{**},c_{k_a})q,a^{**}c_{k_a})\\
&= (k_a^\circ+k_{a^*}+\Phi(c_{k_a},a^{**})q,c_{k_a}a^{**}).
\end{align*}

\pfitem{\ref{inv7} and \ref{inv8}}  For $a,b\in S$ we use \eqref{eq:kaa*} to calculate
\[
((0,a)(0,b))^\oa = (0,b)^\oa(0,a)^\oa = (k_b,b^*)(k_a,a^*) = (k_a+k_b+\Phi(b^*,a^*)q,b^*a^*).
\]
Instead using \eqref{eq:iaca1} we have
\[
((0,a)(0,b))^\oa = (\Phi(a,b)q,ab)^\oa = ((\Phi(a,b)q)^\circ+k_{ab}+\Phi((ab)^*,c_{\Phi(a,b)q})q,(ab)^*c_{\Phi(a,b)q}).
\]

\pfitem{\ref{inv9}}  For $i\in M$ and $a\in S$, it follows from \ref{inv5}, \eqref{eq:a*ci} and \ref{inv1} that
\[
ac_i = (a^{**}c_{k_a})c_i = c_i(a^{**}c_{k_a}) = c_ia.  \qedhere
\]
\epf

\subsection{Consequences}\label{subsect:consequences}

Now that we have completed the proof of Theorem \ref{thm:inv2}, we record some useful consequences.  For the next statement, recall that $G_M$ and $G_S$ denote the groups of units of $M$ and $S$.

\newpage

\begin{prop}\label{prop:group}
Let ${}^\oa$ be an involution of a unital twisted product $T = M\times_\Phi^q S$, as in Theorem~\ref{thm:inv2}.  Then 
\ben
\item \label{group1} $S\to S:a\mt a^*$ restricts to a group anti-morphism $G_S\to G_S$,
\item \label{group2} $M\to M:i\mt i^\circ$ restricts to a group morphism $G_M\to G_M$,
\item \label{group3} $S\to M:a\mt k_a$ restricts to a group morphism $G_S\to G_M$,
\item \label{group4} $M\to S:i\mt c_i$ restricts to a group morphism $G_M\to G_S$.
\een
\end{prop}

\pf
First, substituting $a=1$ and $i=0$ in \eqref{eq:kaa*} and \eqref{eq:i0ci}, respectively, and remembering that $(0,1)^\oa=(0,1)$, we have
\begin{equation}\label{eq:0circ}
0^\circ=0=k_1 \AND 1^*=1=c_0.
\end{equation}
This says that all four maps under consideration map the identity of $S$ or $M$ to the identity of~$S$ or $M$, as appropriate.  Consequently, it suffices to show that the restriction to each domain group is an \mbox{(anti-)morphism}, as it will automatically follow that the image of each such restriction is contained in the appropriate target group.  Indeed, it is easy to show that the image of a monoid \mbox{(anti-)morphism} $G\to U$, with $G$ a group, is contained in the group of units $G_U$.

\pfitem{\ref{group1}}  If $a,b\in G_S$, then by \ref{inv7}, \ref{U1} and \eqref{eq:0circ}, we have
\[
(ab)^* = (ab)^*c_0 = (ab)^*c_{\Phi(a,b)q} = b^*a^*.
\]

\pfitem{\ref{group3}}  If $a,b\in G_S$, then since also $a^*,b^*,(ab)^*\in G_S$ by part \ref{group1}, it follows from \ref{U1} that
\[
\Phi(b^*,a^*) = \Phi(a,b) = \Phi((ab)^*,c_{\Phi(a,b)q}) = 0.
\]
Combining this with \ref{inv8} and $0^\circ=0$ (cf.~\eqref{eq:0circ}), we then obtain $k_a+k_b=k_{ab}$.

\pfitem{\ref{group4}}  By \ref{inv1}, $M\to S:i\mt c_i$ is a morphism, and so too therefore is its restriction to $G_M$.

\pfitem{\ref{group2}}  If $i,j\in G_M$, then by part \ref{group4} and \ref{U1} we have $\Phi(c_i,c_j)=0$.  Combining this with \ref{inv2} gives ${(i+j)^\circ = i^\circ+j^\circ}$.
\epf

An important special case occurs when the elements $c_i$ from Theorem \ref{thm:inv2} are all units:

\begin{prop}\label{prop:unit}
Let ${}^\oa$ be an involution of a unital twisted product $T = M\times_\Phi^q S$, as in Theorem~\ref{thm:inv2}, and suppose $c_i\in G_S$ for all $i\in M$.  Then for all $i,j\in M$ and $a,b\in S$ we have
\ben
\item \label{unit1} $(i,a)^\oa=(i^\circ+k_a,a^*c_i)$,
\item \label{unit2} $(i+j)^\circ=i^\circ+j^\circ$,
\item \label{unit3} $k_a^\circ+k_{a^*}=0$,
\item \label{unit4} $k_a+k_b+\Phi(b^*,a^*)q = k_{ab} + (\Phi(a,b)q)^\circ$,
\item \label{unit5} $a\R^Sb \implies a^*\L^Sb^*$, 
\item \label{unit6} $a\L^Sb \implies a^*\R^Sb^*$,
\item \label{unit7} $a\H^Sb \implies a^*\H^Sb^*$.
\een
\end{prop}

\pf
\firstpfitem{\ref{unit1}}  This follows from the rule for $(i,a)^\oa$ in Theorem \ref{thm:inv2}, given that $\Phi(a^*,c_i)=0$ by~\ref{U1}.

\pfitem{\ref{unit2}--\ref{unit4}}  These follow respectively from \ref{inv2}, \ref{inv6} and \ref{inv8}, since ${\Phi(c_i,s)=\Phi(s,c_i)=0}$ for any $i\in M$ and $s\in S$, by \ref{U1}.

\pfitem{\ref{unit5}--\ref{unit7}}  It suffices to prove \ref{unit5}, as \ref{unit6} is symmetrical, while \ref{unit7} follows by combining \ref{unit5} and~\ref{unit6}.  First note that for any $a,b\in S$, it follows from \ref{inv7} and \ref{inv9}, and the assumption that the $c_i$ are units, that
\[
(ab)^* = b^*a^*g \qquad\text{for some central unit $g$ of $S$.}
\]
Now suppose $a\R^Sb$, so that $a=bs$ and $b=at$ for some $s,t\in S$.  Then for some central units $g$ and $h$, we have
\[
a^* = (bs)^* = s^*b^*g = gs^*\cdot b^*, \ANDSIM b^* = ht^*\cdot a^*,
\]
so that indeed $a^* \L^S b^*$.
\epf

\begin{rem}\label{rem:unit}
Note that elements of the form $c_{i^\circ}$ are always units by \ref{inv3}.  Consequently, if the map ${}^\circ$ is surjective, then every $c_i$ is a unit, in which case Proposition \ref{prop:unit} applies.  For example, if we had $k_{c_i}=0$ for all $i\in M$, then surjectivity of ${}^\circ$ would follow from \ref{inv4}.  By Proposition~\ref{prop:group}\ref{group4}, $M$ being a group would also imply that each $c_i$ is a unit.
\end{rem}

We have just considered the situation in which each $c_i$ is a unit of $S$.  The stronger condition that each $c_i=1$ is also very important; in fact, it is equivalent to purity of the involution ${}^\oa$:

\begin{prop}\label{prop:ext}
Let ${}^\oa$ be an involution of a unital twisted product $T = M\times_\Phi^q S$, as in Theorem~\ref{thm:inv2}.  Then ${}^\oa$ is pure if and only if $c_i=1$ for all $i\in M$, in which case ${}^\oa$ extends the involution ${}^*$ of $S$.  In this case, ${}^\circ$ is also an involution of $M$.
\end{prop}

\pf
Suppose first that ${}^\oa$ is pure, and hence extends some involution $^\times$ of $S$.  Looking at~\eqref{eq:i0ci}, and keeping in mind that $1^\times=1$, we see that for any $i\in M$,
\[
(i^\circ,c_i) = (i,1)^\oa = (j,1^\times) = (j,1) \qquad\text{for some $j\in M$,}
\]
so indeed $c_i=1$.

Conversely, suppose $c_i=1$ for all $i\in M$.  Then for any $i\in M$ and $a\in S$, the second coordinate of $(i,a)^\oa$ is $a^*c_i = a^*$.  Since this depends only on $a$, it follows that ${}^\oa$ is pure.  Since~${}^\oa$ extends~${}^*$, Lemma \ref{lem:pure} tells us that the latter is an involution of $S$.

For the final assertion, note that $c_i=1$ implies $k_{c_i} = k_1 = 0$ for all ${i\in M}$; cf.~\eqref{eq:0circ}.  It then follows from \ref{inv2} and \ref{inv4} that ${}^\circ$ is an involution of $M$, again remembering ${\Phi(c_i,c_j)=\Phi(1,1)=0}$.
\epf

Here is a useful special case in which \emph{every} involution of a twisted product is pure.

\begin{cor}\label{cor:inv}
If $T=M\times_\Phi^qS$ is a unital twisted product, and if $S$ has trivial centre, then every involution of $T$ is pure, and hence has the form described in Theorem \ref{thm:inv1}.
\end{cor}

\pf
Suppose ${}^\oa$ is an involution of $T$, as in Theorem \ref{thm:inv2}.  Since the~$c_i$ are central by \ref{inv9}, it follows by assumption that $c_i=1$ for all~${i\in M}$.  The result now follows by Proposition \ref{prop:ext}.
\epf

\subsection{Examples}\label{subsect:examples}

We conclude this section with some examples of \emph{non-pure} involutions of a (tight) twisted product ${T=M\times_\Phi^qS}$.  For such involutions to exist, we require $S$ to have non-trivial centre (cf.~Corollary~\ref{cor:inv}), and in fact we will use a commutative monoid $S$.  The stated involutions of~$T$ are not pure because we do not have $c_i=1$ for all $i\in M$ (cf.~Proposition \ref{prop:ext}).  On the other hand, the maps $M\to M:i\mt i^\circ$ and $S\to S:a\mt a^*$ are themselves involutions.  

\begin{eg}\label{eg:nonext}
Consider the Brauer monoid $S=\B_2$, and write $\B_2 = \{1,s,e\}$, where 
\[
1 = \custpartn{1,2}{1,2}{\stline11\stline22} \COMMA s=s_1 = \custpartn{1,2}{1,2}{\stline12\stline21} \AND e=e_1 = \custpartn{1,2}{1,2}{\uarc12\darc12}.
\]
Also let $T = \Z\times_\Phi^1\B_2$ be the (tight) twisted product arising from the canonical float-counting twisting~$\Phi$, the additive group of integers $M=\Z$, and with $q=1$.  We claim that $T$ has an involution ${}^\oa$ given by
\[
(i,a)^\oa = (i,s^ia) \qquad\text{for all $i\in\Z$ and $a\in\B_2$.}
\]
Indeed, this arises in the manner described in Theorem \ref{thm:inv2} upon setting
\[
i^\circ=i \COMMA c_i = s^i \COMMA a^*=a \AND k_a = 0 \qquad\text{for all $i\in\Z$ and $a\in\B_2$.}
\]
Verification of properties \ref{inv1}--\ref{inv10} is straightforward.  Note that \ref{inv7} requires separate cases, depending on whether stated elements of $\B_2$ are units, or equal to $e$.  Indeed,~\ref{inv7} reduces to $ab\cdot c_{\Phi(a,b)}=ba$ for all $a,b\in\B_2$.  If either of $a$ or $b$ is a unit, then this follows from commutativity of~$\B_2$, and $c_{\Phi(a,b)}=c_0=1$.  Otherwise, $a=b=e$, and $ab\cdot c_{\Phi(a,b)}=ba$ holds because $e$ is the zero of $\B_2$.
\end{eg}

\begin{eg}\label{eg:nonext2}
Consider again the product $T = \Z\times_\Phi^1\B_2$, and again write $\B_2 = \{1,s,e\}$.  We claim that $T$ has an involution ${}^\oa$ given by
\[
(i,a)^\oa = (-i-\non(a),s^ia),
\]
where $\non(a)$ is the number of non-transversals of $a$.  Indeed, this also arises from Theorem \ref{thm:inv2}, this time with
\[
i^\circ=-i \COMMA c_i = s^i \COMMA a^*=a \AND k_a = -\non(a) = \begin{cases}
0 &\text{if $a=1$ or $s$}\\
-2 &\text{if $a=e$.}
\end{cases}
\]
Verification of properties \ref{inv1}--\ref{inv10} is again straightforward; this time \ref{inv5}, \ref{inv7} and \ref{inv8} require separate cases.  Observe that for this involution, the map $a\mt k_a$ is \emph{not} a morphism.
\end{eg}

We will encounter similar non-pure twistings of $\Z\times_\Psi^1\B_2$, arising from the rank-based twisting~$\Psi$, in Section \ref{sect:Psi}.

\section{Star-regularity and regular-starity of twisted products}\label{sect:*reg}

Now that we have classified the involutions of a unital twisted product $T=M\times_\Phi^qS$, we wish to characterise those that give $T$ the structure of a star-regular or regular star-monoid.  The former is done in Theorems \ref{thm:*reggen} and \ref{thm:*reg} (in Section \ref{subsect:*reg}), and the latter in Theorem \ref{thm:reg*} (in Section~\ref{subsect:reg*}).  The first theorem for star-regularity concerns general involutions, and the second is for pure ones.  There is only one theorem for regular-starity, as it turns out that $T$ can only be a regular star-monoid with respect to a pure involution.

\subsection{Star-regularity}\label{subsect:*reg}

We first consider star-regularity of twisted products $T=M\times_\Phi^qS$.  As in Section \ref{sect:*}, we will again give two statements, one for pure involutions (Theorem \ref{thm:*reg}), and one for the general case (Theorem~\ref{thm:*reggen}).  The pure theorem assumes only semi-unitality, whereas full unitality is required in the general theorem.  This time, however, we prove the general theorem first, for reasons explained in Remark \ref{rem:*reggenproof}.

In the following results, note that if the commutative monoid $M$ is regular, then it is necessarily a union of (abelian) groups.  In this case, for any $i\in M$ we will write~$-i$ for the inverse of~$i$ in the group $\H^M$-class $H_i^M$ containing $i$.  Since ${\D^M}={\H^M}$ (by commutativity), this element~$-i$ is in fact the unique inverse of $i$ in $M$.  Consequently, $M$ is an \emph{inverse monoid}, indeed a commutative \emph{Clifford monoid}, or equivalently a \emph{strong semilattice of abelian groups}; see \cite[Section~4.2]{Howie1995} for the definitions of these terms.  For any $i\in M$, $i-i$ is the identity of the group $H_i^M$.  Consequently, we have
\begin{equation}\label{eq:HM}
i\H^Mj \Iff i-i = j-j \qquad\text{for all $i,j\in M$.}
\end{equation}

\begin{thm}\label{thm:*reggen}
Let $T = M\times_\Phi^qS$ be a unital twisted product, and suppose~${}^\oa$ is an involution of~$T$, so that
\[
(i,a)^\oa = (i^\circ+k_a+\Phi(a^*,c_i)q,a^*c_i) \qquad\text{for all $i\in M$ and $a\in S$,}
\]
as in Theorem \ref{thm:inv2}.  Then $T$ is $\oa$-regular if and only if the following are all satisfied:
\begin{enumerate}[label=\textup{\textsf{(SR\arabic*)}},leftmargin=12mm]
\item \label{*reggen1} for all $a\in S$ we have $a^*\H^S a^+$ for some inverse $a^+$ of $a$,
\item \label{*reggen2} $k_a$ and $\Phi(a,a^+,a)q$ are units of $M$ for all $a\in S$,
\item \label{*reggen3} $M$ is regular,
\item \label{*reggen4} $i^\circ \H^M i$ for all $i\in M$,
\item \label{*reggen5} $c_i$ is a unit of $S$ for all $i\in M$.
\een
When $T$ is $\oa$-regular, the involution and Moore--Penrose inverses are given by the rules
\[
(i,a)^\oa = (i^\circ+k_a,a^*c_i)
\AND
(i,a)^\op = (-i-\Phi(a,a^+,a)q,a^+) \qquad\text{for all $i\in M$ and $a\in S$.}
\]
\end{thm}

\begin{rem}\label{rem:*reggenproof}
Before we begin the proof of the theorem, we observe that unitality is assumed in the statement to ensure that the involution ${}^\oa$ necessarily has the form described in Theorem~\ref{thm:inv2}.  However, in the argument below, we only use the weaker assumption of semi-unitality.  This point will be important when we come to prove the next theorem.
\end{rem}

\pf[\bf Proof of Theorem \ref{thm:*reggen}.]
First assume that \ref{*reggen1}--\ref{*reggen5} all hold.  Since each $c_i$ is a unit, the simplified expression for $(i,a)^\oa$ follows from Proposition \ref{prop:unit}\ref{unit1}.  Looking at~\ref{*reggen2}, we see that either $q$ is a unit of $M$, or else $\Phi(a,a^+,a)=0$ for all $a\in S$.  We will soon consider these two cases separately.  But first we observe that $(i,a)^\op$ (as given in the statement) is an inverse of~$(i,a)$ for all $i\in M$ and $a\in S$.  Indeed, we have
\[
(i,a)(i,a)^\op(i,a) = (2i+(-i-\Phi(a,a^+,a)q) + \Phi(a,a^+,a)q,aa^+a) = (i,a).
\]
A similar calculation gives $(i,a)^\op(i,a)(i,a)^\op = (i,a)^\op$, keeping in mind that $\Phi(a^+,a,a^+) = \Phi(a,a^+,a)$ by Lemma \ref{lem:aba}\ref{aba2}.

\pfcase1  Suppose first that $q$ is a unit of $M$.  In this case we complete the proof that $T$ is $\oa$-regular by showing that $(i,a)^\op \H^T (i,a)^\oa$ for all $i\in M$ and $a\in S$.  Since $T$ is Green-deomposable by Lemma \ref{lem:Green}, this amounts to showing that
\[
-i-\Phi(a,a^+,a)q \ \H^M\  i^\circ+k_a
\AND
a^+ \H^S a^*c_i.
\]
For the former we have
\begin{align*}
-i-\Phi(a,a^+,a)q &\ \ \H^M \ \ -i &&\text{as $q$ is a unit}\\
&\ \ \H^M \ \ i^\circ &&\text{by \ref{*reggen4} and $i\H^M -i$}\\
&\ \ \H^M \ \ i^\circ+k_a &&\text{as $k_a$ is a unit by \ref{*reggen2}.}
\end{align*}
Since $a^+ \H^S a^*$ by \ref{*reggen1}, we can prove the latter by showing that $a^* \H^S a^*c_i$.  Now, $c_i$ being a unit (by \ref{*reggen5}) means that $a^* \R^S a^*c_i$.  And since $c_i$ is central in $S$ by \ref{inv9}, we also have ${a^* \L^S c_ia^* = a^*c_i}$.  The last two conclusions give $a^* \H^S a^*c_i$.

\pfcase2  Now we assume that $q$ is not a unit.  By \ref{*reggen2}, this implies that $\Phi(a,a^+,a) = 0$ for all $a\in S$.  Since
\[
\Phi(a,a^+,a) = \Phi(a,a^+)+\Phi(aa^+,a) = \Phi(a,a^+a) + \Phi(a^+,a),
\]
each $\Phi$ term in the above sums is also $0$.  It also follows from Lemma \ref{lem:aba}\ref{aba1} that
\[
\Phi(aa^+,aa^+) = \Phi(aa^+,a) = 0, \ANDSIM \Phi(a^+a,a^+a) = 0.
\]
In this case we can complete the proof that $T$ is $\oa$-regular by showing that $(i,a)(i,a)^\op$ and $(i,a)^\op(i,a)$ are both fixed by ${}^\oa$ for all $i\in M$ and $a\in S$; cf.~\eqref{eq:axa}.  We just do this for $(i,a)(i,a)^\op$, as the other is analogous.  To do so, and keeping $\Phi(a,a^+,a)=0=\Phi(a,a^+)$ in mind, we write
\[
(j,e) = (i,a)(i,a)^\op = (i,a)(-i,a^+) = (i-i,aa^+).
\]
So we need to show that $(j,e) = (j,e)^\oa = (j^\circ+k_e,e^*c_j)$, i.e.~that
\begin{equation}\label{eq:je}
j^\circ+k_e = j \AND e^*c_j = e.
\end{equation}
Proposition \ref{prop:unit}\ref{unit2}, the map ${}^\circ:M\to M$ is a morphism.  Keeping in mind that $M$ is an inverse monoid, it follows that
\[
j^\circ = (i-i)^\circ = i^\circ - i^\circ = i-i = j,
\]
where in the second last step we used \eqref{eq:HM} and \ref{*reggen4}.  By Proposition \ref{prop:unit}\ref{unit4}, and remembering that $e=aa^+$ is an idempotent, and that $\Phi(e,e)=\Phi(aa^+,aa^+)=0$ and $0^\circ=0$ (cf.~\eqref{eq:0circ}), we have
\[
2k_e + \Phi(e^*,e^*)q = k_e +(\Phi(e,e)q)^\circ = k_e + 0^\circ = k_e.
\]
Subtracting the unit $k_e$ (cf.~\ref{*reggen2}) yields $k_e + \Phi(e^*,e^*)q = 0$, which implies that $\Phi(e^*,e^*)q$ is a unit.  Since~$q$ is not a unit, it follows that $\Phi(e^*,e^*)=0$, and so $0 = k_e+\Phi(e^*,e^*)q = k_e$.  Combined with $j^\circ=j$ (established above), this completes the proof of the first equality in \eqref{eq:je}.  

Next note that since $j=i-i$ is an idempotent of $M$, it follows from \ref{inv1} that $c_j$ is an idempotent of $S$.  But $c_j$ is also a unit by \ref{*reggen5}, so in fact $c_j=1$.  Thus, to establish the second equality in \eqref{eq:je} we need to show that
\[
e^* = e.
\]
We do this by showing that $e^*$ is an idempotent in the $\H^S$-class $H_e^S$ of $e$; since $e=aa^+$ is itself an idempotent, it will follow that indeed $e^*=e$.  To see that $e^*$ is an idempotent, we use $\Phi(e,e)=0$, $c_0=1$ (cf.~\eqref{eq:0circ}), and \ref{inv7} to calculate
\[
e^* = (ee)^*c_0 = (ee)^*c_{\Phi(e,e)q} = e^*e^*.
\]
Next note that since $a^+$ is an inverse of $a$, we have $a\L^S a^+a\R^Sa^+$, which says that the $\H^S$-class $L_a^S \cap R_{a^+}^S$ contains the idempotent $a^+a$.  It follows from standard semigroup theory (see for example \cite[Proposition 2.3.7]{Howie1995}) that
\begin{equation}\label{eq:He}
H_a^S \cdot H_{a^+}^S = H_{aa^+}^S = H_e^S.
\end{equation}
Next we use \ref{inv7}, together with $\Phi(a,a^+)=0$ and $c_0=1$, to calculate
\[
e^* = (aa^+)^* = (aa^+)^*c_{\Phi(a,a^+)q} = (a^+)^*a^*.
\]
Thus, given \eqref{eq:He}, we can complete the proof that $e^*\in H_e^S$ by showing that
\[
(a^+)^*\in H_a^S \AND a^* \in H_{a^+}^S.
\]
The latter follows from \ref{*reggen1}.  For the former, it follows from Proposition \ref{prop:unit}\ref{unit7} and $a^+ \H^S a^*$ that $(a^+)^* \H^S a^{**}$, and it follows from~\ref{inv5} (and $c_{k_a}$ being a central unit) that $a^{**} \H^S a$.  Combining the last two conclusions, we obtain $(a^+)^* \H^S a$, which completes the proof that $(a^+)^* \in H_a^S$, and hence (as explained above) the proof of $\oa$-regularity of $T$.

\aftercases

Conversely, suppose $T$ is $\oa$-regular, and denote the Moore--Penrose inverse of $(i,a)$ by $(i,a)^\op$.  For each $i\in M$ and $a\in S$ we have
\[
(i,1)^\op = (i^\di,d_i) \AND (0,a)^\op = (m_a,a^+) \qquad\text{for some $i^\di,m_a\in M$ and $d_i,a^+\in S$.}
\]
We now establish items \ref{*reggen1}--\ref{*reggen5}.

\pfitem{\ref{*reggen1} and \ref{*reggen2}}  For any $a\in S$ we have
\[
(0,a) = (0,a)(0,a)^\op(0,a) = (m_a+\Phi(a,a^+,a)q,aa^+a).
\]
This tells us that $\Phi(a,a^+,a)q$ is a unit (with inverse $m_a$), and that $a=aa^+a$.  We obtain $a^+=a^+aa^+$ from expanding $(0,a)^\op = (0,a)^\op(0,a)(0,a)^\op$.  Since
\[
(m_a,a^+) = (0,a)^\op \H^T (0,a)^\oa = (k_a,a^*),
\]
using \eqref{eq:kaa*} in the last step, it follows that $a^+ \H^S a^*$ and $k_a \H^M m_a$ (see Remark \ref{rem:Green}).  The former completes the proof of \ref{*reggen1}.  Since $m_a$ is a unit, so too therefore is $k_a$, which completes the proof of \ref{*reggen2}.  

\pfitem{\ref{*reggen3}--\ref{*reggen5}}  For any $i\in M$ we have $\Phi(1,d_i,1) = \Phi(1,d_i)+\Phi(1d_i,1) = 0$ by \ref{U1'}, and so
\[
(i,1) = (i,1)(i,1)^\op(i,1) = (2i+i^\di+\Phi(1,d_i,1)q,d_i) = (2i+i^\di,d_i).
\]
It follows that $d_i=1$, and also that $i=2i+i^\di$.  Looking instead at $(i,1)^\op = (i,1)^\op(i,1)(i,1)^\op$, and using $d_i=1$ and $\Phi(1,1,1)=0$, we obtain $i^\di=2i^\di+i$.  Together, $i=2i+i^\di$ and $i^\di=2i^\di+i$ show that $i^\di$ is an inverse of $i$ (giving \ref{*reggen3}), and also that $i \H^M i^\di$.  Since
\[
(i^\di,1) = (i,1)^\op \H^T (i,1)^\oa = (i^\circ,c_i),
\]
using \eqref{eq:i0ci} in the last step, it follows that $i^\di\H^M i^\circ$ and $1\H^S c_i$ (again see Remark \ref{rem:Green}).  The latter says that~$c_i$ is a unit of $S$ (giving \ref{*reggen5}), and the former combines with $i\H^M i^\di$ to show that $i^\circ\H^M i$ (giving~\ref{*reggen4}).
\epf

\newpage

\begin{rem}
As we observed in the above proof, the condition that the elements $\Phi(a,a^+,a)q$ are all units of $M$ (in \ref{*reggen2}) is equivalent to the disjunction:
\bit
\item $q$ is a unit of $M$, or else
\item $\Phi(a,a^+,a) = 0$ for all $a\in S$.
\eit
It was shown in \cite[Theorem 7.4]{EGPR2026} that a \emph{tight} twisted product $T = M\times_\Phi^qS$ is regular if and only if both $M$ and $S$ are regular, and either $q$ is a unit of $M$ or else for every $\D^S$-class $D$ of $S$ we have $\Phi(e,e)=0$ for some idempotent $e\in E(D)$.  When $T$ is $\oa$-regular, it is of course regular, and it is easy to see that this latter condition from \cite{EGPR2026} is implied by the above disjunction.  Indeed, fix a $\D^S$-class $D$, and let $a\in D$ be arbitrary.  Then since $a$ and $a^+$ are mutual inverses the idempotents $e=aa^+$ and $f=a^+a$ both belong to $D$.  Moreover, we showed in the above proof that $\Phi(a,a^+,a)=0$ implies $\Phi(e,e)=0=\Phi(f,f)$.
\end{rem}

Theorem \ref{thm:*reggen} simplifies in the special case of a pure involution, and we note that in the following statement, we relax the assumption of unitality to semi-unitality.

\begin{thm}\label{thm:*reg}
Let $T = M\times_\Phi^qS$ be a semi-unital twisted product, and suppose~${}^\oa$ is a pure involution of~$T$, so that
\[
(i,a)^\oa = (i^\circ+k_a,a^*) \qquad\text{for all $i\in M$ and $a\in S$,}
\]
as in Theorem \ref{thm:inv1}.  Then $T$ is $\oa$-regular if and only if the following are all satisfied:
\ben
\item \label{*reg1} $S$ is $*$-regular,
\item \label{*reg2} $M$ is $\circ$-regular,
\item \label{*reg3} $k_a$ and $\Phi(a,a^+,a)q$ are units of $M$ for all $a\in S$ (where ${}^+$ is the Moore--Penrose inverse in~$S$).
\een
When $T$ is $\oa$-regular, Moore--Penrose inverses are given by the rule
\[
(i,a)^\op = (-i-\Phi(a,a^+,a)q,a^+) \qquad\text{for all $i\in M$ and $a\in S$.}
\]
\end{thm}

\pf
As observed in Remark \ref{rem:*reggenproof}, the proof of Theorem \ref{thm:*reggen} works under the assumptions of the current theorem to show that $T$ is $\oa$-regular (for the pure involution ${}^\oa$) if and only if~\mbox{\ref{*reggen1}--\ref{*reggen5}} hold.
By Proposition \ref{prop:ext}, ${}^*$ and ${}^\circ$ are involutions of~$S$ and~$M$ respectively.  Thus, \ref{*reggen1} is equivalent to $S$ being $*$-regular.  Conditions~\ref{*reggen3} and~\ref{*reggen4} are equivalent to~$M$ being $\circ$-regular, as $i \H^M -i$ for all $i\in M$ when $M$ is regular.  Finally,~\ref{*reggen5} is trivial in this case, as each $c_i=1$ by Proposition \ref{prop:ext}.
\epf

It is convenient to record the following special case, which concerns the involution given in Corollary \ref{cor:*sym}.

\begin{cor}\label{cor:*sym*reg}
Let $T = M\times_\Phi^qS$ be a semi-unital twisted product, and suppose $\Phi$ is $*$-symmetric with respect to an involution ${}^*$ of $S$.  Then $T$ is $\oa$-regular with respect to the involution ${(i,a)^\oa = (i,a^*)}$ if and only if the following are all satisfied:
\ben
\item $S$ is $*$-regular,
\item $M$ is regular,
\item $\Phi(a,a^+,a)q$ is a unit of $M$ for all $a\in S$ (where ${}^+$ is the Moore--Penrose inverse in $S$).
\een
\end{cor}

\pf
As in the proof of Corollary \ref{cor:*sym}, here we have $i^\circ=i$ and $k_a=0$ for all $i\in M$ and $a\in S$.  The result now follows from Theorem \ref{thm:*reg}, since each $k_a=0$ is a unit, and since $M$ is regular if and only if it is $\circ$-regular (as each $i^\circ = i \H^M -i$ when $M$ is regular).
\epf

\subsection{Regular-starity}\label{subsect:reg*}

We can also characterise the stronger property of $M\times_\Phi^qS$ being a regular star-monoid.  We first show that regular-starity can only arise from a \emph{pure} involution.

\begin{prop}\label{prop:reg*}
Let $T = M\times_\Phi^qS$ be a unital twisted product, and let ${}^\oa$ be an involution of~$T$.  If $T$ is a regular $\oa$-monoid, then ${}^\oa$ is pure.
\end{prop}

\pf
As ever, we may assume that ${}^\oa$ has the form described in Theorem \ref{thm:inv2}.  By Proposition~\ref{prop:ext}, it suffices to show that $c_i=1$ for all $i\in M$, so fix some such $i$.  By \eqref{eq:i0ci} and \ref{U1} we have
\[
(i,1) = (i,1)(i,1)^\oa(i,1) = (i,1)(i^\circ,c_i)(i,1) = (2i+i^\circ,c_i),
\]
which gives $c_i=1$.
\epf

Thus we only need to characterise regular-starity for pure involutions, and it turns out that we can do this under the weaker assumption of semi-unitality:

\begin{thm}\label{thm:reg*}
Let $T = M\times_\Phi^qS$ be a semi-unital twisted product, and let ${}^\oa$ be a pure involution of~$T$, so that
\[
(i,a)^\oa = (i^\circ+k_a,a^*) \qquad\text{for all $i\in M$ and $a\in S$,}
\]
as in Theorem \ref{thm:inv1}.  Then $T$ is a regular $\oa$-monoid if and only if the following are all satisfied:
\begin{enumerate}[label=\textup{\textsf{(RS\arabic*)}},leftmargin=12mm]
\item \label{reg*a} $S$ is a regular $*$-monoid,
\item \label{reg*b} $M$ is a regular $\circ$-monoid,
\item \label{reg*c} $k_a+\Phi(a,a^*,a)q=0$ for all $a\in S$.
\een
\end{thm}

\pf
Suppose first that $T$ is a regular $\oa$-monoid.  Now, ${}^*$ and ${}^\circ$ are involutions of $S$ and $M$, respectively, by Theorem \ref{thm:inv1}.  As in the proof of Proposition \ref{prop:reg*}, we have $(i,1)=(2i+i^\circ,1)$ for all $i\in M$; the first coordinates say that $M$ is a regular $\circ$-monoid.
If $a\in S$, then by \eqref{eq:kaa*} we have
\[
(0,a) = (0,a)(0,a)^\oa(0,a) = (0,a)(k_a,a^*)(0,a) = (k_a+\Phi(a,a^*,a)q,aa^*a).
\]
The second coordinates tell us that $S$ is a regular $*$-monoid, and the first coordinates give the identity $k_a+\Phi(a,a^*,a)q=0$.

Conversely, if ${}^\oa$ satisfies conditions~\ref{reg*a}--\ref{reg*c}, then for any $i\in M$ and $a\in S$ we have
\[
(i,a)(i,a)^\oa(i,a) = (i,a)(i^\circ+k_a,a^*)(i,a) = (2i+i^\circ+k_a+\Phi(a,a^*,a)q,aa^*a) = (i,a).  \qedhere
\]
\epf

It is again convenient to explicitly treat the special case concerning the involution given in Corollary \ref{cor:*sym}.
For the next statement, recall that a group $G$ (written additively) has \emph{exponent}~$m$ if $mg$ is the identity for every $g\in G$, and if $m$ is the least positive integer with respect to this property.  We have already noted that if our (commutative) monoid~$M$ is regular, then it is a union of (commutative) groups.  When this is the case we say~$M$ has exponent $m$ if for every~${i\in M}$,~$mi$ is the identity of the subgroup of $M$ containing~$i$ (which is equivalent to $M$ satisfying the identity~$(m+1)i=i$), with~$m$ being the least positive integer with respect to this property.  Note then that every subgroup of $M$ has exponent dividing~$m$.
Note also that $M$ being regular and having exponent at most $2$ is equivalent to it satisfying the identity $3i=i$ for all~$i\in M$.

\newpage

\begin{prop}\label{prop:*symreg*}
Let $T = M\times_\Phi^qS$ be a semi-unital twisted product, and suppose $\Phi$ is $*$-symmetric with respect to an involution ${}^*$ of $S$.  Then $T$ is a regular $\oa$-monoid with respect to the involution $(i,a)^\oa = (i,a^*)$ if and only if the following are all satisfied:
\ben
\item \label{*symreg*1} $S$ is a regular $*$-monoid,
\item \label{*symreg*3} $M$ is regular and has exponent at most $2$,
\item \label{*symreg*2} $\Phi(a,a^*,a)q = 0$ for all $a\in S$.
\een
\end{prop}

\pf
For this involution we have $i^\circ=i$ and $k_a=0$ for $i\in M$ and $a\in S$.  So \ref{reg*b} is equivalent to $M$ satisfying the identity $3i=i$, and we noted above that this is equivalent to~$M$ being regular and having exponent at most $2$.  Since $k_a=0$, \ref{reg*c} is equivalent to having $\Phi(a,a^*,a)q = 0$ for all $a\in S$.
\epf

\section{Twisted diagram monoids}\label{sect:T}

We now consider involutions, and star-regularity and regular-starity, of a twisted diagram monoid $T = M\times_\Phi^qS$, where $S$ is one of~$\P_n$,~$\PP_n$, $\PB_n$, $\M_n$, $\B_n$ or~$\TL_n$, and~$\Phi$ is the canonical float-counting twisting.  (We will consider the rank-based twisting in Section \ref{sect:Psi}.)
We have already observed that tightness of $\Phi$ for~$\P_n$ implies unitality for all of our diagram monoids, as follows from Lemma \ref{lem:TUGA} and Remark \ref{rem:TUGA}.  Thus, the general results of Sections \ref{sect:*} and \ref{sect:*reg} apply to all of our twisted diagram monoids.  

The monoids $\P_0=\PP_0=\PB_0=\M_0=\B_0=\TL_0$ and $\B_1=\TL_1$ are trivial; in these cases our twisted diagram monoid $M\times_\Phi^qS$ is isomorphic to $M$, so we generally exclude these trivial cases from our considerations.  The remaining `edge cases' of 
\[
\P_1=\PP_1=\PB_1=\M_1=\left\{ \custpartn11{\stline11}, \custpartn11{}\right\} 
\COMMa 
\TL_2 = \left\{ \custpartn{1,2}{1,2}{\stline11\stline22}, \custpartn{1,2}{1,2}{\uarc12\darc12}\right\} 
\ANd
\B_2 = \left\{ \custpartn{1,2}{1,2}{\stline11\stline22}, \custpartn{1,2}{1,2}{\stline12\stline21}, \custpartn{1,2}{1,2}{\uarc12\darc12}\right\}
\]
are all commutative, which makes them somewhat different to the `generic cases' of $\P_n$, $\PP_n$, $\PB_n$, $\M_n$ ($n\geq2$), or $\B_n$,~$\TL_n$  ($n\geq3$).
The monoids~$\P_1$ and~$\TL_2$ are both isomorphic to the two element multiplicative semilattice~$\{1,0\}$, while $\B_2$ is the symmetric group $\S_2$ with a zero element attached.  We will see that $\B_2$ is unique among our diagram monoids, in that involutions of its twisted products need not be pure in general.

In the generic case, however, our diagram monoid $S$ has trivial centre (Theorem \ref{thm:Z}), in which case every involution of $T = M\times_\Phi^qS$ is pure (Corollary \ref{cor:inv}).  For $T$ to be star-regular with respect to such an involution, we need $S$ to be star-regular with respect to its underlying involution (cf.~Theorems \ref{thm:Pn*}--\ref{thm:PPn*}).  Beyond this, we need a star-regular involution of~$M$, a suitable map $a\mt k_a$, and for $\Phi(a,a^+,a)q$ to be a unit for all $a\in S$, where ${}^+$ is the Moore--Penrose inverse with respect to the given involution of $S$ (see Theorem \ref{thm:*reg}).
With respect to this last point, we record the following.

\begin{lemma}\label{lem:aa+a}
Let $S$ be any of $\P_n$, $\PP_n$, $\PB_n$, $\M_n\ (n\geq1)$, or $\B_n$, $\TL_n\ (n\geq2)$, and suppose~$S$ is $\tau$-regular with respect to an involution $\tau$ (see Theorems~\mbox{\ref{thm:Pn*}--\ref{thm:PPn*}}).  
\ben
\item  \label{aa+a1} If $S=\TL_3$ and $\tau = \rho$, then $\Phi(a,a^+,a) = 0$ for all $a\in \TL_3$.
\item \label{aa+a2} If $S=\B_n$, $\TL_n$ or $\P_1=\PP_1=\PB_1=\M_1$, and if $\tau = \io$, then $\Phi(a,a^+,a) \equiv 0 \pmod2$ for all $a\in S$, and we have $\Phi(a,a^+,a)=2$ for some $a\in S$.
\item \label{aa+a3} If $S=\PB_n$ or $\M_n$ for $n\geq2$, and if $\tau = \io$, then $\Phi(a,a^+,a)=3$ for some $a\in S$.
\item \label{aa+a4} Otherwise, we have $\Phi(a,a^+,a)=1$ for some $a\in S$.
\een
\end{lemma}

\pf
\firstpfitem{\ref{aa+a1}}  Since $\TL_3$ is a regular $\rho$-monoid by Theorem \ref{thm:TLn*}, we have $a^+=a\rho$ here.  The claim is now easily checked (using $\Phi(a,a^+,a)=\Phi(a,a^+)+\Phi(a^+,a)$; cf.~Lemma \ref{lem:aba}\ref{aba3}), as $\rho$ fixes each of $\custpartn{1,2,3}{1,2,3}{\stline11\stline22\stline33}  $, $  \custpartn{1,2,3}{1,2,3}{\stline13\uarc23\darc12} $ and $ \custpartn{1,2,3}{1,2,3}{\stline31\uarc12\darc23}$, and swaps $\custpartn{1,2,3}{1,2,3}{\stline11\uarc23\darc23} $ and $ \custpartn{1,2,3}{1,2,3}{\stline33\uarc12\darc12}$.

\pfitem{\ref{aa+a2}}  Here we have $a^+=a^*$, so by Lemma \ref{lem:Phiaa*a} we have $\Phi(a,a^+,a) = \Phi(a,a^*,a) = \non(a)$.  The first assertion then follows quickly from the fact that any Brauer partition has the same number of upper and lower blocks, as do the two elements of $\P_1$.  For the second, we can take $a=e_1= \custpartn{1,2,3,5}{1,2,3,5}{\stline33\stline55\uarc12\darc12\udotted35\ldotted35}\in\TL_n\sub\B_n$ or $a=t_1=\custpartn{1}{1}{}\in\P_1$.

\pfitem{\ref{aa+a3}}  Here we again have $\Phi(a,a^+,a) = \non(a)$, and we can take $a = \custpartn{1,2,3,5}{1,2,3,5}{\stline33\stline55\uarc12\udotted35\ldotted35}\in S$, which has three non-transversals.

\pfitem{\ref{aa+a4}}  If $S=\P_n$ or $\PP_n$ ($n\geq2$) and $\tau = \io$, then we still have $\Phi(a,a^+,a) = \non(a)$, and this time we can take $a = \custpartn{1,2,3,5}{1,2,3,5}{\stline33\stline55\uline12\stline12\stline22\udotted35\ldotted35}\in S$.

The only remaining case is when $S=\B_n$, and $\tau = \io_g$ for some transposition $g$.  By symmetry we can assume that $g=(1,2)$, and we take $a = \custpartn{0,1,2,3,5}{0,1,2,3,5}{\stline21\uarc01\darc02\stline33\stline55\udotted35\ldotted35}\in\B_n$.  Now,
\[
a^+ \H a\tau = ga^*g = ga^* = \custpartn{0,1,2,3,5}{0,1,2,3,5}{\stline02\uarc12\darc01\stline33\stline55\udotted35\ldotted35},
\]
so it follows from \ref{G1} and \ref{G2}, together with Lemma \ref{lem:aba}\ref{aba3}, that
\[
\Phi(a,a^+,a) = \Phi(a,a^+) + \Phi(a^+,a) = \Phi(a,ga^*) + \Phi(ga^*,a) = 0+1 = 1.  \qedhere
\]
\epf

Consequently, with the exception of $(S,\tau) = (\TL_3,\rho)$, the condition that $\Phi(a,a^+,a)q$ is a unit for all $a\in S$ is equivalent to $q$ itself being a unit.  Thus, excluding this special case, and looking at Theorem \ref{thm:*reg}, the \emph{pure} star-regular involutions of $T=M\times_\Phi^qS$ require $M$ to be $\circ$-regular,~$q$ to be a unit of~$M$, and have the form $(i,a)^\oa = (i^\circ+k_a,a\tau)$, where $S\to M:a\mt k_a$ satisfies the conditions of Theorem~\ref{thm:inv1}, with each $k_a$ a unit.  To say any more than this would require knowledge of the monoid $M$ itself.  
Our remaining considerations therefore take two separate directions:
\bit
\item In Section \ref{subsect:can} we look at the `canonical' involution of $M\times_\Phi^qS$ given by $(i,a)^\oa = (i,a\tau)$, where~$\tau$ is a (star-regular) involution of the underlying diagram monoid $S$.
\item In Section~\ref{subsect:Z} we look at the twisted monoids $\Z\times_\Phi^1S$ over the integers.
\eit

\subsection{The canonical involution}\label{subsect:can}

Consider again a non-trivial diagram monoid $S$, being any of $\P_n$, $\PP_n$, $\PB_n$, $\M_n$ ($n\geq1$), or $\B_n$, $\TL_n$~($n\geq2$).  By Corollary \ref{cor:*symD}, the canonical float-counting twisting $\Phi$ is $\tau$-symmetric with respect to every involution $\tau$ of $S$.  It then follows from Corollary \ref{cor:*sym} that any such involution $\tau$ extends to an involution of $T=M\times_\Phi^qS$, given by $(i,a)^\oa=(i,a\tau)$.  The next result characterises the cases in which $T$ is a $\oa$-regular monoid, or a regular $\oa$-monoid, with respect to this involution.

\newpage

\begin{thm}\label{thm:ia*}
Let $S$ be any of $\P_n$, $\PP_n$, $\PB_n$, $\M_n\ (n\geq1)$, or $\B_n$, $\TL_n\ (n\geq2)$, suppose $S$ is $\tau$-regular with respect to an involution $\tau$ (see Theorems~\mbox{\ref{thm:Pn*}--\ref{thm:PPn*}}), and let ${}^\oa$ be the involution of $T=M\times_\Phi^qS$ given by
\[
(i,a)^\oa = (i,a\tau) \qquad\text{for all $i\in M$ and $a\in S$.}
\]
\ben
\item \label{ia*1} If $S=\TL_3$ and $\tau=\rho$, then $T$ is $\oa$-regular if and only if $M$ is regular.
\item \label{ia*2} Otherwise, $T$ is $\oa$-regular if and only if $M$ is regular and $q$ is a unit.
\item \label{ia*3} If $S=\TL_3$ and $\tau=\rho$, then $T$ is a regular $\oa$-monoid if and only if $M$ is regular and has exponent at most $2$.
\item \label{ia*4} If $S=\B_n$, $\TL_n$ or $\P_1=\PP_1=\PB_1=\M_1$, and if $\tau=\io$, then $T$ is a regular $\oa$-monoid if and only if $M$ is regular and has exponent at most $2$, and $q$ is a unit.
\item \label{ia*5} If $S=\P_n$, $\PP_n$, $\PB_n$ or $\M_n\ (n\geq2)$, and if $\tau=\io$, then $T$ is a regular $\oa$-monoid if and only if $M$ is regular and has exponent at most $2$, and $q=0$.
\item \label{ia*6} Otherwise, $T$ is not a regular $\oa$-monoid.
\een
\end{thm}

\pf
We have already mentioned that $\Phi$ is $\tau$-symmetric with respect to each involution~$\tau$ of~$S$.

\pfitem{\ref{ia*1} and \ref{ia*2}}  These now follow from Corollary \ref{cor:*sym*reg} and Lemma \ref{lem:aa+a}.  The latter is used to interpret the condition `$\Phi(a,a^+,a)q$ is a unit for all $a\in S$' in Corollary \ref{cor:*sym*reg}; if $\Phi(a,a^+,a)\not=0$ for some $a\in S$, then this condition is equivalent to $q$ being a unit.

\pfitem{\ref{ia*3}--\ref{ia*6}}  These follow from Proposition \ref{prop:*symreg*} and Lemma \ref{lem:aa+a}, combined with the relevant parts of Theorems \ref{thm:Pn*}--\ref{thm:PPn*}.  (When considering $S=\PB_n$ or $\M_n$, note that $3q=q$ when $M$ has exponent $\leq2$.)
\epf

\subsection{Twisted diagram monoids over the integers}\label{subsect:Z}

We now move on to study the monoids $\Z\times_\Phi^1S$, where $S$ is one of $\P_n$, $\PP_n$, $\PB_n$, $\M_n$, $\B_n$ or~$\TL_n$.  As we have already mentioned, these are the most-studied \emph{regular} twisted diagram monoids in the literature.  It turns out that we can explicitly classify the involutions of these (Theorem \ref{thm:invT}), and characterise their star-regularity or regular-starity (Theorem \ref{thm:*regZD}).

We first show that in the generic case, involutions of twisted diagram monoids over $\Z$ are all pure.

\begin{lemma}\label{lem:Dext}
If $S$ is any of $\P_n$, $\PP_n$, $\PB_n$, $\M_n$ or $\TL_n$ \emph{(}for any $n$\emph{)}, or $\B_n$ \emph{(}for $n\not=2$\emph{)}, then every involution of $\Z\times_\Phi^1S$ is pure.
\end{lemma}

\pf
The result follows from Corollary \ref{cor:inv} if $S$ has trivial centre.  Given Theorem \ref{thm:Z} (and the exclusion of $\B_2$ in the statement), it remains only to consider $S=\P_1=\PP_1=\PB_1=\M_1$ and $S=\TL_2$.  In both cases we have $S = \{1,z\}$ for some non-unit $z$.  Assuming the involution has the form in Theorem \ref{thm:inv2}, it suffices by Proposition \ref{prop:ext} to show that $c_i=1$ for all $i\in \Z$.  But, as explained in Remark \ref{rem:unit}, $M=\Z$ being a group means that each $c_i$ is a unit of $S$, and the only unit of $S=\{1,z\}$ is $1$.
\epf

The monoid $\B_2$ was excluded from Lemma \ref{lem:Dext} because of the existence of the non-pure involutions in Examples \ref{eg:nonext} and \ref{eg:nonext2}.  We will return to these shortly.

Now let $S$ be any of $\P_n$, $\PP_n$, $\PB_n$, $\M_n$, $\B_n$ or $\TL_n$, and let ${}^\oa$ be an arbitrary \emph{pure} involution of $T = \Z\times_\Phi^1S$, extending some involution $\tau$ of $S$.  (By Lemma \ref{lem:Dext}, this includes \emph{all} involutions of $T$ unless $S=\B_2$.)  So by Theorem \ref{thm:inv1}, ${}^\oa$ has the form
\[
(i,a)^\oa = (i^\circ+k_a,a\tau) \qquad\text{for each $i\in\Z$ and $a\in S$,}
\]
where $\Z\to\Z:i\mt i^\circ$ is an involution, and $S\to\Z:a\mt k_a$ is a map satisfying:
\bena
\item \label{Da} $k_a^\circ+k_{a\tau} = 0$ for all $a\in S$, and
\item \label{Db} $k_a+k_b + \Phi(a,b) = k_{ab} + \Phi(a,b)1^\circ$ for all $a,b\in S$.
\een
(In \ref{Db}, recall that $q=1$ here, and also that $\Phi(b\tau,a\tau)=\Phi(a,b)$ by Corollary \ref{cor:*symD}.)  There are exactly two involutions of $\Z$, the identity map ($i^\circ=i$) and the negation map ($i^\circ=-i$), in which case we say that ${}^\oa$ is \emph{positive} or \emph{negative}, respectively.  We address these two kinds of involutions in the next two lemmas, in which we keep the current notation, including the involution $\tau$ of $S$, and the twisted product $T=\Z\times_\Phi^1S$.

\begin{lemma}\label{lem:pos}
The only positive (pure) involution of $T$ extending~$\tau$ is given by $(i,a)^\oa = (i,a\tau)$.
\end{lemma}

\pf
For a positive involution we have $i^\circ=i$ for all $i\in\Z$, so items \ref{Da} and \ref{Db} above become:
\begin{enumerate}[label=(\alph*)$'$,leftmargin=10mm]
\item \label{Da1} $k_a+k_{a\tau} = 0$, and 
\item \label{Db1} $k_a+k_b + \Phi(a,b) = k_{ab} + \Phi(a,b)$. 
\een
Cancelling $\Phi(a,b)$ from \ref{Db1} gives $k_a+k_b=k_{ab}$, meaning that $S\to\Z:a\mt k_a$ is a morphism.  It is easy to show that the only morphism from a finite semigroup to $\Z$ is the zero map, so it follows that $k_a=0$ for all $a\in S$.  Thus, in this case, we have
\[
(i,a)^\oa = (i^\circ+k_a,a\tau) = (i,a\tau) \qquad\text{for all $i\in\Z$ and $a\in S$.}  \qedhere
\]
\epf

\begin{lemma}\label{lem:neg}\leavevmode
\ben
\item \label{neg1} If $S$ is any of $\P_n$, $\PP_n$, $\PB_n$, $\M_n\ (n\geq2)$, or $\B_n$, $\TL_n\ (n\geq3)$, then there is no negative (pure) involution of $T$ extending $\tau$.
\item \label{neg2} If $S$ is any of $\P_1=\PP_1=\PB_1=\M_1$, $\B_2$ or $\TL_2$, then the only negative (pure) involution of $T$ extending $\tau$ \emph{(}which in this case is the identity map $\tau=\id$\emph{)} is given by ${(i,a)^\oa = (-i-\non(a),a)}$.
\een
\end{lemma}

\pf
For a negative involution we have $i^\circ=-i$ for all $i\in\Z$, so items~\ref{Da} and~\ref{Db} above become:
\begin{enumerate}[label=(\alph*)$''$,leftmargin=10mm]
\item \label{Da2} $k_a=k_{a\tau}$, and 
\item \label{Db2} $k_a+k_b + 2\Phi(a,b) = k_{ab}$.
\een
In particular, it follows from \ref{Db2} that
\begin{equation}\label{eq:ke}
k_e = -2\Phi(e,e) \qquad\text{for any idempotent $e\in E(S)$.}
\end{equation}

\pfitem{\ref{neg1}}
Suppose first that $S=\P_n$, $\PP_n$, $\PB_n$ or $\M_n$ (for $n\geq2$), and consider the idempotents
\[
e = e_1 = \custpartn{1,2,3,6}{1,2,3,6}{\stline33\stline66\uarc12\darc12\udotted36\ldotted36} \COMMA
f = t_1 = \custpartn{1,2,3,6}{1,2,3,6}{\stline22\stline33\stline66\udotted36\ldotted36} \AND
ef = \custpartn{1,2,3,6}{1,2,3,6}{\stline33\stline66\uarc12\udotted36\ldotted36} .
\]
By \eqref{eq:ke}, we have $k_e=k_f=k_{ef}=-2$.  Together with $\Phi(e,f)=0$, it follows from \ref{Db2} that
\[
-2 = k_{ef} = k_e + k_f + 2\Phi(e,f) = -4.
\]
This contradiction shows that no such negative involution exists in this case.

Now suppose $S=\B_n$ or $\TL_n$ (for $n\geq3$), and consider the idempotents
\[
e = e_1 = \custpartn{0,1,2,3,6}{0,1,2,3,6}{\uarc01\darc01\stline22\stline33\stline66\udotted36\ldotted36} \COMMA
f = e_2 = \custpartn{0,1,2,3,6}{0,1,2,3,6}{\stline00\uarc12\darc12\stline33\stline66\udotted36\ldotted36} \AND
ef = \custpartn{0,1,2,3,6}{0,1,2,3,6}{\uarc01\darc12\stline20\stline33\stline66\udotted36\ldotted36} .
\]
This time we have $k_e=k_f=-2$, $k_{ef}=0$ and $\Phi(e,f)=0$, which leads to
\[
0 = k_{ef} = k_e + k_f + 2\Phi(e,f) = -4,
\]
another contradiction.

\pfitem{\ref{neg2}}  If $S$ is one of $\P_1$ or $\TL_2$, then every element $a\in S$ satisfies $a^2=a=a^*$, so by \eqref{eq:ke} and Lemmas \ref{lem:aba}\ref{aba3} and \ref{lem:Phiaa*a} we have
\[
k_a = -2\Phi(a,a) = -(\Phi(a,a^*)+\Phi(a^*,a)) = -\Phi(a,a^*,a) = -\non(a).
\]
It is routine to verify that this (uniquely determined) mapping $a\mt k_a=-\non(a)$ satisfies conditions \ref{Da2} and \ref{Db2} above.  It leads to the stated involution of $T$, as
\[
(i,a)^\oa = (i^\circ+k_a,a\tau) = (-i-\non(a),a) \qquad\text{for all $i\in\Z$ and $a\in S$.}
\]

Finally, suppose $S=\B_2=\{1,s,e\}$, where $s=s_1$ and $e=e_1$.  Again we have ${k_a=-\non(a)}$ for $a\in E(S)=\{1,e\}$.  Setting $a=b=s$ in \ref{Db2} leads to $k_s=0=-\non(s)$.  It is again routine to verify that this mapping satisfies conditions \ref{Da2} and \ref{Db2} above, and leads to the stated involution.
\epf

We are now in a position to give the following classification of the involutions of our twisted diagram monoids.

\begin{thm}\label{thm:invT}
Let $S$ be any of $\P_n$, $\PP_n$, $\PB_n$, $\M_n\ (n\geq1)$, or $\B_n$, $\TL_n\ (n\geq2)$, and let~${T = \Z\times_\Phi^1S}$.
\ben
\item \label{invT1} For any involution $\tau$ of $S$ (see Theorem \ref{thm:invD}), there is an involution of $T$ given by
\[
(i,a)^\oa = (i,a\tau).
\]
\item \label{invT2} If $S$ is any of $\P_1=\PP_1=\PB_1=\M_1$, $\B_2$ or $\TL_2$, then there is another involution of $T$ given by
\[
(i,a)^\oa = (-i-\non(a),a).
\]
\item \label{invT3} If $S = \B_2$, then with $s=s_1$, there are two further involutions of $T$ given by
\[
(i,a)^\oa=(i,s^ia)
\OR 
(i,a)^\oa = (-i-\non(a),s^ia).
\]
\een
Moreover, every involution of $T$ has one of the above forms.
\end{thm}

\pf
By Lemmas \ref{lem:pos} and \ref{lem:neg}, the pure involutions of $T$ are precisely those listed in parts \ref{invT1} and \ref{invT2}.  By Lemma \ref{lem:Dext}, this completes the proof when $S\not=\B_2$.

For the rest of the proof we assume that $S=\B_2=\{1,s,e\}$, where $s=s_1$ and $e=e_1$.  The maps in \ref{invT3} are involutions, as we saw in Examples \ref{eg:nonext} and \ref{eg:nonext2}.  It remains to show that any non-pure involution ${}^\oa$ of $T$ is one of these two.  By Theorem~\ref{thm:inv2}, Proposition \ref{prop:unit}\ref{unit7} and Remark \ref{rem:unit} (remembering that $M=\Z$ is a group), we have
\[
(i,a)^\oa = (i^\circ+k_a,a^*c_i) \qquad\text{for all $i\in\Z$ and $a\in\B_2$,}
\]
for maps $i\mt i^\circ$, $i\mt c_i$, $a\mt a^*$ and $a\mt k_a$ satisfying \ref{inv1}--\ref{inv10}, and with each $c_i$ a unit.  If we had $c_1=1$, then \ref{inv1} would give $c_i=1$ for all $i\in\Z$; but this would contradict Proposition~\ref{prop:ext}, as ${}^\oa$ is not pure.  Since the only other unit of $\B_2$ is $s$, it follows that $c_1=s$, and \ref{inv1} now gives
\begin{equation}\label{eq:cixi}
c_i = s^i \qquad\text{for all $i\in\Z$.}
\end{equation}
Substituting this for the rule for ${}^\oa$, we have
\begin{equation}\label{eq:oa}
(i,a)^\oa = (i^\circ+k_a,s^ia^*) \qquad\text{for all $i\in\Z$ and $a\in\B_2$.}
\end{equation}
Next note that by \eqref{eq:0circ} we have
\[
k_1 = 0 \AND 1^* = 1.
\]
Since
\[
(0,1) = (0,1)^\oa = \big((0,s)(0,s)\big)^\oa = (0,s)^\oa(0,s)^\oa = (k_s,s^*)(k_s,s^*) = (2k_s+\Phi(s^*,s^*),s^*s^*),
\]
it follows that $s^*s^*=1$, so that $s^*=s$ or $1$.  In either case $s^*$ is a unit, so $\Phi(s^*,s^*)=0$, and continuing from above we have $(0,1) = (2k_s,s^*s^*)$, which gives
\[
k_s=0.
\]
Now, if we had $s^*=1$, then combining this with \ref{inv5} we would obtain ${s = s^{**}c_{k_s} = 1^*c_0 = 1}$, a contradiction.  Consequently,
\[
s^*=s.
\]
Condition \ref{inv5} again gives $e=c_{k_e}e^{**}$; since $e$ is not a unit, but $c_{k_e}$ is (by \eqref{eq:cixi}), it follows that~$e^{**}$ is not a unit, and hence $e^{**}=e$.  If $e^*$ was a unit, then $e=(e^*)^*$ would be a unit as well, as $1^*=1$ and $s^*=s$.  Thus, $e^*$ is not a unit, and so
\[
e^*=e.
\]
Since $\B_2=\{1,s,e\}$, we have therefore shown that 
\begin{equation}\label{eq:a*a}
a^* = a \qquad\text{for all $a\in\B_2$.}
\end{equation}
Next, for any $i\in\Z$ we have $c_i=s^i = s$ or $1$, so it follows that $k_{c_i} = 0$ (as $k_1=k_s=0$).  Condition~\ref{inv4} then gives $i^{\circ\circ} = i$ for all $i\in\Z$; combined with \ref{inv2} and $\Phi(s^i,s^j)=0$, it follows that $i\mt i^\circ$ is an involution of~$\Z$.  Thus, we either have
\bena
\item \label{icirca} $i^\circ=i$  for all $i\in\Z$, or
\item \label{icircb} $i^\circ=-i$ for all $i\in\Z$.
\een
Finally, substituting $a=b=e$ in \ref{inv8}, we obtain $2k_e+1=k_e+1^\circ$, so that
\[
k_e = 1^\circ-1 = \begin{cases}
0 &\text{in case \ref{icirca}}\\
-2=-\non(e) &\text{in case \ref{icircb}.}
\end{cases}
\]
Combined with $k_1=k_s=0$, and using \eqref{eq:cixi}--\eqref{eq:a*a}, it follows that ${}^\oa$ is the first involution stated in part \ref{invT3} in case \ref{icirca}, and the second involution in case \ref{icircb}.
\epf

Finally, we can characterise star-regularity and regular-starity of our twisted diagram monoids.

\newpage

\begin{thm}\label{thm:*regZD}
Let $S$ be any of $\P_n$, $\PP_n$, $\PB_n$, $\M_n\ (n\geq1)$, or $\B_n$, $\TL_n\ (n\geq2)$, let $T = \Z\times_\Phi^1S$, and let ${}^\oa$ be an involution of $T$.  Then $T$ is $\oa$-regular if and only if one of the following hold:
\ben
\item \label{*regZD1} $(i,a)^\oa = (i,a\tau)$, where $S$ is $\tau$-regular (see Theorems \ref{thm:Pn*}--\ref{thm:PPn*}),
\item \label{*regZD2} $S$ is one of $\P_1=\PP_1=\PB_1=\M_1$, $\B_2$ or $\TL_2$, and $(i,a)^\oa = (-i-\non(a),a)$,
\item \label{*regZD3} $S = \B_2$, and $(i,a)^\oa=(i,s^ia)$ or $(i,a)^\oa = (-i-\non(a),s^ia)$, where $s=s_1$.
\een
Moreover, $T$ is a regular $\oa$-monoid only in case \ref{*regZD2}.
\end{thm}

\pf
One simply goes through the involutions listed in Theorem \ref{thm:invT}, and checks them against the conditions of Theorems \ref{thm:*reggen}, \ref{thm:*reg} or \ref{thm:reg*} (or Corollary \ref{cor:*sym*reg} or Proposition \ref{prop:*symreg*}), as appropriate.
\epf

\section{The rank-based twisting}\label{sect:Psi}

We now follow the program of the previous section, this time focussing on twisted diagram monoids of the form $M\times_\Psi^qS$, where $S$ is one of $\P_n$, $\PP_n$, $\PB_n$, $\M_n$, $\B_n$ or $\TL_n$, and where $\Psi$ is the rank-based twisting introduced in \cite{EGPR2026}.  This twisting is denoted here by $\Psi$ to differentiate it from the float-counting twisting $\Phi$ considered earlier, and is defined by
\[
\Psi(a,b) = n - \rank(a) - \rank(b) + \rank(ab) \qquad\text{for $a,b\in\P_n$.}
\]
It is an example of a so-called \emph{rigid} twisting, as defined in \cite[Section 4.2]{EGPR2026}.  It is easy to show that \ref{T1} holds for this $\Psi$, but quite non-trivial to show that we always have $\Psi(a,b)\geq0$; this is the substance of \cite[Theorem~4.17]{EGPR2026}.  

As always, $\Psi$ restricts to a twisting of each of the submonoids $\PP_n$, $\PB_n$, $\M_n$, $\B_n$ and~$\TL_n$, and we denote these restrictions by $\Psi$ as well.  It was shown in \cite[Proposition 4.18]{EGPR2026} that $\Psi$ is tight for $\B_n$, but loose for all of the other monoids (except for trivially small $n$).  However, we always have unitality and Green-invariance:

\begin{lemma}\label{lem:Psi1}
If $S$ is any of $\P_n$, $\PP_n$, $\PB_n$, $\M_n$, $\B_n$ or $\TL_n$, then $\Psi$ is unital and Green-invariant.
\end{lemma}

\pf
Given Remark \ref{rem:TUGA}, we can assume that $S=\P_n$.

For unitality, let $a,b\in\P_n$, and suppose $a\in\S_n$ is a unit, so that $\rank(a) = n$.  Since $a$ is a unit, we have $ab\L b$, which implies $ab\J b$, and so $\rank(ab)=\rank(b)$ by \eqref{eq:RJ}.  This gives
\[
\Psi(a,b) = [n - \rank(a)] + [\rank(b) - \rank(ab)] = 0+0 = 0.
\]
An analogous calculation gives $\Psi(b,a) = 0$.  

For Green-invariance, we just establish \ref{G1}, as \ref{G2} is symmetrical.  So let $a,b,c\in\P_n$, and suppose $a\L b$.  Since $\L$ is a right congruence we also have $ac\L bc$.  Since ${\L}\sub{\J}$ we have $a\J b$ and $ac\J bc$, so it follows from \eqref{eq:RJ} that $\rank(a)=\rank(b)$ and $\rank(ac)=\rank(bc)$.  But then
\[
\Psi(a,c) = n - \rank(a) - \rank(c) + \rank(ac) = n - \rank(b) - \rank(c) + \rank(bc) = \Psi(b,c).  \qedhere
\]
\epf

We also have the following analogue of Corollary \ref{cor:*symD}.

\begin{lemma}\label{lem:Psi2}
If $S$ is any of $\P_n$, $\PP_n$, $\PB_n$, $\M_n$, $\B_n$ or $\TL_n$, and if $\tau$ is any involution of~$S$, then $\Psi$ is $\tau$-symmetric.
\end{lemma}

\pf
First note that $\Psi$ is $*$-symmetric for the canonical involution ${}^*$; this was observed in \cite[Theorem 4.17]{EGPR2026}, and follows from the obvious identity $\rank(a^*)=\rank(a)$.  Combined with unitality, this gives the result, exactly as in the proof of Corollary \ref{cor:*symD}.
\epf

For a star-regular involution $\tau$ of $S$, the values of $\Psi(a,a^+,a)$ are more constrained than those arising from the float-counting twisting $\Phi$ (cf.~Lemma \ref{lem:aa+a}):

\begin{lemma}\label{lem:Psi3}
Let $S$ be any of $\P_n$, $\PP_n$, $\PB_n$, $\M_n\ (n\geq1)$, or $\B_n$, $\TL_n\ (n\geq2)$, and suppose~$S$ is $\tau$-regular with respect to an involution $\tau$ (see Theorems~\mbox{\ref{thm:Pn*}--\ref{thm:PPn*}}).  
\ben
\item \label{Psi31} If $S=\B_n$ or $\TL_n$, then $\Psi(a,a^+,a) \equiv 0 \pmod4$ for all $a\in S$, and we have $\Psi(a,a^+,a)=4$ for some $a\in S$.
\item \label{Psi32} Otherwise, $\Psi(a,a^+,a) \equiv 0 \pmod2$ for all $a\in S$, and we have $\Psi(a,a^+,a)=2$ for some $a\in S$.
\een
\end{lemma}

\pf
Since $\Psi$ is Green-invariant by Lemma \ref{lem:Psi1}, we have $\Psi(a,a^+,a) = \Psi(a,a^+) + \Psi(a^+,a)$ by Lemma \ref{lem:aba}\ref{aba3}.  Since $a^+\L aa^+ \implies a^+\J aa^+ \implies \rank(a^+)=\rank(aa^+)$, we have
\[
\Psi(a,a^+) = n - \rank(a) - \rank(a^+) + \rank(aa^+) = n - \rank(a).
\]
A similar calculation gives $\Psi(a^+,a) = n - \rank(a)$, and so
\begin{equation}\label{eq:Psiaa+a}
\Psi(a,a^+,a) = 2n - 2\rank(a).
\end{equation}
This immediately gives $\Psi(a,a^+,a) \equiv 0\pmod2$.  If $a\in\B_n$, then since each non-transversal has size $2$, we must have $\rank(a) = n-2k$ for some $k$, and then $2n-2\rank(a) = 4k\equiv0\pmod4$.

Finally, note that $\Psi(a,a^+,a) = 2$ for any $a\in\M_n$ of rank $n-1$, while $\Psi(a,a^+,a) = 4$ for any $a\in\TL_n$ of rank $n-2$.
\epf

By Corollary \ref{cor:*sym} and Lemma \ref{lem:Psi2}, any involution $\tau$ of $S$ extends to a (pure) involution $(i,a)^\oa=(i,a\tau)$ of $T=M\times_\Psi^qS$.  We then have the following analogue of Theorem \ref{thm:ia*}, the proof of which is a simple application of Corollary \ref{cor:*sym*reg}, Proposition \ref{prop:*symreg*} and Lemma \ref{lem:Psi3}.

\begin{thm}\label{thm:Psi1}
Let $S$ be any of $\P_n$, $\PP_n$, $\PB_n$, $\M_n\ (n\geq1)$, or $\B_n$, $\TL_n\ (n\geq2)$, suppose $S$ is $\tau$-regular with respect to an involution $\tau$ (see Theorems~\mbox{\ref{thm:Pn*}--\ref{thm:PPn*}}), and let ${}^\oa$ be the involution of $T=M\times_\Psi^qS$ given by
\[
(i,a)^\oa = (i,a\tau) \qquad\text{for all $i\in M$ and $a\in S$.}
\]
\ben
\item \label{Psi11} $T$ is $\oa$-regular if and only if $M$ is regular and $q$ is a unit.
\item \label{Psi12} $T$ is a regular $\oa$-monoid if and only if $S$ is a regular $\tau$-monoid (see Theorems~\mbox{\ref{thm:Pn*}--\ref{thm:PPn*}}),~$M$ is regular and has exponent at most $2$, and $q$ is a unit.  \epfres
\een
\end{thm}

For the rest of the section we focus on the twisted diagram monoids $T = \Z\times_\Psi^1S$, so the analogues of Volkov's regular twisted monoids, but using the rank-based twisting $\Psi$ in place of the float-counting $\Phi$.  Some of the results and proofs of Section \ref{subsect:Z} are essentially identical for the current set-up; for these we simply reference the relevant result, and omit the proof.  But in other cases the situation is substantially different.  For example, while regular-starity holds only for three low-degree cases of float-based twisted monoids $\Z\times_\Phi^1S$, it holds for \emph{every} rank-based monoid $\Z\times_\Psi^1S$; compare Theorems \ref{thm:*regZD} and \ref{thm:*regZDPsi}.

\begin{lemma}[cf.~Lemma \ref{lem:Dext}]
If $S$ is any of $\P_n$, $\PP_n$, $\PB_n$, $\M_n$ or $\TL_n$ \emph{(}for any $n$\emph{)}, or $\B_n$ \emph{(}for $n\not=2$\emph{)}, then every involution of $\Z\times_\Psi^1S$ is pure.  \epfres
\end{lemma}

We will see below that there are two non-pure involutions of $\Z\times_\Psi^1\B_2$.  

Now let $S$ be any of $\P_n$, $\PP_n$, $\PB_n$, $\M_n$, $\B_n$ or $\TL_n$, and let ${}^\oa$ be a \emph{pure} involution of $T=\Z\times_\Psi^1S$, extending an involution $\tau$ of $S$.  So by Theorem \ref{thm:inv1}, ${}^\oa$ has the form
\[
(i,a)^\oa = (i^\circ+k_a,a\tau) \qquad\text{for each $i\in\Z$ and $a\in S$,}
\]
where $\Z\to\Z:i\mt i^\circ$ is an involution, and $S\to\Z:a\mt k_a$ is a map satisfying:
\bena
\item \label{DPsia} $k_a^\circ+k_{a\tau} = 0$ for all $a\in S$, and
\item \label{DPsib} $k_a+k_b + \Psi(a,b) = k_{ab} + \Psi(a,b)1^\circ$ for all $a,b\in S$.
\een
(Again recall that $q=1$, and that $\Psi(b\tau,a\tau)=\Psi(a,b)$ by Lemma \ref{lem:Psi2}.)  We again call ${}^\oa$ \emph{positive} or \emph{negative} when $i^\circ=i$ or $i^\circ=-i$, respectively.

\begin{lemma}[cf.~Lemma \ref{lem:pos}]\label{lem:posPsi}
The only positive (pure) involution of $T$ extending~$\tau$ is given by $(i,a)^\oa = (i,a\tau)$.  \epfres
\end{lemma}

In contrast to Lemma \ref{lem:neg}, $\tau$ \emph{always} extends to a negative involution of $T$.  The next proof makes use of the variant $\Z^{2n} = (\Z,\star)$ of $\Z$ with respect to the sandwich element $2n$, so that $\star$ is given by $i\star j = i+2n+j$ for $i,j\in\Z$.

\begin{lemma}\label{lem:negPsi}
The only negative (pure) involution of $T$ extending~$\tau$ is given by
\[
(i,a)^\oa = (-i+2\rank(a)-2n,a\tau).
\]
\end{lemma}

\pf
For a negative involution we have $i^\circ=-i$ for all $i\in\Z$, so items~\ref{Da} and~\ref{Db} above become:
\begin{enumerate}[label=(\alph*)$'$,leftmargin=10mm]
\item \label{Da2Psi} $k_a=k_{a\tau}$, and 
\item \label{Db2Psi} $k_a+k_b + 2\Psi(a,b) = k_{ab}$.
\een
Expanding $\Psi(a,b)$ in \ref{Db2Psi} and rearranging, we obtain
\[
[k_a - 2\rank(a)] + 2n + [k_b - 2\rank(b)] = k_{ab} - 2\rank(ab) \qquad\text{for all $a,b\in S$.}
\]
This says that the map $a\mt k_a - 2\rank(a)$ is a morphism $S\to\Z^{2n}$.  As we noted at the end of Section \ref{subsect:M}, we have~${\Z^{2n}\cong\Z}$.  As in the proof of Lemma \ref{lem:pos}, there is only one morphism from a finite semigroup to~$\Z$, and hence only one morphism $S\to\Z^{2n}$, namely the map sending every element of $S$ to the identity of $\Z^{2n}$, which is $-2n$.  Thus, we have $k_a - 2\rank(a) = -2n$, and so
\[
k_a = 2\rank(a) - 2n \qquad\text{for all $a\in S$.}
\]
Note that \ref{Da2Psi} holds for this $k_a$, by Corollary \ref{cor:rankatau}.  We then have
\[
(i,a)^\oa = (i^\circ+k_a,a\tau) = (-i + 2\rank(a) - 2n, a\tau) \qquad\text{for all $i\in M$ and $a\in S$.}  \qedhere
\]
\epf

Here is the classification of the involutions of $T=\Z\times_\Psi^1S$.  In part \ref{invTPsi2}, note that the `$-4$' term is `$-2n$' for $n=2$.

\newpage

\begin{thm}\label{thm:invTPsi}
Let $S$ be any of $\P_n$, $\PP_n$, $\PB_n$, $\M_n\ (n\geq1)$, or $\B_n$, $\TL_n\ (n\geq2)$, and let~${T = \Z\times_\Psi^1S}$.
\ben
\item \label{invTPsi1} For any involution $\tau$ of $S$ (see Theorem \ref{thm:invD}), there are two involutions of $T$ given by
\[
(i,a)^\oa = (i,a\tau)
\OR
(i,a)^\oa = (-i+2\rank(a)-2n,a\tau).
\]
\item \label{invTPsi2} If $S = \B_2$, then with $s=s_1$, there are two further involutions of $T$ given by
\[
(i,a)^\oa=(i,s^ia)
\OR 
(i,a)^\oa = (-i+2\rank(a)-4,s^ia).
\]
\een
Moreover, every involution of $T$ has one of the above forms.
\end{thm}

\pf
As in the proof of Theorem \ref{thm:invT}, we are left to deal with part \ref{invTPsi2}, for which we must show that the non-pure involutions of $T = \Z\times_\Psi^1\B_2$ are precisely the two stated maps.  We again write $\B_2=\{1,s,e\}$, where $s=s_1$ and $e=e_1$.
First, one can check using Theorem \ref{thm:inv2} that the stated maps are involutions, for which we respectively have
\bit
\item $a^*=a$, $i^\circ=i$, $k_a=0$ and $c_i=s^i$, or
\item $a^*=a$, $i^\circ=-i$, $k_a=2\rank(a)-4$ and $c_i=s^i$.
\eit
(Some of items \ref{inv1}--\ref{inv10} are trivial, in light of the fact that ${ab=ba \implies \Psi(a,b)=\Psi(b,a)}$.)  

Conversely, suppose ${}^\oa$ is a fixed non-pure involution of $T$.  Following the proof of Theorem~\ref{thm:invT}, we have $a^*=a$ for all $a\in\B_2$, $c_i=s^i$ for all $i\in \Z$, and $k_1=k_s=0$, so that ${}^\oa$ is given by the rule
\[
(i,a) = (i^\circ+k_a,s^ia).
\]
The map ${}^\circ$ is also an involution of $\Z$, so that
\bena
\item \label{icircaPsi} $i^\circ=i$  for all $i\in\Z$, or
\item \label{icircbPsi} $i^\circ=-i$ for all $i\in\Z$.
\een
Again substituting $a=b=e$ in \ref{inv8}, but this time using $\Psi(e,e)=2$, we obtain $2k_e+2=k_e+2^\circ$, so that
\[
k_e = 2^\circ-2 = \begin{cases}
0 &\text{in case \ref{icirca}}\\
-4 = 2\rank(e)-4 &\text{in case \ref{icircb}.}
\end{cases}
\]
Combined with $k_1=k_s=0$, we see that ${}^\oa$ is the first stated involution in case \ref{icircaPsi}, and the second in case \ref{icircbPsi}.
\epf

Finally, we can also characterise star-regularity and regular-starity of $T=\Z\times_\Psi^1S$.  We refer to Theorems \ref{thm:Pn*}--\ref{thm:PPn*} for the corresponding properties of involutions $\tau$ of $S$ itself.

\begin{thm}\label{thm:*regZDPsi}
Let $S$ be any of $\P_n$, $\PP_n$, $\PB_n$, $\M_n\ (n\geq1)$, or $\B_n$, $\TL_n\ ({n\geq2})$, let ${T = \Z\times_\Psi^1S}$, and let ${}^\oa$ be an involution of $T$.  Then $T$ is $\oa$-regular if and only if one of the following hold:
\ben
\item \label{*regZDPsi1} $(i,a)^\oa = (i,a\tau)$, where $S$ is $\tau$-regular,
\item \label{*regZDPsi2} $(i,a)^\oa = (-i+2\rank(a)-2n,a\tau)$, where $S$ is $\tau$-regular,
\item \label{*regZDPsi3} $S = \B_2$, and $(i,a)^\oa=(i,s^ia)$ or $(i,a)^\oa = (-i+2\rank(a)-4,s^ia)$, where $s=s_1$.
\een
Moreover, $T$ is a regular $\oa$-monoid only in case \ref{*regZDPsi2}, when $S$ is a regular $\tau$-monoid.
\end{thm}

\pf
The proof is analogous to that of Theorem \ref{thm:*regZD}.  For the final statement, note that for the involution in part \ref{*regZDPsi2} we have
\[
k_a + \Phi(a,a\tau,a) = [2\rank(a)-2n] + [2n-2\rank(a)] = 0,
\]
using \eqref{eq:Psiaa+a}, and remembering that $a^+=a\tau$ when $S$ is a regular $\tau$-monoid
\epf

\footnotesize
\def\bibspacing{-1.1pt}
\bibliography{biblio}
\bibliographystyle{abbrv}

\end{document}